\documentclass[11pt, leqno]{article}
\usepackage[nocompress]{cite}

\usepackage[utf8]{inputenc}
\usepackage[english]{babel}

\usepackage{amsmath}
\usepackage{geometry}
\usepackage{color}
\usepackage{xcolor}
\usepackage{enumitem}
\usepackage[linktoc=all]{hyperref}

\usepackage{graphicx}
\usepackage{wrapfig}
\usepackage{caption}
\usepackage{float}
\usepackage{tikz}

\usepackage{mathtools}
\usepackage{amssymb}
\usepackage{wasysym}
\usepackage{textcomp}
\usepackage{stmaryrd}
\usepackage{esint}

\usepackage{amsthm}
\usepackage{thmtools}

\usepackage{xspace}

\allowdisplaybreaks[4]

\colorlet{myred}{red!50!black}
\colorlet{mylightblue}{blue!50!black}
\colorlet{mydarkblue}{blue!80!black}
\colorlet{mygreen}{green!50!black}

\hypersetup{%
colorlinks = true,
linkcolor = myred,
citecolor = mylightblue,
urlcolor = mydarkblue
}

\usetikzlibrary{matrix}
\usetikzlibrary{arrows,decorations.markings}
\usetikzlibrary{intersections}
\tikzstyle{fleche}=[->, >=latex]
\tikzset{
xa/.store in=\xa, xa/.default=0, xa=0,%
xb/.store in=\xb, xb/.default=0, xb=0,%
xc/.store in=\xc, xc/.default=0, xc=0,%
}

\newcommand*{\N}{\mathbb{N}}

\newcommand*{\R}{\mathbb{R}}
\newcommand*{\C}{\mathbb{C}}

\newcommand*{\deq}{\overset{\scriptscriptstyle{\mathrm{def}}}{=}}

\newcommand{\V}{{\widehat{V}}}
\newcommand{\W}{{\widehat{W}}}

\newcommand{\ri}{\mathrm{i}}
\newcommand{\re}{\mathrm{e}}
\newcommand{\de}{\mathrm{d}}

\newcommand{\U}{\widehat{U}}
\newcommand{\Zc}{\widehat{Z}_{\mathrm{c}}}
\newcommand{\bigo}{\mathcal{O}}
\newcommand{\tl}{\tilde{\lambda}}

\def\eps {\varepsilon}

\renewcommand{\Re}{\operatorname{Re}}

\newcommand{\F}{\mathcal{F}}

\newcommand{\hL}{\hat{\Lambda}}

\newcommand{\rc}{\widehat{\rho}_{\mathrm{c}}}
\newcommand{\sic}{\varsigma_{\mathrm{c}}}

\newcommand*{\Norm}[1]{\left\lVert#1\right\rVert}

\newcommand*{\absolute}[1]{\lvert#1\rvert}

\newcommand*{\Set}[1]{\left\lbrace #1 \right\rbrace}

\newtheorem{theorem}{Theorem}[section]

\newtheorem{proposition}[theorem]{Proposition}
\newtheorem{lemma}[theorem]{Lemma}

\theoremstyle{definition}

\theoremstyle{remark}
\newtheorem{remark}[theorem]{Remark}

\numberwithin{equation}{section}

\title{Long-time behavior of small solutions in the viscoelastic Klein-Gordon equation}

\author{Louis Gar\'enaux$^*$ and Bj\"orn de Rijk\thanks{Department of Mathematics, Karlsruhe Institute of Technology, Englerstra\ss e 2, 76131 Karlsruhe, Germany; \texttt{louis.garenaux@kit.edu}, \texttt{bjoern.rijk@kit.edu}}}

\date{\today}

\newcommand{\calN}{\mathcal{N}}

\newcommand{\calR}{\mathcal{R}}
\newcommand{\calE}{\mathcal{E}}

\newcommand{\Uc}{\widehat{U}_\mathrm{c}}
\newcommand{\Us}{\widehat{U}_\mathrm{s}}

\begin{document}

\maketitle

\begin{abstract}
We investigate the long-time behavior of solutions with small initial data to the viscoelastic Klein-Gordon equation with general smooth nonlinearity. Our analysis relies on the space-time resonances method to eliminate all nonresonant quadratic and cubic terms. We identify a sign condition for the remaining critical resonant term to be of absorption type, leading to global-in-time existence and diffusive decay of solutions with small initial data. Even when this condition fails, our analysis shows existence and diffusive decay of small solutions on exponentially long time intervals.
\smallskip

\smallskip

\end{abstract}

\section{Introduction}

We study the long-time behavior of solutions with small initial data to the viscoelastic Klein-Gordon equation on the real line:
\begin{align}
u_{tt} - c^2 \partial_x^2 u - \alpha \partial_x^2 u_t + f(u) = 0, 
\qquad
u(x,t) \in \R, \, x \in \R, \, t \geq 0,
\label{e:vKG}
\end{align}
where $c > 0$ is the wave speed, $\alpha > 0$ is the viscosity coefficient, and $f \colon \R \to \R$ is a smooth nonlinearity with $f(0) = 0$ and $f'(0) > 0$. This damped nonlinear wave equation models the dynamics of an extended one-dimensional viscoelastic medium, a so-called Kelvin-Voigt solid, which exhibits both elastic behavior (instantaneous response to deformation) and viscous damping (dissipation of energy). Here, $u(x,t)$ represents the displacement of the viscoelastic solid at position $x$ and time $t$. The linear term $-\alpha \partial_x^2 u_t$ accounts for the viscous damping\footnote{This contrasts with a first-order damping term $\gamma u_t$ with $\gamma > 0$ on the left-hand side of~\eqref{e:vKG}, which dampens the amplitude of solutions without smoothing out oscillations, see~\cite{AVI,Potier-Ferry} and Remark~\ref{rem:damping}.}, while the nonlinearity $f(u)$ governs the material's elastic response. Since $f'(0) > 0$, the material behaves as a stiff elastic solid near equilibrium, with $f(u)$ acting as a restoring force. In the absence of damping and the nonlinearity,~\eqref{e:vKG} reduces to the standard wave equation $u_{tt} - c^2 \partial_x^2 u =0$. For further background, we refer to~\cite{Kawashima2,Potier-Ferry} and references therein.

By rescaling time, space, the displacement $u$, and the viscosity coefficient $\alpha > 0$, we can arrange for $f'(0) = 1$ and $c = 1$. This simplifies~\eqref{e:vKG} to
\begin{align}
u_{tt} + u - \partial_x^2 u - \alpha \partial_x^2 u_t = N(u), \label{e:vKG2}
\end{align}
where $N \colon \R \to \R$ is a smooth nonlinearity with $N(0) = N'(0) = 0$.

In this paper, we study the impact of the nonlinearity on the long-time behavior of solutions with small initial data. More precisely, we derive a sign condition on $N(u)$ ensuring that solutions with small initial data to~\eqref{e:vKG} exist for all positive times, and decay diffusively. Notably, our result implies the nonlinear stability of the equilibrium state $u(x,t) = 0$ against $L^2$-localized perturbations. We view the present contribution as a first step towards the nonlinear stability analysis of more complex solutions such as time- or space-periodic waves. 

\begin{figure}[h]
\begin{center}
\includegraphics[scale=0.9]{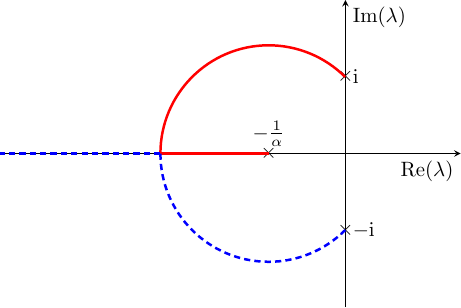}
\end{center}
\caption{Depiction of the spectrum of the operator $\Lambda$, which consists of the half line $\smash{(-\infty,-\tfrac{1}{\alpha}]}$ and the intersection of the closed left-half plane with the circle with center $\smash{-\tfrac{1}{\alpha}}$ and radius $\smash{\sqrt{1+\alpha^{-2}}}$. At frequency $k = 0$, the curves $\lambda_\pm(k)$ (depicted in red and blue) touch the imaginary axis in a quadratic tangency at the points $\pm \ri$. \label{fig:spectrum}}
\end{figure}

The linearized equation
\begin{align}
u_{tt} + u - \partial_x^2 u - \alpha \partial_x^2 u_t = 0 \label{e:vKG2_linear}
\end{align}
admits the solutions $u(x,t) = \re^{\pm \ri t}$. Consequently, expressing~\eqref{e:vKG2} as an evolution system in $(u,u_t)$, one finds that its linearization has continuous $L^2$-spectrum touching the imaginary axis at $\pm \ri$, see Figure~\ref{fig:spectrum}. Therefore, the associated semigroup exhibits at most algebraic decay rates, complicating the closure of a nonlinear iteration argument. In fact, a detailed analysis of the linear equation~\eqref{e:vKG2_linear} in~\cite{DAbbicco-Ikehata-19} shows that the optimal decay rate of its solutions is diffusive, which is typically insufficient to control quadratic or cubic nonlinearities. Indeed, all nonnegative nontrivial solutions to the nonlinear heat equation 
\begin{align} \label{e:nonlinear_heat}
u_t = \partial_x^2 u + u^p,
\end{align}
with $p = 2$ or $p = 3$, blow up in finite time~\cite{FUJI,HAYA}. 

Nevertheless, our setting differs fundamentally from that of the nonlinear heat equation, as the touchings with the imaginary axis occur at nonzero temporal frequencies. The associated time-oscillatory behavior of the critical modes carries over to the nonlinear terms in the variation-of-constants formula. Provided the nonlinear terms are not time-resonant, this results in oscillatory integrals with a nonstationary phase exhibiting enhanced temporal decay, which can be harnessed by integrating by parts in time, effectively following the space-time resonances method of Germain, Masmoudi, and Shatah~\cite{GERMAIN,SHAT10,GERM, GMS, GMS2}. 

Although originally developed for purely dispersive systems, the space-time resonances method has been successfully extended to dissipative settings. In~\cite{DRS1}, it is used to establish global existence and decay of solutions with small initial data in large classes of reaction-diffusion-advection systems where components exhibit different velocities, resulting in an absence of space resonances. The present work serves as a further extension of the method to dissipative systems, where one instead exploits an absence of time resonances.\footnote{We note that all nonlinear terms in~\eqref{e:vKG2} are space-resonant, since the two critical modes possess the same zero group velocity, see~\eqref{e:eigenvalue expansion}.}

As explained in the expository paper~\cite{GERMAIN}, the treatment of nonlinear terms that are not time-resonant is closely related to Shatah's normal form transform~\cite{Shatah-85}. The normal form transform has been adopted by Hayashi and Naumkin in~\cite{Hayashi-Naumkin-10,HANA,HANA3} to investigate the long-term behavior of solutions with small initial data to the Klein-Gordon equation~\eqref{e:vKG2} in the purely dispersive setting, where there is no damping ($\alpha = 0$). Using this method, global existence and decay of small solutions is obtained for the specific nonlinearities $N(u) = \nu u^2$ and $N(u) = \nu u^3$, where $\nu$ is real-valued. 

The current setting with viscous damping is significantly different from the one in~\cite{Hayashi-Naumkin-10,HANA,HANA3}. First, the damping term instantaneously regularizes solutions, so that we can afford to work in low regularity spaces. Second, at $\alpha = 0$ the linear equation~\eqref{e:vKG2} admits the solutions $u(x,t) = \re^{\ri (k x \pm \omega(k) t)}$ with $\omega(k) = \smash{\sqrt{1+k^2}}$ and $k \in \R$. As a result, the spectrum of the linearization occupies the imaginary axis, reflecting the lack of damping. For $\alpha > 0$ all frequencies but those at $k = 0$ are damped, which reduces the number of critical modes to two. Thus, after applying mode filters, only the occurrence of time resonances at Fourier frequency $0$ needs to be inspected. This enables us to handle \emph{general} smooth nonlinearities. In particular, we find that a nonlinear term $c_n u^n$ with $c_n \in \R$ and $n \in \N_{\geq 2}$ can only be time-resonant if $n$ is odd. Finally, both in~\cite{Hayashi-Naumkin-10,HANA,HANA3} and in our work, a critical resonant cubic term remains that cannot be handled with the space-time resonances method or normal form transform. At $\alpha = 0$, this resonant term has a purely imaginary coefficient after a change of variables, allowing it to be eliminated by an integrating factor with a purely imaginary phase. However, for $\alpha > 0$, the coefficient acquires a nonzero real part, whose sign determines whether it enhances decay (acting as an absorption mechanism) or obstructs the closure of a global nonlinear iteration argument. We refer to~\S\ref{sec:technical} for further details. 

Before presenting our main results, we review existing results on the global existence of solutions to the viscoelastic Klein-Gordon equation~\eqref{e:vKG2}. These results were obtained on \emph{bounded} domains\footnote{Notably, the linearization of~\eqref{e:vKG2} has discrete spectrum on bounded domains yielding a spectral gap, which is absent in our situation, see Figure~\ref{fig:spectrum}.} with the aid of energy estimates\footnote{The method of energy estimates is inherently different from our approach, as it accommodates solutions with large initial data, but is naturally more restrictive on the type of nonlinearities that can be handled.}. The first result~\cite{WE80} considers~\eqref{e:vKG2} for nonlinearities $N$ for which there exists a constant $C > 0$ such that $N'(u) \leq C$ for all $u \in \R$. It asserts that solutions are global and converge to a stationary solution as $t \to \infty$, see also~\cite{Dinh88} for explicit temporal decay rates. Existence of global solutions on bounded domains in the specific case of the power-law nonlinearity 
\begin{align*}
N(u) = -\nu |u|^{p-1} u 
\end{align*}
with $\nu > 0$ and integer $p \geq 1$ was established in~\cite{AVI}. These global solutions converge in the vanishing viscosity limit $\alpha \downarrow 0$ by the results in~\cite{AVR2,AVR}.

\subsection{Main results}

Our first result establishes global existence and diffusive decay of solutions to~\eqref{e:vKG2} with small initial data, allowing for general smooth nonlinearities $N$ obeying the sign condition
\begin{align}
3N'''(0) + 5 N''(0)^2 < 0. \label{e:signcondition}
\end{align} 
Specifically, we consider initial data $u_0 \in H^2(\R)$ whose Fourier transform 
\begin{align*}
\hat{u}_0(k) \deq \int_\R \re^{-\ri k x} u(x) \de x
\end{align*}
lies in the Sobolev space $W^{1,1}(\R) \cap W^{1,\infty}(\R)$ and is small in that space.

\begin{theorem}[Global existence and diffusive decay] \label{thm:mainresult}
Let $\alpha > 0$. Take $N \in C^4(\R)$ such that $N(0) = 0$, $N'(0) = 0$, and the inequality~\eqref{e:signcondition} holds. Then, there exist positive constants $M_0$ and $\eps$ such that, whenever $u_0 \in H^2(\R)$ and $w_0 \in L^2(\R)$ satisfy $\hat{u}_0,\hat{w}_0 \in W^{1,\infty}(\R) \cap W^{1,1}(\R)$ and
\begin{align} \label{e:cond_initial}
E_0 \deq \|\hat{u}_0\|_{W^{1,1} \cap W^{1,\infty}} + \|\hat{w}_0\|_{W^{1,1} \cap W^{1,\infty}} < \eps,
\end{align}
there exists a unique global classical solution
\begin{align*}
u &\in C\big([0,\infty),H^2(\R)\big) \cap C^1\big([0,\infty),L^2(\R)\big) \cap C^1\big((0,\infty),H^2(\R)\big) \cap C^2\big((0,\infty),L^2(\R)\big)
\end{align*}
of the viscoelastic Klein-Gordon equation~\eqref{e:vKG2} with initial conditions $u(0) = u_0$ and $u_t(0) = w_0$, which enjoys the diffusive estimates
\begin{align} \label{e:diffest}
\left\|u(t)\right\|_{L^\infty} \leq \frac{M_0E_0}{\sqrt{1+t}}, \qquad \left\|u(t)\right\|_{L^2} \leq \frac{M_0E_0}{(1+t)^{\frac14}}
\end{align}
and the enhanced pointwise estimate
\begin{align} 
\label{e:diffest2}
\left|u(x,t)\right| \leq \frac{M_0}{\sqrt{(1+t)\log(2+t)}} \exp\left(-\frac{\alpha x^2}{2(1+\alpha^2)(1+t)}\right) + \frac{M_0}{\sqrt{1+t} \log(2+t)^{\frac23}}
\end{align}
for all $x \in \R$ and $t \geq 0$. 
\end{theorem}

Theorem~\ref{thm:mainresult} establishes nonlinear asymptotic stability of the equilibrium state $u(x,t) = 0$ in~\eqref{e:vKG2} against $H^2$-perturbations whose Fourier transform is small in $W^{1,1}(\R) \cap W^{1,\infty}(\R)$. We note that such perturbations must be algebraically localized, see~\S\ref{sec:notation} for further details.

The estimate~\eqref{e:diffest} shows that solutions to~\eqref{e:vKG} with small initial data decay at the same rates as solutions to the heat equation $u_t = \partial_x^2 u$. Identical diffusive decay rates are obtained in the purely dispersive setting ($\alpha = 0$) in~\cite{HANA,HANA3} for the nonlinearities $N(u) = \nu u^p$, with $\nu \in \R$ and $p = 2,3$. However, the Gaussian principal part of the pointwise bound~\eqref{e:diffest2}, reflecting the viscous nature of equation~\eqref{e:vKG2}, implies that the long-term asymptotics of the solution $u(x,t)$ must be fundamentally different from the purely dispersive case, cf.~\cite{HANA}. The proof of Theorem~\ref{thm:mainresult} reveals that, as long as the nonlinearity obeys the sign condition~\eqref{e:signcondition}, the critical resonant term is cubic and of absorption type, contributing an additional logarithmic decay factor $\smash{\sqrt{\log(2+t)}}$ in~\eqref{e:diffest2}. The same phenomenon has been observed in the cubic heat equation $u_t = \partial_x^2 u - u^3$ with absorption, cf.~\cite{GALAK-85}.

If the nonlinearity does \emph{not} satisfy the sign condition~\eqref{e:signcondition}, then the critical resonant term obstructs the closure of a global nonlinear iteration argument. However, since the quadratic terms are not time-resonant and can be eliminated, we are still able to establish existence and diffusive decay on time intervals which are exponentially large with respect to the size of the initial data.

\begin{theorem}[Existence and diffusive decay on exponentially long time scales] \label{thm:mainresult2}
Let $\alpha > 0$. Take $N \in C^4(\R)$ such that $N(0) = 0$ and $N'(0) = 0$. Then, there exist positive constants $M_0$ and $\eps$ such that, whenever $u_0 \in H^2(\R)$ and $w_0 \in L^2(\R)$ satisfy $\hat{u}_0, \hat{w}_0 \in L^1(\R) \cap L^\infty(\R)$ and
\begin{align} \label{e:cond_initial2}
E_0 \deq \|\hat{u}_0\|_{L^1 \cap L^\infty} + \|\hat{w}_0\|_{L^1 \cap L^\infty} < \eps,
\end{align}
there exists a unique classical solution
\begin{align*}
u &\in C\big([0,T_\eps],H^2(\R)\big) \cap C^1\big([0,T_\eps],L^2(\R)\big) \cap C^1\big((0,T_\eps],H^2(\R)\big) \cap C^2\big((0,T_\eps],L^2(\R)\big)
\end{align*}
of the viscoelastic Klein-Gordon equation~\eqref{e:vKG2} on an interval of length
\begin{align*}
T_\varepsilon \deq \re^{\eps/E_0}-2
\end{align*}
with initial data $u(0) = u_0$ and $u_t(0) = w_0$. Moreover, $u(t)$ obeys the diffusive estimates
\begin{align*}
\left\|u(t)\right\|_{L^\infty} \leq \frac{M_0E_0}{\sqrt{1+t}}, \qquad \left\|u(t)\right\|_{L^2} \leq \frac{M_0E_0}{(1+t)^{\frac14}}
\end{align*}
 for all $t \in [0,T_\eps]$.
\end{theorem}

A result similar to Theorem~\ref{thm:mainresult2} holds for the nonlinear heat equation~\eqref{e:nonlinear_heat} with cubic nonlinearity ($p=3$). Specifically, it is shown in~\cite[Theorem~3.21]{LENI92} and~\cite[Theorem~2.1]{SUP03} that solutions with initial data of size $E_0$ in $L^1(\R) \cap L^\infty(\R)$ exist and decay diffusively on time intervals that are exponentially long in $E_0$. In contrast, in case of a quadratic nonlinearity ($p=2$), solutions to~\eqref{e:nonlinear_heat} with small initial data of size $E_0$ in $L^1(\R) \cap L^\infty(\R)$ can blow up within time $\leq CE_0^{-2}$, where $C>0$ is some $E_0$-independent constant, cf.~\cite[Theorem~3.15]{LENI92}. 

The question of whether solutions exhibit blow-up when~\eqref{e:signcondition} is not satisfied remains open. In the viscous regime with $\alpha > 0$, we do not expect solutions to preserve compact support, which complicates the adaptation of known blow-up results~\cite{KETA,Delort-99}. Nevertheless, as discussed in Remark~\ref{rem:blowup}, we expect that the failure of~\eqref{e:signcondition} leads to the instability of the rest state $u(x,t) = 0$ in~\eqref{e:vKG2}.

\begin{remark}
Another approach to obtaining global existence of solutions with small initial data is to consider so-called \emph{transparent} nonlinearities~\cite{Lannes-12}, which vanish at resonant frequencies, thereby eliminating the singularity that typically arises when integrating by parts. For the Klein-Gordon equation, known transparent terms involve a spatial derivative~\cite{Moriyama-97,Katayama-99}. Our analysis shows that the critical quadratic and cubic terms in~\eqref{e:vKG2} are not transparent. However, in case $3N'''(0) + 5 N''(0)^2 = 0$ the coefficient in front of the critical remaining resonant term vanishes and we expect that Theorem~\ref{thm:mainresult2} can be extended to a global result. We leave this boundary case as an open question for future research.
\end{remark}

\subsection{Technical summary} \label{sec:technical}

In this section, we provide an outline of our analysis which eventually leads to the proofs of Theorems~\ref{thm:mainresult} and~\ref{thm:mainresult2}. The main idea is to eliminate nonresonant critical terms using the space-time resonances method and arrive at a reduced nonlinear ordinary differential equation governing the leading-order dynamics of solutions to~\eqref{e:vKG2} with small initial data. The reduction process involves several steps that we summarize below.

We start by recasting equation~\eqref{e:vKG2} as a first-order system in time. Using the change of variable  
\begin{equation*}
U = 
\begin{pmatrix}
u \\ \big(1-\partial_{x}^2\big)^{-1}\big(u_t - \frac{\alpha}{2} \partial_x^2 u\big)
\end{pmatrix}
\end{equation*}
whose components have balanced regularity, we arrive at a system of the form  
\begin{equation}
\label{e:heuristic-1}
U_t = \Lambda U + \calN_2(U) + \calN_3(U) + \calR(U).
\end{equation}
Here, the linear operator $\Lambda$ is sectorial, and has spectrum that touches the imaginary axis in a quadratic tangency at $\pm \ri$, see Figure~\ref{fig:spectrum}. This indicates that the semigroup $\re^{\Lambda t}$ obeys the same diffusive estimates as the heat semigroup: $\Norm{\smash{\re^{\partial_x^2 t}}}_{L^p \to L^\infty} \leq C t^{-1/(2p)}$, $1 \leq p \leq \infty$. The quadratic and cubic terms, $\calN_2(U)$ and $\calN_3(U)$, pose the main challenge in proving global existence, whereas the residual term $\calR(U)$ is at least quartic in $U$ and thus irrelevant for the long-time dynamics. Indeed, proving global existence of small solutions via the mild formulation  
\begin{equation*}
U(t) = \re^{\Lambda t} U_0 + \int_0^t \re^{\Lambda (t - \tau)} \big(\calN_2(U) + \calN_3(U) + \calR(U) \big)(\tau) \, \de \tau
\end{equation*}
of~\eqref{e:heuristic-1} through $L^1$-$L^\infty$-estimates is only possible if the cumulative nonlinear effects remain integrable over time, cf.~\cite[Section~14]{SUBook}. Given the aforementioned diffusive temporal decay rates, quadratic and cubic nonlinear terms are critical, while the decay rates associated with quartic or higher-order terms terms are integrable in time.

Our first step in analyzing the dynamics of \eqref{e:heuristic-1} is to decompose the solution into low- and high-frequency components in order to distinguish between critical oscillatory modes and exponentially damped ones. Writing $\widehat{U} = \Uc + \Us = \chi \widehat{U} + (1-\chi) \widehat{U}$, where $\chi$ is a cut-off function centered at the origin, we see that problematic nonlinear terms are low-frequency interactions only. That is, the dynamics in Fourier space is given by
\begin{equation*}
\begin{cases}
\partial_t \Uc = \widehat{\Lambda} \Uc + \chi \widehat{\calN}_2(\Uc) + \chi \widehat{\calN}_3(\Uc) + \calE_1,\\
\partial_t \Us = \widehat{\Lambda} \Us + \calE_2.
\end{cases}
\end{equation*}
where all terms denoted by $\calE$, here and in the following, represent irrelevant nonlinear terms that do not affect the long-time behavior of solutions with small initial data. Next, we introduce a near-identity change of variables
\begin{equation} \label{e:transformation}
V = \Uc + K_2(\Uc) + K_3(\Uc),
\end{equation}
for the critical low-frequency component, where $K_2(\Uc)$ is quadratic in $\Uc$ and $K_3(\Uc)$ is cubic in $\Uc$. This normal form transformation is invertible and designed to remove all quadratic terms and all non-resonant cubic terms, leading to the evolution equation
\begin{equation*}
V_t = \widehat{\Lambda} V + Q_{\mathrm{res}}(V) + \calE_3,
\end{equation*}
where $Q_{\mathrm{res}}(V)$ is a resonant cubic term. The transformation~\eqref{e:transformation} arises from the space-time resonances method by integrating by parts in time. At this stage, we can already establish existence and diffusive decay of small solutions on exponentially long time scales, as in~\cite{Moriyama-Tonegawa-Tsutsumi-97}, leading to the proof of Theorem~\ref{thm:mainresult2}. 

The proof of Theorem~\ref{thm:mainresult} requires a final reduction step in which we identify the leading-order behavior of $V$ as a diffusive Gaussian profile $\hat{g}$ with a complex-valued amplitude $A$. Proceeding as in~\cite[Theorem~1.5]{DRS}, we write
\begin{equation*} V(t,k) = A(t) \hat{g}(t, k) + \hat{\rho}(t,k), \end{equation*}
where $\hat{\rho}$ denotes a residual term exhibiting higher-order decay. We arrive at an evolution system of the form
\begin{equation*} \begin{cases} \partial_t \hat{\rho} = \widehat{\Lambda} \hat{\rho} + \calE_4, \\ r'(t) = \frac{\omega}{1+t} r(t)^3 + \calE_5. \end{cases} \end{equation*}
for $r = \absolute{A}$ and $\hat{\rho}$. The sign of the parameter $\omega$ governs the long-time behavior of the ODE for $r(t)$. In particular, if the sign condition~\eqref{e:signcondition} holds, then we find $\omega < 0$ and solutions exist globally in time. Exploiting this, we are able to close a global nonlinear iteration argument and prove Theorem~\ref{thm:mainresult}. 

\begin{remark} \label{rem:blowup}
If $\omega > 0$, then all solutions to the separable ODE $r'(t) = \tfrac{\omega}{1+t} r(t)^3$ with initial data $r(0) = r_0\neq 0$ blow up at $t_0 = \smash{\re^{1/(2 \omega r_0^2)}-1}$. This suggests, as in the nonlinear heat equation~\eqref{e:nonlinear_heat} with $p = 3$, cf.~\cite{HAYA,LENI92}, that the rest state $u(x,t) = 0$ in~\eqref{e:vKG2} is unstable if $3N'''(0) + 5 N''(0)^2 > 0$, with the instability only manifesting itself on exponentially long time scales. 
\end{remark}

\paragraph{Organization.} In~\S\ref{sec:notation} we introduce some notation and define the function spaces in which we consider the solutions to~\eqref{e:vKG2}. Section~\ref{sec:locwellp} is devoted to the local existence analysis of solutions to~\eqref{e:vKG2}. In~\S\ref{sec:linear} we study the linear dynamics of~\eqref{e:vKG2}. Subsequently, in~\S\ref{sec:modefilters} we apply mode filters to~\eqref{e:vKG2} which separate low-frequency from high-frequency modes. In~\S\ref{sec:isolate} we isolate the critical quadratic and cubic terms and estimate the irrelevant residual terms. We eliminate the quadratic and nonresonant cubic terms with the aid of the space-time resonances method in~\S\ref{sec:eliminate2} and~\ref{sec:eliminate3}, respectively. Section~\ref{sec:rescubic} contains estimates on the resonant cubic terms, whereas in~\S\ref{s:reduced-eq} we analyze the reduced system governing the leading-order dynamics. Finally, we close the nonlinear iteration argument in~\S\ref{s:theorem-proofs}, which finishes the proofs of Theorems~\ref{thm:mainresult} and~\ref{thm:mainresult2}. 

\paragraph*{Acknowledgment.} Funded by the Deutsche Forschungsgemeinschaft (DFG, German Research Foundation) -- Project-ID 258734477 -- SFB 1173.

\section{Notation and function spaces} \label{sec:notation}

In this section, we introduce some notation and define the function spaces in which we construct solutions to the viscoelastic Klein-Gordon equation~\eqref{e:vKG2}. 

First of all, given a set $S$ and maps $A, B \colon S \to \R$, we write ``$A(x) \lesssim B(x)$ for $x \in S$'' to express that there exists a constant $C>0$, independent of $x$, such that $A(x) \leq CB(x)$ holds for all $x \in S$. 

Second, we employ the nonunitary Fourier transform $\F \colon L^2(\R) \to L^2(\R)$ and its inverse $\F^{-1} \colon L^2(\R) \to L^2(\R)$ throughout this paper. They are determined by their action on the dense subspace $L^1(\R) \cap L^2(\R)$ of $L^2(\R)$, which is given by
\begin{align*}
\F(u)(k) = \int_\R \re^{-\ri k x} u(x) \de x, \qquad \F^{-1}(v)(k) = \frac{1}{2\pi} \int_\R \re^{\ri k x} v(k) \de k.
\end{align*}
As usual, we abbreviate $\hat{u} = \F(u)$. 

Next, we introduce the algebraically weighted $L^2$-space
\begin{align*}
L^2_1(\R) = \left\{f \in L^2(\R) : \rho f \in L^2(\R)\right\},
\end{align*}
where $\rho \colon \R \to \R$ is the weight $\rho(k) = \sqrt{1+k^2}$. We equip $L^2_1(\R)$ with the norm $\smash{\|f\|_{L^2_1}} = \smash{\|\rho f\|_{L^2}}$. It is well-known that the Fourier transform maps $L^2_1(\R)$ isomorphically onto $H^1(\R)$. Moreover, we define the Banach spaces 
\begin{align*}
Y_m &= \big\{f \in L^2(\R) : \rho^m \hat{f} \in L^1(\R) \cap L^\infty(\R)\big\},\\
X_m &= \big\{f \in L_1^2(\R) : \rho^m \hat{f}, \rho^m \hat{f}' \in L^1(\R) \cap L^\infty(\R)\big\},
\end{align*}
for $m \in \mathbb N_0$, by their norms 
\begin{align*}
\|f\|_{Y_m} &= \big\|\rho^m \hat{f}\big\|_{L^1} + \big\|\rho^m \hat{f}\big\|_{L^\infty},\\
\|f\|_{X_m} &= \big\|\rho^m \hat{f}\big\|_{L^1} + \big\|\rho^m \hat{f}\big\|_{L^\infty} + \big\|\rho^m \hat{f}'\big\|_{L^1} + \big\|\rho^m \hat{f}'\big\|_{L^\infty},
\end{align*}
respectively. For $m \in \mathbb N_0$, we have the continuous embeddings 
\begin{align*}
W^{m,1}(\R) \cap W^{m,\infty}(\R) \hookrightarrow Y_m \hookrightarrow H^m(\R) \cap W^{m,\infty}(\R),
\end{align*}
and, similarly,
\begin{align*}
&\left\{f \in L^2(\R) : f, \rho f \in W^{m,1}(\R) \cap W^{m,\infty}(\R)\right\}  \\
&\qquad \ \hookrightarrow X_m \hookrightarrow \left\{f \in L^2(\R) : f, \rho f \in H^m(\R) \cap W^{m,\infty}(\R)\right\}.
\end{align*}

\section{Local existence and uniqueness} \label{sec:locwellp}

We write the viscoelastic Klein-Gordon equation~\eqref{e:vKG2} as a semilinear evolution problem by introducing the tailor made variable
\begin{align}
\label{e:diff-variable}
v = \left(1-\partial_x^2\right)^{-1} \left(u_t - \frac{\alpha}{2} \partial_x^2 u\right).
\end{align}
Thus, we obtain the system
\begin{align} \label{e:vKGsys}
U_t = \Lambda U + \mathcal{N}(U),
\end{align}
in $U = (u,v)^\top$, where the linear operator $\Lambda$ is given by
\begin{align*}
\Lambda = \begin{pmatrix} \frac{\alpha}{2} \partial_x^2 & 1-\partial_x^2 \\ -1 + \frac{\alpha^2}{4} \partial_x^4 \left(1-\partial_x^2\right)^{-1} & \frac{\alpha}{2} \partial_x^2\end{pmatrix},
\end{align*}
and the nonlinearity $\mathcal{N}(U)$ is defined by
\begin{align*}
\mathcal{N}(U) = \left(1-\partial_x^2\right)^{-1} N\left(U_1\right) \mathbf{e}_2,
\end{align*}
where $\mathbf{e}_2$ is the unit vector $\mathbf{e}_2 = (0, 1)^\top$ and $U_1$ denotes the first coordinate of the vector $U$.

We will establish that the linear operator $\Lambda$ in~\eqref{e:vKGsys} is sectorial on the spaces $Y_0$ and $L^2(\R)$, implying that~\eqref{e:vKGsys} is a semilinear parabolic problem. We note that this is a direct consequence of the viscous dissipation, modeled by the term $-\alpha \partial_x^2 u_t$ in~\eqref{e:vKG2}. Thus, standard parabolic semigroup theory~\cite{LUN} yields local existence of a maximal mild solution to \eqref{e:vKGsys} in the spaces $X_0$ and $Y_0$, required for the proof of Theorems~\ref{thm:mainresult} and~\ref{thm:mainresult2}, respectively, as well as in the space $L^2(\R)$. In the subsequent result we then show, under the additional regularity assumption $u(0) \in H^2(\R)$, that, if $U(t) = (u(t),v(t))^\top$ is a mild solution in $L^2(\R)$ of system~\eqref{e:vKGsys}, then its first coordinate $u(t)$ is a classical solution of the viscoelastic Klein-Gordon equation~\eqref{e:vKG2}.

\begin{remark}
Upon introducing the traditional variable $z = u_t$, equation~\eqref{e:vKG2} can be written as the semilinear evolution problem
\begin{align} \label{e:vKGsysAlter}
Z_t = \begin{pmatrix} 0 & 1 \\ \partial_x^2 - 1 & \alpha \partial_x^2\end{pmatrix} Z + \begin{pmatrix} 0 \\ N(u)\end{pmatrix},
\end{align}
in $Z = (u,z)^\top$. It is well-known~\cite{ITFY13,WE80} that the linear operator 
\begin{align} \label{e:defA}
\mathcal{A} = \begin{pmatrix} 0 & 1 \\ \partial_x^2 - 1 & \alpha \partial_x^2\end{pmatrix},
\end{align}
generates a $C^0$-semigroup on the space $H^1(\R) \times L^2(\R)$ and an analytic semigroup on the space $H^2(\R) \times L^2(\R)$. These facts are used in the upcoming Proposition~\ref{p:wellposedness2}. The disadvantage of the first-order formulation~\eqref{e:vKGsysAlter} with respect to~\eqref{e:vKGsys} is that the components of the vector $Z = (u,z)^\top$ do not have the same regularity in the spaces $H^1(\R) \times L^2(\R)$ and $H^2(\R) \times L^2(\R)$, which is avoided by the preconditioner $(1-\partial_x^2)^{-1}$ in~\eqref{e:diff-variable}. Moreover, the preconditioner induces additional localization on the nonlinearity in~\eqref{e:vKGsys} in Fourier space. For these reasons we adopt the first-order formulation~\eqref{e:vKGsys} in our nonlinear analysis. 
\end{remark}

We prove local existence of a maximal mild solution to~\eqref{e:vKGsys} in one of the spaces $X_0, Y_0$ or $L^2(\R)$.

\begin{proposition} \label{p:wellposedness}
Let $\mathcal{X}$ be one of the spaces $X_0$, $Y_0$ or $L^2(\R)$. Let $U_0 \in \mathcal{X}$. Then, there exist $T_{\max} \in (0,\infty]$ and a unique, maximally defined, mild solution $U \in C\big([0,T_{\max}),\mathcal{X}\big)$ of~\eqref{e:vKGsys} with initial condition $U(0) = U_0$. If $T_{\max} < \infty$, then it holds 
\begin{align} \label{e:blowup}
\limsup_{t \uparrow T_{\max}} \|U(t)\|_{\mathcal{X}} = \infty.
\end{align}
\end{proposition}
\begin{proof}
We first consider the cases $\mathcal{X} = Y_0$ or $\mathcal{X} = L^2(\R)$. Let $Z_m$ denote either the space $Y_m$ or $H^m(\R)$ for $m \in \mathbb{N}_0$. We show that $\Lambda$ is a sectorial operator on $Z_0$ and $\mathcal{N}$ is a locally Lipschitz continuous nonlinearity on $Z_0$. Then, the existence of a maximal mild solution to~\eqref{e:vKGsys} follows from standard semigroup theory for semilinear parabolic problems.

We observe that the elliptic operator $\partial_x^2$ acts on $Z_0$ with dense domain $Z_2$. Thus, the preconditioner $(1-\partial_x^2)^{-1}$ is a bounded linear operator from $Z_2$ into $Z_0$. We aim to show that $\Lambda$ is a sectorial operator on $Z_0$ by regarding $\Lambda$ as a bounded perturbation of the operator $\Lambda_0$ given by
\begin{align*}
\Lambda_0 = \begin{pmatrix} \frac{\alpha}{2} \partial_x^2 & 1-\partial_x^2 \\ \frac{\alpha^2}{4} \partial_x^4 \left(1-\partial_x^2\right)^{-1} & \frac{\alpha}{2} \partial_x^2\end{pmatrix}.
\end{align*}

The spectrum of the constant-coefficient operator $\Lambda_0$ is determined by the eigenvalues $\lambda_{0,\pm}(k)$ of its Fourier symbol 
\begin{align*} \widehat{\Lambda}_0(k) = \begin{pmatrix} -\frac{\alpha}{2} k^2 & 1+k^2 \\ \frac{\alpha^2}{4} \frac{k^4}{1+k^2} &- \frac{\alpha}{2} k^2\end{pmatrix},
\end{align*}
which are given by
\begin{align*}
\lambda_{0,-}(k) = -\alpha k^2, \qquad \lambda_{0,+}(k) = 0.
\end{align*}
In particular, for $k \in \R \setminus \{0\}$ the matrix $\hL_0(k)$ is diagonalizable. For later use, we note that the associated change of basis is represented by a matrix $S(k)$, whose columns are comprised of the eigenvectors of $\hL_0(k)$, and its inverse, which are given by
\begin{align*}
S(k) = \begin{pmatrix}
\frac{2 \left(k^2+1\right)}{\alpha k^2} & -\frac{2
\left(k^2+1\right)}{\alpha k^2} \\
1 & 1 
\end{pmatrix}, \qquad S(k)^{-1} = \begin{pmatrix}
\frac{\alpha k^2}{4 k^2+4} & \frac{1}{2} \\
-\frac{\alpha k^2}{4 k^2+4} & \frac{1}{2} 
\end{pmatrix},
\end{align*}
for $k \in \R \setminus \{0\}$. That is, we have $S(k)^{-1} \hL_0(k) S(k) = \mathrm{diag}(0,-\alpha k^2)$ for $k \in \R \setminus \{0\}$. One readily observes that the coefficients of $S(\cdot)$ and $S(\cdot)^{-1}$ are bounded on $\R \setminus (-1,1)$. 

Thus, we find $\sigma(\Lambda_0) = (-\infty,0]$. So, the resolvent set $\rho(\Lambda_0)$ contains the sector $\Sigma_0 = \{\lambda \in \C : \lambda \neq 1, |\mathrm{arg}(\lambda - 1)| \leq \frac{3\pi}{4}\}$. The resolvent $(\Lambda_0 - \lambda)^{-1}$ possesses the Fourier symbol 
\begin{align*}
\frac{1}{\lambda}\begin{pmatrix}
-1 + \frac{\alpha k^2}{2(\alpha k^2 + \lambda)} & -\frac{1+k^2}{\alpha k^2 + \lambda} \\ 
-\frac{\alpha^2 k^4}{4\left(1+k^2\right)\left(\alpha k^2 + \lambda\right)} & -1 + \frac{\alpha k^2}{2(\alpha k^2 + \lambda)}
\end{pmatrix},
\end{align*}
for $\lambda \in \Sigma_0$. For $\lambda \in \Sigma_0$ and $k \in \R$ we have the basic inequalities
\begin{align*}
\left|1 + \frac{\lambda}{\alpha k^2}\right| \geq \frac{1}{\sqrt{2}}, \quad \ \ \left|\alpha k^2 + \lambda\right| \geq \frac{1}{\sqrt{2}}, \quad\ \ \frac{k^2}{1+k^2} \leq 1, \quad\ \ |\lambda - 1| \leq |\lambda| + 1 \leq |\lambda|(1+\sqrt{2}).
\end{align*}
Hence, there exists a constant $M > 0$ such that 
\begin{align*} \left\|(\Lambda_0 - \lambda)^{-1}\right\|_{Y_0} \leq \frac{M}{|\lambda - 1|}, 
\end{align*}
for all $\lambda \in \Sigma_0$. We conclude that $\Lambda_0$ is a sectorial operator on $Z_0$. By standard perturbation theory of sectorial operators~\cite[Proposition~2.4.1]{LUN}, it follows that $\Lambda$ is a sectorial operator on $Z_0$, since $\Lambda$ is a bounded perturbation of $\Lambda_0$. In addition, since the domain $D(\Lambda) = \{U \in Z_0 : \Lambda U \in Z_0\}$ obviously contains the dense subspace $C_c^\infty(\R) \subset Z_0$ of all test functions, $\Lambda$ is densely defined. 

Furthermore, $\mathcal{N}$ is locally Lipschitz continuous on $Z_0$, because $N$ is smooth, the space $Z_1$ continuously embeds into $L^\infty(\R)$ and $\smash{(1-\partial_x^2)^{-1}}$ is a bounded linear operator from $Z_2$ into $Z_0$. Combining the latter with the fact that $\Lambda$ is sectorial and densely defined on $Z_0$, it follows by standard local existence theory~\cite[Theorem~7.1.2 and Proposition~7.1.7]{LUN} for semilinear parabolic equations that there exist $T_{\max} \in (0,\infty]$ and a unique, maximally defined, mild solution $U \in C\big([0,T_{\max}),Z_0\big)$ of~\eqref{e:vKGsys} with initial condition $U(0) = U_0$. If $T_{\max} < \infty$, then we have~\eqref{e:blowup}. This establishes the result for the case $\mathcal{X} = Y_0$ or $\mathcal{X} = L^2(\R)$.

Next, we consider the case $\mathcal{X} = X_0$. We introduce the tailor made variable $V = x U$ and observe that $(U,V)$ satisfies the $4$-component system
\begin{align} \label{e:vKGsys2}
\begin{split}
U_t &= \Lambda U + \mathcal{N}(U),\\
V_t &= \Lambda V + \widetilde{\Lambda} U + \widetilde{\mathcal{N}}(U,V)
\end{split}
\end{align}
where $\widetilde{\Lambda}$ is the linear operator with Fourier symbol
\begin{align*}
\ri \hL'(k) = \begin{pmatrix} -\alpha \ri k & 2\ri k \\ \frac{\alpha^2 k^3\left(2+k^2\right)}{2\left(1+k^2\right)^2} & -\alpha \ri k
\end{pmatrix},
\end{align*}
and we denote
\begin{align*}
\widetilde{\mathcal{N}}(U,V) = \left(1-\partial_x^2\right)^{-1}\left(V_1\frac{N(U_1)}{U_1} \mathbf{e}_2\right) - 2 \partial_x \left(1-\partial_x^2\right)^{-2}\left(N(U_1) \mathbf{e}_2\right).
\end{align*}

Since $\Lambda$ is sectorial and densely defined on $Y_0$, the operator $\mathcal{L}$ on the product space $Y_0 \times Y_0$ given by $(U,V) \mapsto (\Lambda U, \Lambda V)$ is also sectorial and densely defined. Next, we show that, for any $\epsilon > 0$, the operator $(U,V) \mapsto (0,\widetilde{\Lambda}U)$ is relatively bounded with respect to $\mathcal{L}$ with $\mathcal{L}$-bound $\epsilon > 0$. First, we compute
\begin{align*}
S(k)^{-1} \hL'(k) S(k) = \begin{pmatrix}
0 & -\frac{\alpha k}{k^2+1} \\
\frac{\alpha k}{k^2+1} & -2 \alpha k \\
\end{pmatrix},
\end{align*}
for $k \in \R \setminus \{0\}$. Hence, abbreviating $S(k)^{-1} W = (w_1,w_2)$ and using that $S(\cdot)$, $S(\cdot)^{-1}$ and $k \mapsto k/(1+k^2)$ are bounded on $\R \setminus (-1,1)$, we establish
\begin{align*}
\left|\hL'(k) W\right| &\leq \left|S(k)\right| \left|S(k)^{-1} \hL'(k) S(k) S(k)^{-1} W\right| \lesssim |k| \left|w_2\right| + \left|S(k)^{-1} W\right|\\
&\lesssim \left|S(k)^{-1} \hL_0(k) S(k) S(k)^{-1} W\right|^{\frac12} \left|S(k)^{-1} W\right|^{\frac12} + \left|S(k)^{-1}\right| \left|W\right| \\
&\lesssim \left|S(k)^{-1}\right| \left|\hL_0(k) W\right|^{\frac12} \left|W\right|^{\frac12} + \left|S(k)^{-1}\right| \left|W\right|\\
&\lesssim \left|\hL_0(k) W\right|^{\frac12} \left|W\right|^{\frac12} + \left|W\right| \lesssim \left|\hL(k) W\right|^{\frac12} \left|W\right|^{\frac12} + \left|W\right|,
\end{align*}
for $W \in \C^2$ and $k \in \R \setminus (-1,1)$. Combining the latter with the fact that $\hL'(\cdot)$ is bounded on $[-1,1]$, we obtain
\begin{align}\label{e:relbounded}
\begin{split}
\left|\hL'(k) W\right|\lesssim \left|\hL(k) W\right|^{\frac12} \left|W\right|^{\frac12} + \left|W\right|,
\end{split}
\end{align}
for all $W \in \C^2$ and $k \in \R$. Thus, by employing Young's inequality we infer that, for any $\epsilon > 0$, the operator $\widetilde{\Lambda}$ is relatively bounded with respect to $\Lambda$ with $\Lambda$-bound $\epsilon$. Consequently, for any $\epsilon > 0$, the operator $(U,V) \mapsto (0,\widetilde{\Lambda}U)$ is relatively bounded with respect to $\mathcal{L}$  with $\mathcal{L}$-bound $\epsilon > 0$. So, it follows by~\cite[Proposition~2.4.2]{LUN} that 
\begin{align*}
\begin{pmatrix}
U \\ V
\end{pmatrix} \mapsto \begin{pmatrix}
\Lambda U \\ \Lambda V + \widetilde{\Lambda} U
\end{pmatrix},
\end{align*}
is a sectorial and densely defined operator on $Y_0 \times Y_0$ and, thus, generates an analytic semigroup on $Y_0 \times Y_0$.

Next, we note that $\widetilde{\mathcal{N}} \colon Y_0 \times Y_0 \to Y_0$ is locally Lipschitz continuous, since the map $u \mapsto N(u)/u$ is smooth and because $\partial_x (1-\partial_x^2)^{-2}$ and $(1-\partial_x^2)^{-1}$ are bounded linear operators on $Y_0$. We conclude that
\begin{align*}
\begin{pmatrix}
U \\ V
\end{pmatrix} \mapsto \begin{pmatrix}
\mathcal{N}(U) \\ \widetilde{\mathcal{N}}(U,V)
\end{pmatrix}
\end{align*}
is locally Lipschitz continuous on $Y_0 \times Y_0$. It follows again by~\cite[Theorem~7.1.2 and Proposition~7.1.7]{LUN} that there exist $T_{\max} \in (0,\infty]$ and a unique, maximally defined, mild solution 
\begin{align*}
\begin{pmatrix}U \\ V \end{pmatrix} \in C\big([0,T_{\max}),Y_0 \times Y_0\big),
\end{align*}
of~\eqref{e:vKGsys2} with initial condition $U(0) = (U_0,V_0)^\top$, where $V_0(x) = x U_0(x)$. If $T_{\max} < \infty$, then it holds 
\begin{align*}
\limsup_{t \uparrow T_{\max}} \|(U(t),V(t))\|_{Y_0 \times Y_0} = \infty.
\end{align*}
Upon recalling $V(x,t) = x U(x,t)$, we conclude that $U \in C\big([0,T_{\max}),X_0\big)$ is a solution of~\eqref{e:vKGsys} with initial condition $U(0) = U_0$. Moreover, if $T_{\max} < \infty$, then we have~\eqref{e:blowup}.
\end{proof}

The next result relates the mild solution $U(t)$ of~\eqref{e:vKGsys}, established in Proposition~\ref{p:wellposedness}, to a classical solution of the viscoelastic Klein-Gordon equation~\eqref{e:vKG2} under the additional regularity assumption that $u_0 \in H^2(\R)$.

\begin{proposition} \label{p:wellposedness2}
Let $\mathcal{X}$ be one of the spaces $X_0,Y_0$ or $L^2(\R)$. Take $u_0 \in \mathcal{X} \cap H^2(\R)$ and $w_0 \in \mathcal{X}$. Set 
$v_0 = (1-\partial_x^2)^{-1} (w_0 - \frac{\alpha}{2} \partial_x^2 u_0) \in \mathcal{X}$. Let $u(t)$ be the first coordinate of the unique maximal mild solution $U \in C\big([0,T_{\max}),\mathcal{X}\big)$, established in Proposition~\ref{p:wellposedness}, with initial condition $U(0) = (u_0,v_0)^\top \in \mathcal{X}$. If $u_0 \in H^2(\R)$, then we have that
\begin{align*}
u &\in C\big([0,T_{\max}),H^2(\R)\big) \cap C^1\big([0,T_{\max}),L^2(\R)\big) \cap C^1\big((0,T_{\max}),H^2(\R)\big)\\ 
&\qquad \cap C^2\big((0,T_{\max}),L^2(\R)\big)
\end{align*}
is a classical solution of the viscoelastic Klein-Gordon equation~\eqref{e:vKG2} with initial condition $u(0) = u_0$ and $u_t(0) = w_0$.
\end{proposition}
\begin{proof}
Using again the variable $Z = (u,u_t)^\top$ we write~\eqref{e:vKG2} as the system
\begin{align} \label{e:semi}
Z_t = \mathcal{A} Z + F(Z),
\end{align}
where $\mathcal{A}$ is the differential operator given by~\eqref{e:defA} defined on the space $H^2(\R) \times L^2(\R)$ with dense domain $H^2(\R) \times H^2(\R)$ and the nonlinear map $F \colon H^2(\R) \times L^2(\R) \to H^2(\R) \times L^2(\R)$ is given by
\begin{align*}
F(z_1,z_2) = \begin{pmatrix} 0 \\ N(z_1)\end{pmatrix}.
\end{align*}
Thanks to the smoothness of $N$ and the fact that $H^1(\R)$ continuously embeds into $L^\infty(\R)$, we find that $F$ is well-defined and locally Lipschitz continuous. Moreover, since $\partial_x^2$ is a sectorial operator on $L^2(\R)$ with dense domain $H^2(\R)$,~\cite[Proposition~2.2]{WE80} yields that $\mathcal{A}$ generates an analytic semigroup on $H^2(\R) \times L^2(\R)$. Standard analytic semigroup theory~\cite[Proposition~7.1.8 and~7.1.10]{LUN} now provides a time $\tau_{\max} \in (0,\infty]$ and a unique, maximally defined, classical solution 
\begin{align} \label{e:regZ}
\begin{split}
Z &\in C\big([0,\tau_{\max}),H^2(\R) \times L^2(\R)\big) \cap  C^1\big((0,\tau_{\max}),H^2(\R) \times L^2(\R)\big)\\ 
&\qquad \cap C\big((0,\tau_{\max}),H^2(\R) \times H^2(\R)\big),
\end{split}
\end{align}
of~\eqref{e:semi} with initial condition $Z(0) = (u_0,w_0) \in H^2(\R) \times L^2(\R)$ such that, if $\tau_{\max} < \infty$, then it must hold
\begin{align} \label{e:Fblowup}
\limsup_{t \uparrow \tau_{\max}} \|F(Z(t))\|_{H^2 \times L^2} = \infty.
\end{align}
Writing $Z(t) = (u(t),w(t))^\top$ one readily observes that $\partial_t u(t) = w(t)$ for all $t > 0$. Combining the latter with~\eqref{e:regZ} we arrive at
\begin{align*}
\begin{split}
u &\in C\big([0,\tau_{\max}),H^2(\R)\big) \cap C^1\big([0,\tau_{\max}),L^2(\R)\big) \cap C^1\big((0,\tau_{\max}),H^2(\R)\big)\\ 
&\qquad \cap C^2\big((0,\tau_{\max}),L^2(\R)\big).
\end{split}
\end{align*}
Finally, due to the continuity of $N$ and the fact that $\|F(z_1,z_2)\|_{H^2 \times L^2} = \|N(z_1)\|_{L^2}$ for $(z_1,z_2) \in H^2(\R) \times L^2(\R)$, identity~\eqref{e:Fblowup} implies
\begin{align} \label{e:L2blowup}
\limsup_{t \uparrow \tau_{\max}} \|u(t)\|_{L^2} = \infty.
\end{align}
Define  $\check{U}(t) = (u(t),v(t))$ with  $v(t) = (1-\partial_x^2)^{-1}(w(t) - \frac{\alpha}{2} \partial_x^2 u(t))$. Then, by construction $\check{U} \in C\big([0,\tau_{\max}), L^2(\R)\big)$ is a mild solution of~\eqref{e:vKGsys} with initial condition $\check{U}(0) = (u_0,v_0)^\top$. On the other hand, $U \in C\big([0,T_{\max}), L^2(\R)\big)$, established in Proposition~\ref{p:wellposedness}, is also a mild solution of~\eqref{e:vKGsys} with $U(0) = (u_0,v_0)^\top$. Recalling from the proof of Proposition~\ref{p:wellposedness} that $\Lambda$ generates a $C_0$-semigroup on $L^2(\R)$ and the nonlinearity $\mathcal{N}$ is locally Lipschitz continuous on $L^2(\R)$, it follows by uniqueness of mild solutions, cf.~\cite[Lemma~4.3.2]{CAHA}, that $\check{U}(t) = U(t)$ for all $t \in [0,\min\{T_{\max},\tau_{\max}\})$. Since~\eqref{e:L2blowup} implies~\eqref{e:blowup}, we must necessarily have $\tau_{\max} \geq T_{\max}$.
\end{proof}

\section{Linear estimates} \label{sec:linear}

In order to exploit dispersive decay exhibited by the critical nonlinear terms in~\eqref{e:vKG2}, we perform our nonlinear argument in Fourier space. In this section, we study the linear dynamics of~\eqref{e:vKG2} in Fourier space. That is, we study the Fourier symbol 
\begin{align} \widehat{\Lambda}(k) = \begin{pmatrix} -\frac{\alpha}{2} k^2 & 1+k^2 \\ -1 + \frac{\alpha^2}{4} \frac{k^4}{1+k^2} &- \frac{\alpha}{2} k^2\end{pmatrix}, \label{e:def Fourier symbol}
\end{align}
of the linear operator $\Lambda$ in~\eqref{e:vKGsys} and obtain estimates on the matrix exponential $\re^{\widehat{\Lambda}(k) t}$, which represents the pointwise action of the semigroup generated by $\Lambda$ in Fourier space. The eigenvalues $\lambda_\pm(k)$ of $\widehat{\Lambda}(k)$ read
\begin{align*}
\lambda_{\pm}(k) = -\dfrac{1}{2}\alpha k^2 \pm \mu(k), \qquad \mu(k) = \sqrt{\dfrac{1}{4}\alpha^2k^4 - 1 - k^2},
\end{align*}
and obey the expansion
\begin{align}
\label{e:eigenvalue expansion}
\lambda_{\pm}(k) = \pm \ri + \dfrac{1}{2}\left(- \alpha \pm \ri\right) k^2 + \bigo\left(k^4\right).
\end{align}

We note that the curves $\lambda_\pm \colon \R \to \C$ are confined to the left-half plane and touch the imaginary axis only at the points $\pm \ri$ for the same frequency $k = 0$. Since $\hL(k)$ is real-valued for $k \in \R$, we find that for small values of $k$, the eigenvalues $\lambda_\pm(k)$ are complex conjugates lying in $\C \setminus \R$, and, thus, the Fourier symbol $\widehat{\Lambda}(k)$ is diagonizable. Hence, there exist $k_0 > 0$ and smooth maps $P_\pm \colon (-k_0,k_0) \to \C^{2 \times 2}$ such that $P_\pm(k)$ is the spectral projection of $\widehat{\Lambda}(k)$ onto the $1$-dimensional eigenspace corresponding to the eigenvalue $\lambda_\pm(k)$ for $k \in (-k_0,k_0)$. Moreover, it holds
\begin{align} \label{e:conjugateids}
\overline{\lambda_\mp(k)} = \lambda_\pm(k) = \lambda_\pm(-k), \qquad \overline{P_\mp(k)} = P_\pm(k) = P_\pm(-k),
\end{align}
for $k \in (-k_0,k_0)$, because $\hL(k)$ is real-valued and analytic in $k^2$. For later use, we explicitly compute the spectral projections at the critical frequency and obtain
\begin{align} \label{e:specprojatcrit}
P_\pm(0) = \frac{1}{2} \begin{pmatrix} 1 & \mp \ri \\ \pm \ri & 1 \end{pmatrix}.
\end{align}

For small Fourier frequencies the linear dynamics is diffusive, which is represented by the estimate
\begin{align} \Re(\lambda_\pm(k)) = -\frac{1}{2} \alpha k^2, \qquad k \in [-k_0,k_0].\label{e:linest2}\end{align}
On the other hand, for Fourier frequencies away from the critical frequency $k = 0$, the linear dynamics is exponentially damped. More precisely, there exists $\theta_1 > 0$ such that
\begin{align} \label{e:linest1}
\sup \Re(\sigma(\widehat{\Lambda}(k))) < -\theta_1, \qquad k \in \R \setminus \left(-\tfrac{k_0}2,\tfrac{k_0}2\right).
\end{align}
We employ~\eqref{e:linest2} and~\eqref{e:linest1} to bound the exponential $\re^{\widehat{\Lambda}(k) t}$ for small and for noncritical frequencies, respectively.

\begin{lemma} \label{l:linear estimates}
There exists $\theta_0 > 0$ such that
\begin{align} 
t^j \left|\partial_t^j \re^{\widehat{\Lambda}(k) t}\right| \lesssim \re^{-\theta_0 t}, \qquad \left|\partial_k \re^{\widehat{\Lambda}(k) t}\right| \lesssim \re^{-\theta_0 t},
\label{e:matrixexp-high}
\end{align}
for $j = 0,1$, $t \geq 0$ and $k \in \R \setminus [-\frac{k_0}{2},\frac{k_0}{2}]$. Moreover, we have
\begin{align}
\left|\re^{\widehat{\Lambda}(k) t}\right| &\lesssim \re^{-\frac12 \alpha k^2 t}, \qquad \left|\partial_k \re^{\widehat{\Lambda}(k) t}\right| \lesssim |k| t \re^{-\frac12 \alpha k^2 t},\label{e:matrixexp-low}
\end{align}
for all $t \geq 0$ and $k \in [-k_0,k_0]$.
\end{lemma}
\begin{proof}
Let
$$k_1 = \sqrt{\frac{2 + 2\sqrt{1 + \alpha^2}}{\alpha^2}}.$$
The eigenvalues $\lambda_+(k)$ and $\lambda_-(k)$ are distinct for $k^2 > k_1^2$ and, thus, $\widehat{\Lambda}(k)$ can be diagonalized for such values. The associated change of basis is represented by a matrix $S(k)$, whose columns are comprised of eigenvectors of $\widehat{\Lambda}(k)$, and its inverse, which are given by
\begin{align*} S(k) = \begin{pmatrix} \frac{1+k^2}{\mu(k)} & -\frac{1+k^2}{\mu(k)}\\ 1 & 1 \end{pmatrix}, \qquad S(k)^{-1} = \frac12\begin{pmatrix} \frac{\mu(k)}{4 k^2+4} & 1 \\
-\frac{\mu(k)}{4 k^2+4} & 1 \end{pmatrix}.
\end{align*}
One readily observes that the coefficients of $S(\cdot)$ and $S(\cdot)^{-1}$ are bounded on $\R \setminus (-2k_1,2k_1)$. Moreover, it holds
\begin{align*}
\lambda_-(k) \leq -\dfrac{1}{2}\alpha k^2, \qquad \lambda_+(k) \leq -\dfrac{1}{\alpha}
\end{align*}
for $|k| > k_1$. We conclude that there exists $\theta_2 > 0$ such that the matrix exponential and its temporal derivative
\begin{align*}
\partial_t^j \re^{\widehat{\Lambda}(k) t} = S(k)^{-1} \mathrm{diag}\left(\lambda_+(k)^j \re^{\lambda_+(k) t}, \lambda_-(k)^j\re^{\lambda_-(k) t}\right) S(k),
\end{align*}
obey the estimate
\begin{align} \label{e:matrixexp}
t^j \left|\partial_t^j \re^{\widehat{\Lambda}(k) t}\right| \lesssim \re^{-\theta_2 t},
\end{align}
for $k \in \R \setminus (-2k_1,2k_1)$, $j = 0,1$ and $t \geq 0$. To bound the matrix exponential $\re^{\widehat{\Lambda}(k) t}$ and its temporal derivative on the compact set $\smash{J \deq [-2k_1,-k_0/2] \cup [k_0/2,2k_1]}$ we collect some facts from~\cite[Chapter A-III, \S7]{ARE84}. First, since $J$ is compact and $\widehat{\Lambda}$ is continuous on $J$, the multiplication operator $A \colon f \mapsto \widehat{\Lambda} f$ generates a strongly continuous semigroup $(T(t))_{t \geq 0}$ on $C(J,\C^2)$, which is given by
$$(T(t)f)(k) = \re^{\widehat{\Lambda}(k) t} f(k), \qquad k \in J.$$
Second, the growth bound of the semigroup $(T(t))_{t \geq 0}$ coincides with the spectral bound of $A$. Third, the spectrum of $A$ is given by
$$\sigma(A) = \bigcup_{k \in J} \sigma(\widehat{\Lambda}(k)).$$
Combining the latter three observations with~\eqref{e:linest1} yields that the growth bound of the semigroup $(T(t))_{t \geq 0}$ is smaller than $-\theta_1$, which implies 
\begin{align*}
\left|\re^{\widehat{\Lambda}(k) t}\right| \lesssim \re^{-\theta_1 t}, \qquad \left|\partial_t\re^{\widehat{\Lambda}(k) t}\right| = \left|\hL(k) \re^{\widehat{\Lambda}(k) t}\right| \lesssim \re^{-\theta_1 t},
\end{align*}
for all $t \geq 0$ and $k \in J$. Combining the latter with~\eqref{e:matrixexp} yields a constant $\theta_3 > 0$ such that 
\begin{align} 
t^j \left|\partial_t^j \re^{\widehat{\Lambda}(k) t}\right| \lesssim \re^{-\theta_3 t},
\label{e:matrixexp-high2}
\end{align}
for $j = 0,1$, $t \geq 0$ and $k \in \R \setminus [-\frac{k_0}{2},\frac{k_0}{2}]$, proving the first estimate in~\eqref{e:matrixexp-high}.

We proceed with the second estimate in~\eqref{e:matrixexp-high}. First, we express the derivative as
\begin{align} \label{e:exprder}
\partial_k \re^{\widehat{\Lambda}(k) t} = t\int_0^1 \re^{\widehat{\Lambda}(k) l t} \hL'(k) \re^{\widehat{\Lambda}(k) (1-l) t} \de l.
\end{align}
Next, we apply the estimate~\eqref{e:relbounded} and obtain
\begin{align*}
\left|\partial_k \re^{\widehat{\Lambda}(k) t} W\right| \leq t \int_0^1 \left|\re^{\widehat{\Lambda}(k) l t} \right| \left|\hL(k) \re^{\widehat{\Lambda}(k) (1-l) t}W\right|^{\frac12} \left|\re^{\widehat{\Lambda}(k) (1-l) t}W\right|^{\frac12} \de l,
\end{align*}
for $W \in \C^2$, $k \in \R$ and $t \geq 0$. So, using~\eqref{e:matrixexp-high2} we arrive at
\begin{align*}
\left|\partial_k \re^{\widehat{\Lambda}(k) t} W\right| \lesssim t \re^{-\theta_0 t} \int_0^1 \frac{1}{\sqrt{(1-l)t}} \de l \left|W\right| \lesssim \sqrt{t}\, \re^{-\theta_3 t} \left|W\right| \lesssim \re^{-\frac{\theta_3}{2} t} \left|W\right|,
\end{align*}
for $W \in \C^2$, $k \in \R \setminus [-\frac{k0}{2},\frac{k0}{2}]$ and $t \geq 0$, which proves the second estimate in~\eqref{e:matrixexp-high}.

Finally, we obtain bounds for $k\in [-k_0, k_0]$. For such values, the matrix $\smash{\widehat{\Lambda}(k)}$ is diagonalizable and both $S(\cdot)$ and $S^{-1}(\cdot)$ are bounded on $[-k_0,k_0]$. Moreover, we have $|\hL'(k)| \lesssim |k|$ for $k \in [-k_0,k_0]$. Combining the latter two observations with~\eqref{e:linest2} and~\eqref{e:exprder} readily yield the last estimate~\eqref{e:matrixexp-low}, which finishes the proof.
\end{proof}

\begin{figure}[h] 
\begin{center}
\includegraphics[scale=0.91,trim={0.1cm 0.01cm 0.1cm 0.1cm},clip]{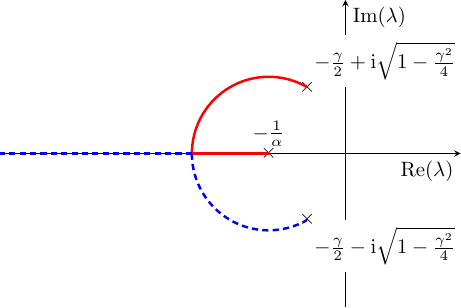} \hspace{0.5cm} \includegraphics[scale=0.91,trim={0.1cm 0.01cm 0.1cm 0.1cm},clip]{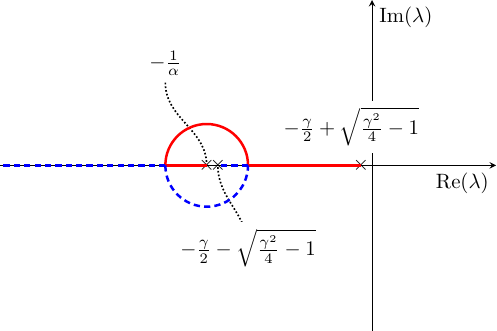}
\end{center}
\caption{Depiction of the set~\eqref{e:speccurve2} for $0 < \gamma < 2$ (left panel) and for $\gamma > 2$ (right panel) under the condition that $\alpha^2 - \alpha \gamma + 1 > 0$. In case $\alpha^2 - \alpha \gamma + 1 \leq 0$, the set \eqref{e:speccurve2} is confined to the negative real line.  \label{fig:speccurve2}}
\end{figure}

\begin{remark} \label{rem:damping}
In addition to the viscous damping term $-\alpha \partial_x^2 u_t$, it is possible to include the first-order damping term $\gamma u_t$ with $\gamma > 0$ on the left-hand side of~\eqref{e:vKG2}, cf.~\cite{AVI,Potier-Ferry}. This dissipative term models internal resistance or frictional forces in the viscoelastic material. In contrast to the viscous damping term $-\alpha \partial_x^2 u_t$, it dampens the amplitude of solutions without smoothing out oscillations. In the presence of the first-order damping term $\gamma u_t$ the spectrum of the linearization is given by
\begin{align} \label{e:speccurve2}
\left\{-\dfrac{1}{2}\left(\gamma + \alpha k^2\right) \pm  \sqrt{\dfrac{1}{4}\left(\gamma + \alpha k^2\right)^2 - 1 - k^2} : k \in \R\right\},
\end{align}
see Figure~\ref{fig:speccurve2}. In contrast to our situation, the spectrum is confined to the open left-half plane for $\gamma > 0$ and it follows from standard arguments, cf.~\cite[Section~10]{CAHA}, that the equilibrium state $u(x,t) = 0$ in~\eqref{e:vKG2} is exponentially stable for $\gamma > 0$. That is, solutions with small initial data decay with exponential rate. 
\end{remark}

\section{Mode filters} \label{sec:modefilters}

We are ready to start the reduction process that is necessary to close the nonlinear argument. Let $\mathcal{X}$ denote the space $X_0$ or $Y_0$. Take $u_0,w_0 \in \mathcal{X}$ satisfying~\eqref{e:cond_initial} (in case $\mathcal{X} = X_0$) or satisfying~\eqref{e:cond_initial2} (in case $\mathcal{X} = Y_0$). Clearly, there exists a constant $K > 0$, independent of $u_0$ and $w_0$, such that
\begin{align*} U_0 \deq \begin{pmatrix}u_0 \\ (1-\partial_x^2)^{-1}\left(w_0 - \frac{\alpha}{2} \partial_x^2 u_0\right)\end{pmatrix} \in \mathcal{X},\end{align*}
obeys
\begin{align} \label{e:cond_initial3}
\|U_0\|_{\mathcal{X}} \leq KE_0.
\end{align}
Moreover, Proposition~\ref{p:wellposedness} yields a time $T_{\max} \in (0,\infty]$ and a unique, maximally defined, mild solution
\begin{align} \label{e:regularity2}
U \in C\big([0,T_{\max}),\mathcal{X}\big),
\end{align}
of~\eqref{e:vKGsys} with initial condition $U(0) = U_0$ such that, if $T_{\max} < \infty$, then~\eqref{e:blowup} holds. 

In Fourier space system~\eqref{e:vKGsys} reads
\begin{align} \label{e:vKGsysFourier}
\partial_t \U = \widehat{\Lambda} \U + \mathcal{F} \mathcal{N}\big(\mathcal{F}^{-1} \U\big),
\end{align}
where $\widehat{\Lambda}$ is the Fourier symbol of $\Lambda$ defined in~\eqref{e:def Fourier symbol}. 

To separate diffusive low-frequency modes from exponentially damped high-frequency modes and diagonalize the system at criticality, we introduce mode filters. Thus, let $k_0 > 0$ be as in~\S\ref{sec:linear} and let $\chi \colon \R \to [0,1]$ be a smooth even cut-off function, whose support is contained in $(-k_0,k_0)$, such that $\chi(k) = 1$ for $k \in [-k_0/2,k_0/2]$. We decompose the solution
\begin{align} \label{e:Fdef}
\U(t) = \chi \U(t) + (1-\chi)\U(t) \deq \Uc(t) + \Us(t), \qquad t \in [0,T_{\max})
\end{align}
of~\eqref{e:vKGsys} into a low- and high-frequency part. We note that, since $U(t)$ is real-valued and $\chi$ is even, it holds
\begin{align} \label{e:conjUC}
\smash{\overline{\Uc(k,t)}} = \Uc(-k,t)
\end{align}
for each $k \in \R$ and $t \in [0,T_{\max})$. Moreover, by~\eqref{e:regularity2}, smoothness of $\chi$ and the fact that $\chi$ is compactly supported, we obtain
\begin{align} \label{e:regularity3}
\U, \Uc, |\cdot| \Uc, \Us \in C\big([0,T_{\max}),W^{1,1}(\R) \cap W^{1,\infty}(\R)\big),
\end{align}
in case $\mathcal{X} = X_0$, and 
\begin{align} \label{e:regularity32}
\U, \Uc, \Us \in C\big([0,T_{\max}),L^1(\R) \cap L^\infty(\R)\big),
\end{align}
in case $\mathcal{X} = Y_0$.

Multiplying~\eqref{e:vKGsysFourier} with $\chi$ and $1-\chi$, we arrive at the system
\begin{align} \label{e:vKGsysFourierCrit}
\partial_t \Uc &= \widehat{\Lambda} \Uc + \chi \mathcal{F} \mathcal{N}\left(\mathcal{F}^{-1}\big(\Uc + \Us\big)\right),\\
\partial_t \Us &= \widehat{\Lambda} \Us + (1-\chi) \mathcal{F} \mathcal{N}\left(\mathcal{F}^{-1}\big(\Uc + \Us\big)\right),
\label{e:vKGsysFourierDamp}
\end{align}
for the new variables $\Uc$ and $\Us$. 

The Duhamel formulation associated to \eqref{e:vKGsysFourierCrit}-\eqref{e:vKGsysFourierDamp} reads
\begin{align} 
\Uc(k,t) &= \re^{\widehat{\Lambda}(k) t} \chi(k) \U(k,0) + \int_0^t \re^{\widehat{\Lambda}(k)(t-s)} \chi(k) \mathcal{F} \mathcal{N}\left(\mathcal{F}^{-1}\big(\Uc(s) + \Us(s)\big)\right)\!(k) \, \de s, \label{e:Duhcrit}\\
\begin{split}
\Us(k,t) &= \re^{\widehat{\Lambda}(k) t} (1-\chi(k)) \U(k,0)\\ 
& \qquad + \, \int_0^t \re^{\widehat{\Lambda}(k)(t-s)} (1-\chi(k)) \mathcal{F} \mathcal{N}\left(\mathcal{F}^{-1}\big(\Uc(s) + \Us(s)\big)\right)\!(k) \, \de s,\end{split} \label{e:Duhdamp}
\end{align}
with $k \in \R$ and $t \in [0,T_{\max})$. 

\section{Separating relevant from irrelevant nonlinear terms} \label{sec:isolate}

As outlined before, our goal is to close a nonlinear argument through iterative estimates on the Duhamel formulation~\eqref{e:Duhcrit}-\eqref{e:Duhdamp}. In this section, we isolate those nonlinear terms in~\eqref{e:Duhcrit}-\eqref{e:Duhdamp}, which cannot be controlled through standard iterative $L^1$-$L^\infty$-estimates, and establish estimates on the irrelevant residual terms. 

Lemma~\ref{l:linear estimates} implies that $\smash{\re^{\widehat{\Lambda}(k) t} (1-\chi(k))}$ is exponentially decaying over time, whereas $\smash{\re^{\widehat{\Lambda}(k) t} \chi(k)}$ only decays diffusively. Therefore, one expects that the decay of $\Us(t)$ is dictated by the slowest decaying nonlinear terms in $\Uc$, which are the quadratic terms. Hence, $\Us(t)$ is expected to decay at higher rate than $\Uc(t)$. Thanks to the higher decay rate of $\Us(t)$ in combination with the fact that it vanishes at the critical frequency $k = 0$, any nonlinear term in~\eqref{e:Duhcrit} with a  $\Us$-contribution can be controlled, as can any quartic or higher-order nonlinear term in~\eqref{e:Duhcrit} as outlined in~\S\ref{sec:technical}. All in all, the only nonlinear terms that cannot controlled through standard iterative $L^1$-$L^\infty$-estimates are the quadratic and cubic terms in~\eqref{e:Duhcrit}-\eqref{e:Duhdamp}. 

We split off these terms by expanding $N \in C^4(\R)$ as a Taylor series
\begin{align}
\label{e:def-N}
N(u) = 2\pi\kappa u^2 + 4\pi^2\beta u^3 + R(u)
\end{align}
with coefficients $\beta,\kappa \in \R$, where there exists a constant $C > 0$ such that the quartic remainder $R\colon \R \to \R$ obeys the estimate 
\begin{align*}
|R(u)| \leq C|u|^4, 
\end{align*}
for $u \in [-1,1]$. The critical quadratic term in~\eqref{e:Duhcrit} is then given by
\begin{align} \label{e:critquadr} 
\int_0^t \chi(k) \re^{\widehat{\Lambda}(k) (t-s)} B_2\big(\Uc(s),\Uc(s)\big)(k) \de l \de s,
\end{align}
where $B_2$ is the symmetric bilinear form on $
\mathcal{X}$ given by
\begin{align*}
B_2(\V, \widehat{W})(k) = \int_\R N_2(k,l)\big(\V,\W) \de l, \qquad N_2(k,l)\big(\V,\W\big) \deq \frac{\kappa}{1+k^2}  \V_1(k-l)  \widehat{W}_1(l) \mathbf{e}_2,
\end{align*}
Moreover, the critical cubic term in~\eqref{e:Duhcrit} reads
\begin{align} \label{e:critcubic} 
\int_0^t \chi(k) \re^{\widehat{\Lambda}(k) (t-s)} B_3\big(\Uc(s),\Uc(s),\Uc(s)\big)(k) \de s,
\end{align}
where $B_3$ is the symmetric trilinear form on $\mathcal{X}$ given by
\begin{align*}
B_3(\V, \widehat{W}, \widehat{Z})(k) &= \int_\R \int_\R N_3(k,l_1,l_2)\big(\V,\W,\widehat{Z}\big) \de l_1 \de l_2, \\
N_3(k,l_1,l_2)\big(\V,\W,\widehat{Z}\big) &\deq \frac{\beta}{1+k^2} \V_1(k-l_1) \widehat{W}_1(l_1-l_2) \widehat{Z}_1(l_2) \mathbf{e}_2.
\end{align*}
Finally, the nonlinear remainder in~\eqref{e:Duhcrit} is
\begin{align*}
\int_0^t \re^{\widehat{\Lambda}(k)(t-s)} \chi(k) \mathcal{E}\big(\Uc(s),\Us(s)\big)(k) \de s,
\end{align*}
where $\mathcal{E} \colon \mathcal{X} \times \mathcal{X} \to \mathcal{X}$ is the nonlinear operator given by
\begin{align*}
\mathcal{E}\big(\V,\W\big) &= \mathcal{F} \mathcal{R}\left(\mathcal{F}^{-1}\big(\V + \W\big)\right) + 2B_2\big(\V,\W\big) + B_2\big(\W,\W\big)\\ 
&\qquad \, + 3B_3\big(\V,\V,\W)\big) + 3B_3\big(\V,\W,\W\big) + B_3\big(\W,\W,\W\big),
\end{align*}
with 
\begin{align*}
\mathcal{R}(V) = \left(1-\partial_x^2\right)^{-1} R(V_1) \mathbf{e}_2.
\end{align*}

The next result establishes estimates on the noncritical nonlinear terms in~\eqref{e:Duhcrit}-\eqref{e:Duhdamp}. 

\begin{lemma} \label{l:noncriticalbounds}
We have
\begin{align*}
\left\|\partial_k^j \mathcal{F} \mathcal{N}\left(\mathcal{F}^{-1}\big(\V\big)\right)\right\|_{L^p} &\lesssim \big\|\V\big\|_{W^{j,1}} \big\|\V\big\|_{L^p},
\end{align*}
for $j = 0,1$, $p = 1,\infty$ and $\V \in L^\infty(\R) \cap W^{j,1}(\R)$ satisfying $\|\V\|_{L^1} \leq 1$. Moreover, it holds
\begin{align*}
\left|\partial_k^j \mathcal{E}\big(\V,(1-\chi)\W\big)(k)\right| &\lesssim \big\|\V\big\|_{L^1}^2 \big\|\V\big\|_{W^{j,1}}  + \big\|(1-\chi)\W\big\|_{W^{j,1}}\big\|(1-\chi)\W\big\|_{L^\infty} \\ 
&\qquad + \, |k| \big\|(1-\chi)\W\big\|_{L^\infty} \big\|\V\big\|_{W^{j,1}} + \big\|(1-\chi)\W\big\|_{L^\infty} \big\||\cdot|\V\big\|_{W^{j,1}}\\ 
&\qquad + \, \big\|(1-\chi)\W\big\|_{L^1} \big\|\V\big\|_{W^{j,1}},
\end{align*}
for $j = 0,1$, $k \in \R$ and $\V,\W \in W^{j,1}(\R) \cap L^\infty(\R)$ satisfying $\|\V\|_{L^1 \cap L^\infty}, \|\W\|_{L^1 \cap L^\infty} \leq 1$ and $\||\cdot| \V\|_{W^{j,1}} < \infty$.
\end{lemma}
\begin{proof}
Recall that $1-\chi$ vanishes on $[-k_0/2,k_0/2]$. So, letting $\upsilon \colon \R \to [0,1]$ be a smooth function supported on $[-k_0/2,k_0/2]$ with $\upsilon(0) = 1$, we obtain 
\begin{align*}
B_2\big(\V,(1-\chi)\W\big)(k) = \int_\R \frac{1-\upsilon(l)}{l} (k + (l-k)) N_2(k,l)\big(\V,(1-\chi)\W\big) \de l.
\end{align*}
for $k \in \R$. Upon noticing that $l \mapsto \frac{1-\upsilon(l)}{l}$ is bounded on $\R$, we deduce
\begin{align} \label{e:resbound0}
\left|\partial_k^j B_2\big(\V,(1-\chi)\W\big)(k)\right| \lesssim |k| \big\|\V\big\|_{W^{j,1}} \big\|\W\big\|_{L^\infty} + \big\||\cdot| \V\big\|_{W^{j,1}} \big\|\W\big\|_{L^\infty},
\end{align}
for $j = 0,1$, $k \in \R$, $\V \in W^{j,1}(\R)$ and $\W \in L^\infty(\R)$ with $\||\cdot|\V\|_{W^{j,1}} < \infty$. On the other hand, using Young's convolution inequality, we infer
\begin{align}
\label{e:resbound1} 
\begin{split}
\left\|\partial_k^j B_2\big(\V,\W\big) \right\|_{L^\infty} &\lesssim \big\|\V\big\|_{W^{j,1}} \big\|\W\big\|_{L^\infty},\\
\left\|\partial_k^j B_2\big(\V,\widehat{Z}\big) \right\|_{L^1} &\lesssim \big\|\V\big\|_{W^{j,1}} \big\|\widehat{Z}\big\|_{L^1},\\
\left\|\partial_k^j B_3\big(\V,\W,\widehat{Z}\big)\right\|_{L^\infty} &\lesssim \big\|\V\big\|_{W^{j,1}} \big\|\W\big\|_{L^\infty} \big\|\widehat{Z}\big\|_{L^1},
\end{split}
\end{align}
for $j = 0,1$, $\V \in W^{j,1}(\R)$, $\W \in L^\infty(\R)$ and $\widehat{Z} \in L^1(\R)$.

Next, we expand
\begin{align} \label{e:decompN}
\mathcal{F} \mathcal{N}\big(\mathcal{F}^{-1}\V\big)\big) = B_2\big(\V,\V\big) +  \mathcal{F} \mathcal{R}_2\left(\mathcal{F}^{-1}\big(\V\big)\right),
\end{align}
with $\mathcal{R}_2(V)$ defined by
\begin{align} \label{e:defR2}
\mathcal{R}_2(V) = (1-\partial_x^2)^{-1} R_0(V_1)\mathbf{e}_2, \qquad R_0(z) \deq N(z) - 2\pi \kappa z^2.
\end{align}
The expansion~\eqref{e:def-N} yields $|R_0(z)| \lesssim |z|^3$ and $|R(z)| \lesssim |z|^4$ for $z \in [-1,1]$. So, using the facts that we have $\|V\|_{L^2}^2 \leq \|V\|_{L^1} \|V\|_{L^\infty}$ for $V \in L^1(\R) \cap L^\infty(\R)$ and the Fourier transform is an isomorphism on $L^2(\R)$ and it maps $L^1(\R)$ continuously into $L^\infty(\R)$, we establish
\begin{align} \label{e:resbound2}
\begin{split}
\left\|\partial_k^j \mathcal{F} \mathcal{R}_2\left(\mathcal{F}^{-1}\big(\V\big)\right)\right\|_{L^\infty} &\lesssim \left\|\frac{1}{1+|\cdot|^{2}}\right\|_{W^{j,\infty}} \left\|(1+|\cdot|^j) R_0\left(\mathcal{F}^{-1}\big(\V\big)_1\right)\right\|_{L^1}\\ 
&\lesssim \left\|\mathcal{F}^{-1}\big(\V\big)\right\|_{L^2}^2 \left\|(1+|\cdot|^j)\mathcal{F}^{-1}\big(\V\big)\right\|_{L^\infty}\\
&\lesssim \big\|\V\big\|_{L^2}^2 \big\|\V\big\|_{W^{j,1}} \lesssim \big\|\V\big\|_{L^1} \big\|\V\big\|_{W^{j,1}},\\
\left\|\partial_k^j \mathcal{F} \mathcal{R}_2\left(\mathcal{F}^{-1}\big(\V\big)\right)\right\|_{L^1} &\lesssim \left\|\frac{1}{1+|\cdot|^{2}}\right\|_{H^j} \left\| (1+|\cdot|^j) R_0\left(\mathcal{F}^{-1}\big(\V\big)_1\right)\right\|_{L^2}\\ 
&\lesssim \big\|\mathcal{F}^{-1}\big(\V\big)\big\|_{L^2} \big\|\mathcal{F}^{-1}\big(\V\big)\big\|_{L^\infty} \big\|(1+|\cdot|^j) \mathcal{F}^{-1}\big(\V\big)\big\|_{L^\infty}\\ 
&\lesssim \big\|\V\big\|_{L^2} \big\|\V\big\|_{L^1} \big\|\V\big\|_{W^{j,1}}\lesssim \big\|\V\big\|_{L^1} \big\|\V\big\|_{W^{j,1}},
\end{split}
\end{align}
for $j = 0,1$ and $\V \in W^{j,1}(\R) \cap L^\infty(\R)$ with $\|\V\|_{L^1 \cap L^\infty} \leq 1$. Similarly, we obtain
\begin{align} \label{e:resbound3}
\begin{split}
\left\|\partial_k^j \mathcal{F} \mathcal{R}\left(\mathcal{F}^{-1}\big(\V\big)\right)\right\|_{L^\infty} &\lesssim \big\|\V\big\|_{L^2}^2 \big\|\V\big\|_{W^{j,1}} \big\|V\big\|_{L^1}  \lesssim \big\|\V\big\|_{L^1}^2 \big\|\V\big\|_{W^{j,1}}, 
\end{split}
\end{align}
for $j = 0,1$ and $\V \in W^{j,1}(\R) \cap L^\infty(\R)$ with $\|\V\|_{L^1 \cap L^\infty} \leq 1$. In summary, the desired bound on $\mathcal{F}\mathcal{N}(\mathcal{F}^{-1}(\V))$ follows by recalling~\eqref{e:decompN} and combining the estimates~\eqref{e:resbound1} and~\eqref{e:resbound2}, whereas the estimate on $\mathcal{E}(\V,(1-\chi)\W)$ follows by~\eqref{e:resbound0},~\eqref{e:resbound1} and~\eqref{e:resbound3}.
\end{proof}

\section{Eliminating quadratic terms} \label{sec:eliminate2}

The linear behavior of~\eqref{e:Duhcrit} at the critical frequency $k = 0$ is time-oscillatory. That is, the matrix $\smash{\re^{\widehat{\Lambda}(k)t}}$ has the eigenvalues $\re^{\pm \ri t}$ at $k = 0$. We can leverage the oscillatory character of the integrals in~\eqref{e:Duhcrit} by integrating by parts with respect to time to eliminate potentially dangerous quadratic or cubic nonlinear terms, as long as the associated phase functions are nonzero, i.e., as long as we avoid time resonances. In this section, we show that quadratic terms are not time-resonant and can therefore be eliminated.

We will establish that the phase functions $\phi_j^2 \colon \R^2 \to \R$ associated with the quadratic nonlinear term~\eqref{e:critquadr} in~\eqref{e:Duhcrit} are given by
\begin{align*}
\phi^2_{j}(k,l) = \lambda_{j_0}(k) - \lambda_{j_1}(k - l) - \lambda_{j_2}(l),
\end{align*}
with $j = (j_0,j_1,j_2) \in \Set{-1, 1}^3$. Clearly, it holds $\phi^2_j(0,0) \in \{\pm \ri\}$ for any $j \in \smash{\Set{-1,1}^3}$, implying that the phase functions do not vanish on the square $[-k_0,k_0]^2$, upon taking $k_0 > 0$ smaller if necessary. Thus, one finds that the quadratic term~\eqref{e:critquadr} is not time-resonant for frequencies close to the critical frequency $k = 0$, which allows for integration by parts with respect to time. All in all, we establish the following result.

\begin{proposition} \label{p:IBP}
The critical quadratic term~\eqref{e:critquadr} in~\eqref{e:Duhcrit} can be expressed as
\begin{align} \label{e:IBPid}
\begin{split}
&\int_0^t \chi(k) \re^{\widehat{\Lambda}(k) (t-s)} B_2\big(\Uc(s),\Uc(s)\big)(k) \, \de s\\
&\quad \ = \int_0^t \chi(k) \re^{\hL(k)(t-s)}\! \left(Q_3^1\big(\Uc(s), \Uc(s), \Uc(s)\big)(k) +  Q_3^2\big(\Uc(s), \Uc(s), \Uc(s)\big)(k) \right)\!\de s \\
&\quad \ \qquad - \, \left[e^{\hL(k)(t-s)}  Q_2\big(\Uc(s), \Uc(s)\big)(k) \right]_0^t + \int_0^t \re^{\widehat{\Lambda}(k)(t-s)} \chi(k) \mathcal{E}_2(k,s) \de s,
\end{split}
\end{align}
for $k \in \R$ and $t \in [0,T_{\max})$, where $Q_2$ is the bilinear form on $C_c(-k_0,k_0) \cap \mathcal{X}$ given by
\begin{align*}
Q_2\big(\widehat{V}, \widehat{W}\big)(k) = \sum_{j \in \Set{-1,1}^3} \int_\R \frac{\chi(k) P_{j_0}(k)}{\phi^2_j(k, l)} N_2(k,l)\big(P_{j_1} \widehat{V}, P_{j_2}\widehat{W}\big) \de l,
\end{align*}
and obeying the bound
\begin{align} \label{e:Q2bound}
\big\|\partial_k^j Q_2\big(\widehat{V}, \widehat{W}\big)\big\|_{L^p} \lesssim \big\|\V\big\|_{W^{j,1}}\big\|\W\big\|_{L^p},
\end{align}
for $\V,\W \in C_c(-k_0,k_0) \cap \mathcal{X}$, $p = 1,\infty$ and $j = 0,1$ (with $j = 0$ in case $\mathcal{X} = Y_0$). Moreover, $Q_3^1, Q_3^2$ are the trilinear forms on $C_c(-k_0,k_0) \cap \mathcal{X}$ defined by
\begin{align*}
Q_3^1\big(\widehat{V}, \widehat{W},\widehat{Z}\big)(k) &= \sum_{j \in \Set{-1,1}^3} \int_\R \frac{P_{j_0}(k)}{\phi^2_{j}(k, l)} N_2(k, l)\big(P_{j_1} \chi B_2\big(\widehat{V}, \widehat{W}\big), P_{j_2}\widehat{Z}\big) \de l\\
Q_3^2\big(\widehat{V}, \widehat{W},\widehat{Z}\big)(k) &= \sum_{j \in \Set{-1,1}^3} \int_\R \frac{P_{j_0}(k)}{\phi^2_{j}(k, l)} N_2(k, l)\big(P_{j_1} \widehat{V}, P_{j_2} \chi B_2\big(\widehat{W}, \widehat{Z}\big)\big) \de l,
\end{align*}
and satisfying the estimate
\begin{align} \label{e:Q3bound}
\big\|\partial_k^j Q_3^i\big(\widehat{V}, \widehat{W},\widehat{Z}\big)\big\|_{L^\infty} \lesssim \big\|\V\big\|_{W^{j,1}}\big\|\W\big\|_{L^\infty}\big\|\widehat{Z}\big\|_{L^1},
\end{align}
for $\V,\W,\widehat{Z} \in C_c(-k_0,k_0) \cap \mathcal{X}$, $i = 1,2$ and $j = 0,1$ (with $j = 0$ in case $\mathcal{X} = Y_0$). In addition, $Q_2\big(\widehat{V},\widehat{V}\big)(0)$ is real for $\V \in C_c(-k_0,k_0)$ satisfying $\smash{\overline{\V(k)} = \V(-k)}$ for $k \in (-k_0,k_0)$. Finally, there exists a $t$-independent constant $C>0$ such that, if we have $t \in [0,T_{\max})$ with $\|\U(t)\|_{L^1 \cap L^\infty} \leq 1$, then the remainder $\mathcal{E}_2(\cdot,t)$ can be bounded as
\begin{align*}
\left\|\partial_k^j\left(\chi \mathcal{E}_2(\cdot,t)\right)\right\|_{L^\infty} &\leq C \big\|\U(t)\big\|_{W^{j,1}} \big\|\U(t)\big\|_{L^1} \left(\big\|\Us(t)\big\|_{L^\infty} + \big\|\U(t)\big\|_{L^1}\right),
\end{align*}
for $j = 0,1$ (with $j = 0$ in case $\mathcal{X} = Y_0$).
\end{proposition}
\begin{proof}
We start by isolating the critical quadratic term in the nonlinearity in the equation~\eqref{e:vKGsysFourierCrit} for $\Uc$. Thus, suppressing the $t$-dependency of $\Us$ and $\Uc$, we write
\begin{align}
\label{e:critnonlexpansion}
\chi \mathcal{F} \mathcal{N}\big(\mathcal{F}^{-1}\big(\Uc + \Us\big)\big) = \chi B_2\big(\Uc,\Uc\big) + \chi R_2\big(\Uc,\Us\big),
\end{align}
where we denote
\begin{align*}
R_2\big(\Uc,\Us\big) &= 2B_2\big(\Uc,\Us\big) + B_2\big(\Us,\Us\big) +  \mathcal{F} \mathcal{R}_2\left(\mathcal{F}^{-1}\big(\Uc + \Us\big)\right),
\end{align*}
with $\mathcal{R}_2(V)$ defined by~\eqref{e:defR2}. Combining~\eqref{e:resbound1} and~\eqref{e:resbound2}, yields a $t$-independent constant $C > 0$ such that, if we have $t \in [0,T_{\max})$ with $\|\U(t)\|_{L^1 \cap L^\infty} \leq 1$, then it holds
\begin{align} \label{e:R2bound}
\begin{split}
\big\|\partial_k^j R_2\big(\Uc(t),\Us(t)\big)\big\|_{L^\infty} &\leq C\big\|\U(t)\big\|_{W^{j,1}}\left(\big\|\Us(t)\big\|_{L^\infty} + \big\|\U(t)\big\|_{L^1} \right),
\end{split}
\end{align}
for $j = 0,1$.

Next, we rewrite the critical quadratic term~\eqref{e:critquadr} using integration by parts with respect to time. First, we recall that the spectral projections $P_\pm(k)$ of $\hL(k)$ are well-defined  on the interval $(-k_0,k_0)$ on which $\Uc(t)$ is supported. Thus, we can decompose~\eqref{e:critquadr} into
\begin{align} \label{e:IBP1}
\begin{split}
&\int_0^t \chi(k) \re^{\hL(k) (t-s)} B_2\big(\Uc(s), \Uc(s)\big)(k) \, \de s \\ 
&\quad \ = \sum_{j \in \Set{-1,1}^3} \int_0^t  \re^{\lambda_{j_0}(k)(t-s)} \chi(k) P_{j_0}(k) \int_\R N_2(k,l)\big(P_{j_1} \Uc(s), P_{j_2}\Uc(s)\big) \de l \de s,
\end{split}
\end{align}
with $k \in \R$ and $t \in [0,T_{\max})$. We will integrate by parts in each summand of~\eqref{e:IBP1} by integrating the exponential $\smash{\re^{\lambda_{j_0}(k)(t-s)}}$ and differentiating the quadratic contribution $B_2(P_{j_1} \Uc(s), P_{j_2}\Uc(s))(k)$ with respect to $s$. Thus, expressing the derivative $\partial_s \Uc(s)$ using the equation~\eqref{e:vKGsysFourierCrit}, recalling the expansion~\eqref{e:critnonlexpansion} and noting that $\hL(k) P_{j_i}(k) = \lambda_{j_i}(k) P_{j_i}(k)$ for $i = 1,2$, we compute
\begin{align} \label{e:derivN2expr}
\begin{split}
&\partial_s N_2(k,l) \big(P_{j_1} \Uc, P_{j_2} \Uc\big)\\ 
&\quad\ = N_2(k,l) \big(P_{j_1} \partial_s \Uc, P_{j_2} \Uc\big) + N_2(k,l) \big(P_{j_1} \Uc, P_{j_2} \partial_s \Uc\big)\\
&\quad\ = \big(\lambda_{j_1}(k-l) + \lambda_{j_2}(l)\big) N_2(k,l)\big(P_{j_1}\Uc, P_{j_2}\Uc\big)\\
&\quad\ \qquad + \, N_2(k,l)\big(\chi P_{j_1} B_2\big(\Uc, \Uc\big), P_{j_2} \Uc\big)  + N_2(k,l)\big(P_{j_1} \Uc, \chi P_{j_2}  B_2\big(\Uc, \Uc\big)\big) \\
&\quad\ \qquad + \, N_2(k,l)\big(\chi P_{j_1} R_2\big(\Uc, \Us\big), P_{j_2} \Uc\big) + N_2(k,l)\big(P_{j_1} \Uc, \chi P_{j_2} R_2\big(\Uc,\Us\big)\big),
\end{split}
\end{align}
for $k, l \in \R$ and $j \in \Set{-1,1}^3$, where we have suppressed the $s$-dependency of $\Us$ and $\Uc$. Since it holds $\phi_j^2(0,0) \in \{\pm \ri\}$ for all $j\in \Set{-1,1}^3$, the phase function $\phi_j^2$ does not vanish on the square $[-k_0,k_0]^2$, upon taking $k_0 > 0$ smaller if necessary. Hence, using integration by parts, we obtain
\begin{align*}
&\int_0^t \re^{\lambda_{j_0}(k)(t-s)} \int_\R \chi(k) P_{j_0}(k) N_2(k, l)\big(P_{j_1} \Uc(s), P_{j_2}\Uc(s)\big) \de l \de s\\ 
&\quad \ = \int_0^t \re^{\lambda_{j_0}(k)(t-s)} \int_\R \chi(k) \frac{\lambda_{j_0}(k)}{\phi_j^2(k,l)} P_{j_0}(k) N_2(k, l)\big(P_{j_1} \Uc(s), P_{j_2}\Uc(s)\big) \de l \de s\\
&\quad \ \qquad - \, \int_0^t \re^{\lambda_{j_0}(k)(t-s)} \int_\R \chi(k)\frac{\lambda_{j_1}(k-l) + \lambda_{j_2}(l)}{\phi_j^2(k,l)} P_{j_0}(k) N_2(k, l)\big(P_{j_1} \Uc(s), P_{j_2}\Uc(s)\big) \de l \de s \\
&\quad \ = -\left[\re^{\lambda_{j_0}(k)(t-s)} \int_\R \frac{\chi(k) }{\phi_j^2(k,l)} P_{j_0}(k) N_2(k, l)\big(P_{j_1} \Uc(s), P_{j_2}\Uc(s)\big) \de l \right]_{s = 0}^t\\
&\quad \ \qquad + \, \int_0^t \re^{\lambda_{j_0}(k)(t-s)} \int_\R \frac{\chi(k)}{\phi_j^2(k,l)} P_{j_0}(k) \partial_s N_2(k, l)\big(P_{j_1} \Uc(s), P_{j_2}\Uc(s)\big) \de l \de s\\
&\quad \ \qquad - \, \int_0^t \re^{\lambda_{j_0}(k)(t-s)} \int_\R \chi(k) P_{j_0}(k) \frac{\lambda_{j_1}(k-l) + \lambda_{j_2}(l)}{\phi_j^2(k,l)} N_2(k, l)\big(P_{j_1} \Uc(s), P_{j_2}\Uc(s)\big) \de l \de s
\end{align*}
for $k \in \R$, $t \in [0,T_{\max})$ and $j\in \Set{-1,1}^3$. So, taking the sum of the latter over all $j \in \Set{-1,1}^3$, we obtain the desired identity~\eqref{e:IBPid} by~\eqref{e:IBP1} and~\eqref{e:derivN2expr} with error function
\begin{align*}
\mathcal{E}_2(k,t) &= \sum_{j \in \Set{-1,1}^3} \int_\R \frac{1}{\phi_j^2(k,l)} P_{j_0}(k) \left(
N_2(k,l)\big(\chi P_{j_1} R_2\big(\Uc(t), \Us(t)\big), P_{j_2} \Uc\big)\right. \\
&\qquad \left. +\, N_2(k,l)\big(P_{j_1} \Uc(t), \chi P_{j_2} R_2\big(\Uc(t),\Us(t)\big)\big)\right) \de l,
\end{align*}
for $t \in [0,T_{\max})$. Applying Young's convolution inequality, while using estimate~\eqref{e:R2bound} and recalling that the phase function $\phi_j^2$ is bounded away from $0$ on the square $[-k_0,k_0]^2$ and the cut-off function $\chi$ is supported on $(-k_0,k_0)$, we obtain $t$-independent constants $C_{1,2} > 0$ such that, if we have $t \in [0,T_{\max})$ with $\smash{\|\U(t)\|_{L^1 \cap L^\infty}} \leq 1$, then it holds
\begin{align*}
&\left\|\partial_k^j\left(\chi \mathcal{E}_2(\cdot,t)\right)\right\|_{L^\infty} \\
&\quad \ \leq C_1 \left(\big\|\Uc(t)\big\|_{L^1} \big\|R_2(\Uc(t),\Us(t)\big\|_{W^{j,\infty}} + \big\|\Uc(t)\big\|_{W^{j,1}} \big\|R_2(\Uc(t),\Us(t)\big\|_{L^\infty}\right) \\
&\quad \ \leq C_2 \big\|\U(t)\big\|_{W^{j,1}} \big\|\U(t)\big\|_{L^1} \left(\big\|\Us(t)\big\|_{L^\infty} + \big\|\U(t)\big\|_{L^1}\right),
\end{align*}
for $j = 0,1$ (with $j = 0$ in case $\mathcal{X} = Y_0$). Similarly, we derive the bounds~\eqref{e:Q2bound} and~\eqref{e:Q3bound}.

Next, let $\V \in C_c(-k_0,k_0)$ with $\overline{\V(k)} = \V(-k)$ for $k \in \R$. We prove $Q_2(\V,\V)(0) \in \R$. First, we note that the maps $\tau, \xi \colon \Set{-1,1}^3 \to \Set{-1,1}^3$ given by $\tau(j_0,j_1,j_2) = (j_0,j_2,j_1)$ and $\xi(j) = -j$ are bijections. We have
\begin{align*}
\phi_j^2(k,l) = \phi_{\tau(j)}^2(k,k-l),
\end{align*}
for $k,l \in \R$ and $j \in \Set{-1,1}^3$. Hence, using~\eqref{e:conjugateids}, we obtain
\begin{align*}
\overline{Q_2(\V,\V)(0)} &= \sum_{j \in \Set{-1,1}^3} \int_\R \frac{P_{-j_0}(0)}{\phi_{-j}^2(0,l)} N_2(0,l)\left(P_{-j_1} \overline{\V}, P_{-j_2} \overline{\V}\right) \de l\\
&= \sum_{j \in \Set{-1,1}^3} \int_\R \frac{P_{j_0}(0)}{\phi_{j}^2(0,l)} N_2(0,l)\left(P_{j_2} \V, P_{j_1} \V \right) \de l\\
&= \sum_{j \in \Set{-1,1}^3} \int_\R \frac{P_{\tau(j)_0}(0)}{\phi_{\tau(j)}^2(0,l)} N_2(0,l)\left(P_{\tau(j)_1} \V, P_{\tau(j)_2} \V \right) \de l\\
&= Q_2(\V,\V)(0),
\end{align*}
which implies $Q_2(\V,\V)(0) \in \R$.
\end{proof}

\section{Eliminating nonresonant cubic terms} \label{sec:eliminate3}

In Proposition~\ref{p:IBP} we have expressed the critical quadratic term~\eqref{e:critquadr} in the Duhamel formula~\eqref{e:Duhcrit} of $\Uc$ through integration by parts in time as a cubic integral term
\begin{align} \label{e:critcubic2}
\int_0^t \chi(k) \re^{\hL(k)(t-s)} \left(Q_3^1\big(\Uc(s), \Uc(s), \Uc(s)\big)(k) + Q_3^1\big(\Uc(s), \Uc(s), \Uc(s)\big)(k)\right) \de s,
\end{align}
a quadratic boundary term
\begin{align} \label{e:quadboundary}
\begin{split}
\left[e^{\hL(k)(t-s)} Q_2\big(\Uc(s), \Uc(s)\big)(k) \right]_0^t &= Q_2\big(\Uc(t), \Uc(t)\big)(k)\\ 
&\qquad - \, \re^{\hL(k)t} Q_2\big(\Uc(0), \Uc(0)\big)(k),
\end{split}
\end{align}
and a remainder integral term
\begin{align} \label{e:remainderduh}
\int_0^t \re^{\widehat{\Lambda}(k)(t-s)} \chi(k) \mathcal{E}_2(k,s) \de s.
\end{align}
The right-hand side of the bound on $\mathcal{E}_2(k,s)$ in Proposition~\ref{p:IBP} consists of terms which have a $\Us$-contribution, or terms which are quartic in $\Uc$. Therefore, as argued before, the remainder integral term~\eqref{e:remainderduh} can be controlled through standard $L^1$-$L^\infty$-estimates. Moreover, since the boundary term~\eqref{e:quadboundary} is quadratic in $\Uc$, it necessarily decays at a higher rate than $\Uc$ itself. So, this boundary terms is irrelevant and, thus, does not obstruct the nonlinear argument.

Therefore, the only remaining critical nonlinear terms in the Duhamel formulas~\eqref{e:Duhcrit}-\eqref{e:Duhdamp} after applying Proposition~\ref{p:IBP} are the cubic terms~\eqref{e:critcubic} and~\eqref{e:critcubic2}. As for the critical quadratic term~\eqref{e:critquadr}, we can exploit the oscillatory character in time of the integrals arising in~\eqref{e:critcubic} and~\eqref{e:critcubic2} as long as the associated phase function is nonzero. We will establish that the phase functions $\phi_j^3 \colon \R^3 \to \R$ associated with the cubic terms are given by
\begin{align*}
\phi^3_{j}(k,l_1,l_2) = \lambda_{j_0}(k) - \lambda_{j_1}(k - l_1) - \lambda_{j_2}(l_1-l_2) - \lambda_{j_3}(l_2),
\end{align*}
for $j = (j_0,j_1,j_2,j_3) \in \Set{-1, 1}^4$. We identify the set of time resonances $\mathcal{T} \subset \Set{-1,1}^4$ as
\begin{align*}
\mathcal{T} = \left\{j \in \Set{-1,1}^4 : j_0-j_1-j_2-j_3 = 0\right\},
\end{align*}
and denote its complement by $\mathcal{T}^c = \Set{-1,1}^4 \setminus \mathcal{T}$. That is, for $j \in \Set{-1,1}^4$ it holds $\phi^3_j(0,0,0) = 0$ if and only if $j \in \mathcal{T}$. Moreover, it holds $\nabla \phi^3_j(0,0,0) = 0$ for all $j \in \{-1,1\}^4$, showing that the cubic terms are all space-resonant at the critical frequency $k = 0$, cf.~\cite{GERMAIN}. Splitting into time-resonant and time-nonresonant terms, the sum of the critical cubic integral terms~\eqref{e:critcubic} and~\eqref{e:critcubic2}  decomposes as
\begin{align} \label{e:IBP2}
\begin{split}
&\int_0^t \chi(k) \re^{\widehat{\Lambda}(k) (t-s)} Z_3\big(\Uc(s),\Uc(s),\Uc(s)\big)(k) \, \de s\\ 
&\quad \ = \int_0^t \chi(k) \re^{\widehat{\Lambda}(k) (t-s)} Z_3^{\mathrm{res}}\big(\Uc(s),\Uc(s),\Uc(s)\big)(k) \, \de s\\ 
&\quad \ \qquad + \, \int_0^t \chi(k) \re^{\widehat{\Lambda}(k) (t-s)} Z_3^{\mathrm{c}}\big(\Uc(s),\Uc(s),\Uc(s)\big)(k) \, \de s,
\end{split}
\end{align}
where $Z_3, Z_3^{\mathrm{res}}, Z_3^{\mathrm{c}}$ are the trilinear forms on $C_c(-k_0,k_0) \cap \mathcal{X}$ given by
\begin{align} \label{e:defZ3}
\begin{split}
Z_3\big(\V,\W,\widehat{Z}\big)(k) &= B_3\big(\V,\W,\widehat{Z}\big)(k) + Q_3^1\big(\V,\W,\widehat{Z}\big)(k) + Q_3^2\big(\V,\W,\widehat{Z}\big)(k)\\ 
&= \int_\R \int_\R \widetilde{N}_3(k,l_1,l_2)\big(\V,\W,\widehat{Z}\big)\de l_1 \de l_2,
\end{split}
\end{align}
and
\begin{align} \label{e:defZ32}
\begin{split}
Z_3^{\mathrm{res}}\big(\V,\W,\widehat{Z}\big) &= \sum_{j \in \mathcal{T}} P_{j_0}(k) Z_3\big(P_{j_1} \V, P_{j_2} \W, P_{j_3} \widehat{Z}\big),\\
Z_3^{\mathrm{c}}\big(\V,\W,\widehat{Z}\big) &= \sum_{j \in \mathcal{T}^c} P_{j_0}(k) Z_3\big(P_{j_1} \V, P_{j_2} \W, P_{j_3} \widehat{Z}\big),
\end{split}
\end{align}
with
\begin{align*} 
\begin{split}
&\widetilde{N}_3(k,l_1,l_2)\big(\V,\W,\widehat{Z}\big) \\ 
&\quad \ = N_3(k,l_1,l_2)\big(\V,\W,\widehat{Z}\big) + \sum_{j \in \Set{-1,1}^3} \frac{P_{j_0}(k)}{\phi^2_{j}(k,l_2)} N_2(k,l_2)\big(\chi P_{j_1} N_2(\cdot,l_1-l_2)  \big(\widehat{V}, \widehat{W}\big), P_{j_2}\widehat{Z}\big)\\
&\quad \ \qquad + \, \sum_{j \in \Set{-1,1}^3} \frac{P_{j_0}(k)}{\phi^2_{j}(k,l_1)} N_2(k,l_1)\big(P_{j_1} \widehat{V}, \chi P_{j_2} N_2(\cdot,l_2)\big(\widehat{W}, \widehat{Z}\big)\big).
\end{split}
\end{align*}
Using the same strategy as in the proof of Proposition~\ref{p:IBP} we express the nonresonant part of the critical cubic term~\eqref{e:IBP2} as a cubic boundary term and a quartic remainder through integration by parts in time. We arrive at the following result.

\begin{proposition} \label{p:IBP2}
The nonresonant part of the critical cubic term~\eqref{e:IBP2} can be expressed as
\begin{align} \label{e:IBPid3}
\begin{split}
&\int_0^t \chi(k) \re^{\widehat{\Lambda}(k) (t-s)} Z_3^{\mathrm{c}}\big(\Uc(s),\Uc(s),\Uc(s)\big)(k) \, \de s \\ 
&\quad \ = \int_0^t \re^{\widehat{\Lambda}(k)(t-s)} \chi(k) \mathcal{E}_3(k,s) \de s - \left[e^{\hL(k)(t-s)} K_3\big(\Uc(s), \Uc(s),\Uc(s)\big)(k) \right]_0^t,
\end{split}
\end{align}
for $k \in \R$ and $t \in [0,T_{\max})$, where $K_3$ is the trilinear form on $C_c(-k_0,k_0) \cap \mathcal{X}$ defined by
\begin{align*}
K_3\big(\widehat{V}, \widehat{W}, \widehat{Z}\big)(k) = \sum_{j \in \mathcal{T}^c} \int_\R \int_\R \frac{\chi(k) P_{j_0}(k)}{\phi^3_j(k, l_1,l_2)} \widetilde{N}_3(k,l_1,l_2)\big(P_{j_1} \widehat{V}, P_{j_2}\widehat{W}, P_{j_3} \widehat{Z}\big) \de l_1 \de l_2,
\end{align*}
and obeying the bound
\begin{align} \label{e:K3bound}
\big\|\partial_k^j K_3\big(\V,\W,\widehat{Z}\big)\big\|_{L^p} \lesssim \big\|\V\big\|_{W^{j,1}} \big\|\W\big\|_{L^p} \big\|\widehat{Z}\big\|_{L^1}, 
\end{align}
for $\V,\W,\widehat{Z} \in C_c(-k_0,k_0) \cap \mathcal{X}$, $p = 1,\infty$ and $j = 0,1$ (with $j = 0$ in case $\mathcal{X} = Y_0$). In addition, $K_3\big(\V,\V,\V\big)(0)$ is real for $\V \in C_c(-k_0,k_0)$ satisfying $\smash{\overline{\V(k)} = \V(-k)}$ for $k \in \R$. Finally, there exists a $t$-independent constant $C>0$ such that, if we have $t \in [0,T_{\max})$ with $\|\U(t)\|_{L^1} \leq 1$, then the remainder $\mathcal{E}_3(\cdot,t)$ can be bounded as
\begin{align*} 
\left\|\partial_k^j\left(\chi \mathcal{E}_3(\cdot,t)\right)\right\|_{L^\infty} &\leq C \big\|\U(t)\big\|_{L^1}^2 \big\|\U(t)\big\|_{W^{j,1}}  \big\|\U(t)\big\|_{L^\infty}.
\end{align*}
for $j = 0,1$ (with $j = 0$ in case $\mathcal{X} = Y_0$).
\end{proposition}
\begin{proof}
We integrate by parts in each summand of
\begin{align} \label{e:IBP3}
\begin{split}
&\int_0^t \chi(k) \re^{\widehat{\Lambda}(k) (t-s)} Z_3^{\mathrm{c}}\big(\Uc(s),\Uc(s),\Uc(s)\big)(k) \, \de s\\ 
&\quad \ = 
\sum_{j \in \mathcal{T}^c} \int_0^t \chi(k) \re^{\lambda_0(k) (t-s)} P_{j_0}(k) Z_3\big(P_{j_1} \Uc(s), P_{j_2} \Uc(s), P_{j_3} \Uc(s)\big)(k) \, \de s\\
&\quad \ = 
\sum_{j \in \mathcal{T}^c} \int_0^t \chi(k) \re^{\lambda_0(k) (t-s)}P_{j_0}(k) \int_\R\int_\R \widetilde{K}_j(k,l_1,l_2,s) \de l_1 \de l_2 \de s,
\end{split}
\end{align}
where we abbreviate
\begin{align*}
\widetilde{K}_j(k,l_1,l_2,s) &= \widetilde{N}_3(k,l_1,l_2)\big(P_{j_1} \Uc(s), P_{j_2}\Uc(s),P_{j_3}\Uc(s)\big).
\end{align*}
We express the derivative $\partial_s \Uc(s)$ through equation~\eqref{e:vKGsysFourierCrit} and note that $\hL(k) P_{j_i}(k) = \lambda_{j_i}(k) P_{j_i}(k)$ for $i = 1,2,3$, which leads to
\begin{align} \label{e:derivN3expr}
\begin{split}
&\partial_s \widetilde{K}_j(k,l_1,l_2,s)\\ 
&\quad \ = \widetilde{N}_3(k,l_1,l_2) \big(P_{j_1} \partial_s \Uc, P_{j_2} \Uc,P_{j_3} \Uc\big) + \widetilde{N}_3(k,l_1,l_2) \big(P_{j_1} \Uc, P_{j_2} \partial_s \Uc,P_{j_3} \Uc\big)\\ &\quad \ \qquad + \, \widetilde{N}_3(k,l_1,l_2) \big(P_{j_1} \Uc, P_{j_2} \Uc,P_{j_3} \partial_s \Uc\big) \\
&\quad \ = \big(\lambda_{j_1}(k-l_1) + \lambda_{j_2}(l_1-l_2) + \lambda_{j_3}(l_2)\big) \widetilde{K}_j(k,l_1,l_2,s) +  \check{K}_j(k,l_1,l_2,s),
\end{split}
\end{align}
for $k, l_1,l_2 \in \R$ and $j \in \Set{-1,1}^4$, where we we have suppressed the $s$-dependency of $\Us$ and $\Uc$ and we
denote
\begin{align*}
\check{K}_j(k,l_1,l_2,s) &= 
\widetilde{N}_3(k,l_1,l_2)\big(\chi P_{j_1} \mathcal{F} \mathcal{N} \big(\mathcal{F}^{-1}\big(\Uc(s) + \Us(s)\big)\big), P_{j_2} \Uc(s), P_{j_3} \Uc(s) \big)\\ 
&\qquad + \, \widetilde{N}_3(k,l_1,l_2)\big(P_{j_1} \Uc(s), \chi P_{j_2}  \mathcal{F} \mathcal{N} \big(\mathcal{F}^{-1}\big( \Uc(s) + \Us(s)\big)\big), P_{j_3} \Uc(s) \big) \\
&\qquad + \, \widetilde{N}_3(k,l_1,l_2)\big(P_{j_1} \Uc(s), P_{j_2} \Uc(s), \chi P_{j_3}  \mathcal{F} \mathcal{N} \big(\mathcal{F}^{-1}\big( \Uc(s) + \Us(s)\big)\big) \big),
\end{align*}
Since we have $\phi_j^3(0,0,0) \neq 0$ for all $j\in \mathcal{T}^c$, the phase function $\phi_j^3$ does not vanish on the cube $[-k_0,k_0]^3$, upon taking $k_0 > 0$ smaller if necessary. Hence, through integration by parts we obtain
\begin{align*}
&\int_0^t \re^{\lambda_{j_0}(k)(t-s)} \int_\R \int_\R \chi(k) P_{j_0}(k) \widetilde{K}_j(k,l_1,l_2,s) \de l_1 \de l_2 \de s\\ 
&\quad \ = \int_0^t \re^{\lambda_{j_0}(k)(t-s)} \int_\R \int_\R \chi(k) \frac{\lambda_{j_0}(k)}{\phi_j^3(k,l_1,l_2)} P_{j_0}(k) \widetilde{K}_j(k,l_1,l_2,s) \de l_1 \de l_2  \de s\\
&\quad \ \qquad + \, \int_0^t \re^{\lambda_{j_0}(k)(t-s)} \int_\R \int_\R \chi(k)\frac{\phi_j^3(k,l_1,l_2) - \lambda_{j_0}(k)}{\phi_j^3(k,l_1,l_2)} P_{j_0}(k) \widetilde{K}_j(k,l_1,l_2,s) \de l_1 \de l_2 \de s \\
&\quad \ = -\left[\re^{\lambda_{j_0}(k)(t-s)} \int_\R \int_\R \frac{\chi(k) }{\phi_j^3(k,l_1,l_2)} P_{j_0}(k) \widetilde{K}_j(k,l_1,l_2,s) \de l_1 \de l_2 \right]_{s = 0}^t\\
&\quad \ \qquad + \, \int_0^t \re^{\lambda_{j_0}(k)(t-s)} \int_\R \int_\R \frac{\chi(k)}{\phi_j^3(k,l_1,l_2)} P_{j_0}(k) \partial_s \widetilde{K}_j(k,l_1,l_2,s) \de l \de s\\
&\quad \qquad - \, \int_0^t \re^{\lambda_{j_0}(k)(t-s)} \int_\R \int_\R \chi(k) P_{j_0}(k) \frac{\lambda_{j_1}(k-l_1) + \lambda_{j_2}(l_1-l_2) + \lambda_{j_3}(l_2)}{\phi_j^3(k,l_1,l_2)}\\
&\quad \ \qquad \cdot \, \widetilde{K}_j(k,l_1,l_2,s) \de l_1 \de l_2 \de s
\end{align*}
for $k \in \R$, $t \in [0,T_{\max})$ and $j\in \mathcal{T}^c$. So, taking the sum of the latter over all $j \in \mathcal{T}^c$ while recalling~\eqref{e:derivN3expr}, we arrive at~\eqref{e:IBPid3} with error function
\begin{align*}
\mathcal{E}_3(k,t) &= \sum_{j \in \mathcal{T}^c} \int_\R \int_\R \frac{1}{\phi_j^3(k,l_1,l_2)} P_{j_0}(k) \check{K}_j(k,l_1,l_2,s) \de l_1 \de l_2,
\end{align*}
for $t \in [0,T_{\max})$. By Young's convolution inequality, Lemma~\ref{l:noncriticalbounds} and the facts that the phase function $\phi_j^3$ is bounded away from $0$ on the cube $[-k_0,k_0]^3$ for $j \in \mathcal{T}^c$ and the cut-off function $\chi$ is supported on $(-k_0,k_0)$, we find $t$-independent constants $C_{1,2} > 0$ such that, if we take $t \in [0,T_{\max})$ with $\smash{\|\U(t)\|_{L^1}} \leq 1$, then we have
\begin{align*}
\left\|\partial_k^j\left(\chi \mathcal{E}_3(\cdot,t)\right)\right\|_{L^\infty} &\leq C_1 \big\|\Uc(t)\big\|_{L^1} \left(\big\|\mathcal{F} \mathcal{N}\big(\mathcal{F}^{-1} \big(\Uc(t) + \Us(t)\big)\big)\big\|_{W^{j,\infty}} \big\|\Uc(t)\big\|_{L^1} \right. 
\\
&\left.\qquad + \, \big\|\Uc(t)\big\|_{W^{j,1}} \big\|\mathcal{F} \mathcal{N}\big(\mathcal{F}^{-1} \big(\Uc(t) + \Us(t)\big)\big)\big\|_{L^\infty}\right)\\
&\leq C_2 \big\|\U(t)\big\|_{L^1}^2 \big\|\U(t)\big\|_{W^{j,1}} \big\|\U(t)\big\|_{L^\infty},
\end{align*}
for $j = 0,1$ (with $j = 0$ in case $\mathcal{X} = Y_0$). Similar considerations afford the bound~\eqref{e:K3bound}. 

Finally, let $\V, \widehat{Y}_{1,2},\W_{1,2},\widehat{Z}_{1,2} \in C_c(-k_0,k_0)$ satisfy
\begin{align*}
\overline{\V(k)} = \V(-k), \qquad \overline{\widehat{Y}_1(k)} = \widehat{Y}_2(-k), \qquad \overline{\W_1(k)} = \W_2(-k), \qquad, \overline{\widehat{Z}_1(k)} = \widehat{Z}_2(-k),
\end{align*}
for $k \in (-k_0,k_0)$. We prove $K_3(\V,\V,\V)(0) \in \R$. First, we note that $\xi \colon \{-1,1\}^3 \to \{-1,1\}^3$ and $\zeta \colon \mathcal{T}^c \to \mathcal{T}^c$ given by $\xi(j) = -j$ and $\zeta(j) = -j$ are bijections. So, using~\eqref{e:conjugateids} we arrive at
\begin{align*}
&\overline{\widetilde{N}_3(0,l_1,l_2)\big(\widehat{Y}_1,\W_1,\widehat{Z}_1\big)}\\ 
&\quad \ = N_3(0,-l_1,-l_2)\big(\widehat{Y}_2,\W_2,\widehat{Z}_2\big)\\ 
&\quad \ \qquad + \sum_{j \in \Set{-1,1}^3} \frac{P_{-j_0}(0)}{\phi^2_{-j}(0,-l_2)} N_2(0,-l_2)\big(\chi P_{-j_1} N_2(\cdot,-l_1+l_2)  \big(\widehat{Y}_2, \widehat{W}_2\big), P_{-j_2}\widehat{Z}_2\big)\\
&\quad \ \qquad + \sum_{j \in \Set{-1,1}^3} \frac{P_{-j_0}(0)}{\phi^2_{-j}(0,-l_1)} N_2(0,-l_1)\big(P_{-j_1} \widehat{Y}_2, \chi P_{-j_2} N_2(\cdot,-l_2)\big(\widehat{W}_2, \widehat{Z}_2\big)\big)\\
&\quad \ = \widetilde{N}_3(0,-l_1,-l_2)\big(\widehat{Y}_2,\W_2,\widehat{Z}_2\big)
\end{align*}
for $l_1,l_2 \in \R$. Therefore, using~\eqref{e:conjugateids} again and applying the substitution rule, we infer
\begin{align*}
\overline{K_3\big(\widehat{V}, \widehat{V}, \widehat{V}\big)(0)} &= \sum_{j \in \mathcal{T}^c} \int_\R \int_\R \frac{P_{-j_0}(0)}{\phi^3_{-j}(0,-l_1,-l_2)} \widetilde{N}_3(0,-l_1,-l_2)\big(P_{-j_1} \widehat{V}, P_{-j_2} \widehat{V}, P_{-j_3} \widehat{V}\big) \de l_1 \de l_2\\
&= \sum_{j \in \mathcal{T}^c} \int_\R \int_\R \frac{P_{-j_0}(0)}{\phi^3_{-j}(0,l_1,l_2)} \widetilde{N}_3(0,l_1,l_2)\big(P_{-j_1} \widehat{V}, P_{-j_2} \widehat{V}, P_{-j_3} \widehat{V}\big) \de l_1 \de l_2\\
&= K_3\big(\widehat{V}, \widehat{V}, \widehat{V}\big)(0)
\end{align*}
implying $K_3\big(\widehat{V}, \widehat{V}, \widehat{V}\big)(0) \in \R$.
\end{proof}

\section{Resonant cubic terms} \label{sec:rescubic}

In the previous section, we first applied Proposition~\ref{p:IBP} to express the quadratic critical term~\eqref{e:critquadr} in the Duhamel formula~\eqref{e:Duhcrit} for $\Uc$ as~\eqref{e:IBPid} and, subsequently, decomposed the remaining critical cubic terms~\eqref{e:critcubic} and~\eqref{e:critcubic2} in~\eqref{e:Duhcrit} in a resonant and nonresonant part in~\eqref{e:IBP2}. Then, we expressed in Proposition~\ref{p:IBP2} the nonresonant critical cubic terms as~\eqref{e:IBP3}. After these manipulations~\eqref{e:Duhcrit} reads
\begin{align} \label{e:DuhZc}
\begin{split}
\Zc(k,t) &= \re^{\hL(k) t} \Zc(k,0) + \int_0^t \re^{\hL(k)(t-s)} \chi(k) Z_3^{\mathrm{res}}\big(\Zc(s),\Zc(s),\Zc(s)\big)(k) \, \de s\\ &\qquad + \, \int_0^t \re^{\hL(k)(t-s)} \chi(k) \mathcal{E}_4(k,s) \de s,
\end{split}
\end{align}
with $k \in \R$ and $t \in [0,T_{\max})$, where we denote 
\begin{align} \label{e:normalform}
\Zc(k,t) = \Uc(k,t) + Q_2\big(\Uc(t),\Uc(t)\big)(k) + K_3\big(\Uc(t),\Uc(t),\Uc(t)\big)(k),
\end{align}
and
\begin{align*}
\mathcal{E}_4(k,s) &\deq \mathcal{E}_2(k,s) + \mathcal{E}_3(k,s) + \mathcal{E}\big(\Uc(s),\Us(s)\big)(k)\\
&\qquad + \, Z_3^{\mathrm{res}}\big(\Uc(s),\Uc(s),\Uc(s)\big)(k) - Z_3^{\mathrm{res}}\big(\Zc(s),\Zc(s),\Zc(s)\big)(k).
\end{align*}
The estimates~\eqref{e:Q2bound} and~\eqref{e:K3bound} on the bi- and trilinear forms $Q_2$ and $K_3$ in Propositions~\ref{p:IBP} and~\ref{p:IBP2}, respectively, in combination with~\eqref{e:regularity3} readily yield 
\begin{align} \label{e:regularity4}
\Zc \in C\big([0,T_{\max}),W^{1,1}(\R) \cap W^{1,\infty}(\R)\big),
\end{align}
in case $\mathcal{X} = X_0$, and
\begin{align} \label{e:regularity42}
\Zc \in C\big([0,T_{\max}),LW^1{1,1}(\R) \cap L^\infty W^{1,\infty}(\R)\big),
\end{align}
in case $\mathcal{X} = Y_0$. Moreover, $\Zc(\cdot,t)$ is supported on $(-k_0,k_0)$ for each $t \in [0,T_{\max})$, since $\chi$ is. Furthermore, the estimates~\eqref{e:Q2bound} and~\eqref{e:K3bound} indicate that the new variable $\Zc$ exhibits the same decay rate as $\Uc$. In fact, the coordinate change~\eqref{e:normalform} can be regarded as a normal form transform for the equation~\eqref{e:vKGsysFourierCrit} of $\Uc$, which eliminates nonresonant critical nonlinearities. Applying the operator $\partial_t - \hL(k)$ to~\eqref{e:DuhZc} we arrive at the equation
\begin{align} \label{e:Zceq}
\partial_t \Zc = \hL \Zc + \chi  Z_3^{\mathrm{res}}\big(\Zc,\Zc,\Zc\big) + \chi \mathcal{E}_4(\cdot,t).
\end{align}
We bound the critical cubic nonlinearity $Z_3^{\mathrm{res}}$ in~\eqref{e:Zceq}.

\begin{lemma} \label{l:Zcnonlinearest}
We have
\begin{align} \label{e:Z_3est}
\big\|\partial_k^j Z_3^{\mathrm{res}}\big(\V,\W,\widehat{Z}\big)\big\|_{L^\infty} \lesssim \big\|\V\big\|_{W^{j,1}}\big\|\W\big\|_{L^\infty}\big\|\widehat{Z}\big\|_{L^1},
\end{align}
for $\V,\W, \widehat{Z} \in C_c(-k_0,k_0) \cap \mathcal{X}$ and $j = 0,1$ (with $j = 0$ in case $\mathcal{X} = Y_0$).
\end{lemma}
\begin{proof}
Young's convolution inequality yields
\begin{align*}
\big\|\partial_k^j B_3\big(\V,\W,\widehat{Z}\big)\big\|_{L^\infty} \lesssim \big\|\V\big\|_{W^{j,1}}\big\|\W\big\|_{L^\infty}\big\|\widehat{Z}\big\|_{L^1},
\end{align*}
for $j = 0,1$ and $\V \in W^{j,1}(\R), \W \in L^\infty(\R)$ and $\widehat{Z} \in L^1(\R)$. Combining the latter with the estimates on $Q_3^1$ and $Q_3^2$ established in Proposition~\ref{p:IBP} we infer~\eqref{e:Z_3est}.
\end{proof}

Next, we obtain a bound on the residual nonlinearity $\mathcal{E}_4(\cdot,t)$ in~\eqref{e:Zceq}. 

\begin{lemma} \label{l:residualnonlZc}
There exists a $t$- and $k$-independent constant $C > 0$ such that, if $k \in \R$ and $t \in [0,T_{\max})$ with $\|\U(t)\|_{L^1 \cap L^\infty} \leq 1$, then we have the bound
\begin{align*}
\left|\partial_k^j\left(\chi(k) \mathcal{E}_4(k,t)\right)\right| &\leq C\left(\big\|\U(t)\big\|_{W^{j,1}} \big\|\U(t)\big\|_{L^1} \left(\big\|\U(t)\big\|_{L^1} + \big\|\Us(t)\big\|_{L^\infty}\right) \right.\\
&\left.\qquad + \, \big\|\Us(t)\big\|_{W^{j,1}}\big\|\Us(t)\big\|_{L^\infty} + |\chi(k) k| \big\|\Us(t)\big\|_{L^\infty} \big\|\U(t)\big\|_{W^{j,1}} \right. \\
&\left.\qquad + \, \big\|\Us(t)\big\|_{L^\infty} \big\||\cdot|\Uc(t)\big\|_{W^{j,1}} + \big\|\Us(t)\big\|_{L^1} \big\|\U(t)\big\|_{W^{j,1}}\right),
\end{align*}
for $j = 0,1$ (with $j = 0$ in case $\mathcal{X} = Y_0$).
\end{lemma}
\begin{proof}
We recall that $Z_3^{\mathrm{res}}$ is a trilinear form. So, the result follows by combining Lemma~\ref{l:Zcnonlinearest} with the estimates on $\mathcal{E}, Q_2, \mathcal{E}_2, K_3$ and $\mathcal{E}_3$ established in Lemma~\ref{l:noncriticalbounds} and Propositions~\ref{p:IBP} and~\ref{p:IBP2}.
\end{proof}

\section{Reduced equations}
\label{s:reduced-eq}

In this section, we take $\mathcal{X} = X_0$. Our aim is to capture the leading-order temporal dynamics of $\Zc(k,t)$ by proceeding as in~\cite[Theorem~1.5]{DRS} and decomposing $\Zc(k,t)$ into an explicit leading-order Gaussian part and a remainder vanishing at the critical frequency $k = 0$. In order to control the leading-order Gaussian part, we derive an ordinary differential equation for $A(t) \,\smash{\deq \re^{-\hL(0)(t+1)} \Zc(0,t)}$. If the sign condition~\eqref{e:signcondition} is fulfilled, the cubic term in this ODE is of absorption type, which induces enhanced diffusive decay of the Gaussian part. Moreover, we derive an equation for the remainder and obtain estimates on the nonlinearities by exploiting that the remainder vanishes at frequency $k = 0$. These findings turn out to be sufficient to close a global nonlinear iteration argument in the next section, yielding the proof of Theorem~\ref{thm:mainresult}.

Thus, we introduce the new variables
\begin{align} \label{e:defAsic}
A(t) \deq \re^{-\hL(0)(t+1)} \Zc(0,t), \qquad \sic(k,t) \deq \re^{\hL(k) (t+1)} \chi(k) A(t),
\end{align}
and decompose
\begin{align} \label{e:decompZc}
\Zc(k,t) = \sic(k,t) + \rc(k,t),
\end{align}
where the residual
\begin{align*}
\rc(k,t) \deq \re^{\hL(k) (t+1)} \chi(k) A(t) - \Zc(k,t),
\end{align*}
vanishes at $k = 0$ for each $t \in [0,T_{\max})$. We note that $A(t)$ is real for each $t \in [0,T_{\max})$ by~\eqref{e:conjUC},~\eqref{e:normalform} and the fact that
\begin{align} \label{e:HL0expr}
\re^{-\hL(0) (t+1)} = \begin{pmatrix} \cos(t+1) & -\sin(t+1)\\ \sin(t+1) & \cos(t+1)\end{pmatrix},
\end{align}
$Q_2(\Uc(t),\Uc(t))(0)$ and $K_2(\Uc(t),\Uc(t),\Uc(t))(0)$ are real by Propositions~\ref{p:IBP} and~\ref{p:IBP2}. So, by~\eqref{e:regularity4}, smoothness of $\hL(\cdot)$ and $\chi$, and the fact that $\chi$ has compact support, we infer
\begin{align} \label{e:regularity5}
A \in C\big([0,T_{\max}),\R\big), \qquad \rc,\sic \in  C\big([0,T_{\max}),W^{1,1}(\R)\cap W^{1,\infty}(\R)\big).
\end{align}

Substituting $\Zc(k,t) = \re^{\hL(k) (t+1)} \chi(k) A(t) + \rc(k,t)$ into~\eqref{e:Zceq} we arrive at
\begin{align} \label{e:Zceq2}
\begin{split}
\re^{\hL(k) (t+1)} \chi(k) A'(t) + \partial_t \rc(k,t) &= \hL(k) \rc(k,t) + \chi(k) Z_3^{\mathrm{res}}\big(\Zc(t),\Zc(t),\Zc(t)\big)(k)\\ \qquad + \, \chi(k) \mathcal{E}_4(k,t),
\end{split}
\end{align}
for $k \in \R$ and $t \in [0,T_{\max})$. Setting $k = 0$ in~\eqref{e:Zceq2}, we find
\begin{align} \label{e:Aeq}
\re^{\hL(0) (t+1)} A'(t) = Z_3^{\mathrm{res}}\big(\Zc(t),\Zc(t),\Zc(t)\big)(0) +  \mathcal{E}_4(0,t),
\end{align}
for $t \in [0,T_{\max})$. Subsequently, inserting~\eqref{e:Aeq} into~\eqref{e:Zceq2}, we obtain
\begin{align} \label{e:rceq}
\partial_t \rc(k,t) = \hL(k) \rc(k,t) + \chi(k) \left(\mathcal{E}_4(k,t) - \re^{\hL(k) (t+1)} \re^{-\hL(0) (t+1)} \mathcal{E}_4(0,t) + \mathcal{E}_5(k,t)\right),
\end{align}
where we denote
\begin{align*}
\mathcal{E}_5(k,t) &= Z_3^{\mathrm{res}}\big(\Zc(t),\Zc(t),\Zc(t)\big)(k) - \re^{\hL(k) (t+1)} \re^{-\hL(0) (t+1)} Z_3^{\mathrm{res}}\big(\Zc(t),\Zc(t),\Zc(t)\big)(0),
\end{align*}
for $t \in [0,T_{\max})$ and $k \in \R$. 

In our nonlinear argument we control the dynamics of the residual $\rc(t)$ through iterative estimates on the Duhamel formulation associated with~\eqref{e:rceq}. The nonlinear term $\mathcal{E}_4(\cdot,t)$ in~\eqref{e:rceq} can be bounded with the aid of Lemma~\ref{l:residualnonlZc}. Moreover, the following result provides a bound on the term $\mathcal{E}_5(\cdot,t)$ in~\eqref{e:rceq}. 

\begin{lemma} \label{l:cubicresidualZc}
There exists a $t$-independent constant $C > 0$ such that
\begin{align*} 
\left\|\partial_k^j \mathcal{E}_5(\cdot,t)\right\|_{L^\infty} &\leq C \left(\frac{|A(t)|}{\sqrt{1+t}} + \big\|\rc(t)\big\|_{L^1}\right)  \left(|A(t)| + \big\|\rc(t)\big\|_{L^\infty}\right)\\ 
&\qquad \cdot \left((1+t)^{-\frac12(1-j)} |A(t)| + \big\|\rc(t)\big\|_{W^{j,1}} + (1+t)^{\frac{j}{2}} \big\|\rc(t)\big\|_{L^1}\right),
\end{align*}
for $j = 0,1$ and $t \in [0,T_{\max})$.
\end{lemma}
\begin{proof}
With the aid of Lemma~\ref{l:linear estimates} we infer
\begin{align} \label{e:Gaussest}
\left\|\partial_k^j \left(\chi \re^{\hL(\cdot) (t+1)}\right)\right\|_{L^p} \lesssim \left\|(1 + j|\cdot|(t+1) ) \re^{-\frac12 \alpha  |\cdot|^2 (t+1)}\right\|_{L^p} \lesssim (1+t)^{-\frac{1}{2p} + \frac{j}{2}},
\end{align}
for $t \geq 0$, $p = 1,\infty$ and $j = 0,1$. Therefore, we find a $t$-independent constant $C_1 > 0$ such that
\begin{align} \label{e:Zcest0}
\big\|\partial_k^j \Zc(t)\big\|_{L^p} \leq C_1 \left(\big\|\partial_k^j \rc(t)\big\|_{L^p} + (1+t)^{-\frac{1}{2p} + \frac{j}{2}} |A(t)| \right),
\end{align}
for $j = 0,1$, $p = 1,\infty$ and $t \in [0,T_{\max})$. On the other hand, combining~\eqref{e:Gaussest} with identity~\eqref{e:HL0expr} and Lemma~\ref{l:Zcnonlinearest} we obtain a $t$-independent constant $C_2 > 0$ such that
\begin{align*} 
\left\|\partial_k^j \mathcal{E}_5(\cdot,t)\right\|_{L^\infty} &\leq C_2\big\|\Zc(t)\big\|_{L^1} \big\|\Zc(t)\big\|_{L^\infty} \left(\big\|\Zc(t)\big\|_{W^{j,1}} + (1+t)^{\frac{j}2}  \big\|\Zc(t)\big\|_{L^1}\right),
\end{align*}
for $j = 0,1$ and $t \in [0,T_{\max})$, which proves the result by inserting~\eqref{e:Zcest0} into the latter. 
\end{proof}

The fact that the residual $\rc(k,t)$ vanishes at the critical frequency $k = 0$ suggests that $\rc(k,t)$ decays at a higher rate than $A(t)$. Indeed, neglecting the nonlinear terms in the equation~\eqref{e:rceq} for $\smash{\rc(k,t)}$, we find that $\smash{\rc(k,t) = \re^{\hL(k) t} \rc(k,0) \approx k \re^{\hL(k) t} \partial_k \rc(0,0)}$ decays at rate $\smash{t^{-\frac12}}$ in $L^\infty(\R)$ and at rate $t^{-1}$ in $L^1(\R)$ by Lemma~\ref{l:linear estimates}. As argued in~\S\ref{sec:rescubic}, the nonlinear terms with an $\mathcal{E}_4(\cdot,t)$-contribution in~\eqref{e:rceq} are controllable using standard $L^1$-$L^\infty$-estimates in combination with the bounds from Lemma~\ref{l:residualnonlZc}. 

Now, let's look at the remaining nonlinear term
\begin{align} \label{e:duh_toy}
\int_0^t \re^{\hL(k)(t-s)} \chi(k) \mathcal{E}_5(k,s) \de s,
\end{align}
arising in the Duhamel formula of $\rc(k,t)$. Applying the mean value theorem and using $\mathcal{E}_5(0,s) = 0$, we can bound
\begin{align} \label{e:MVT}
\left|\mathcal{E}_5(k,s)\right| \leq |k| \left\|\partial_k \mathcal{E}_5(\cdot,s)\right\|_{L^\infty}.
\end{align}
So, we have the choice to bound $\mathcal{E}_5(k,s)$ in~\eqref{e:duh_toy} using the estimate in Lemma~\ref{l:cubicresidualZc} with $j = 0$ or using~\eqref{e:MVT} and the estimate in Lemma~\ref{l:cubicresidualZc} with $j = 1$. By the previous considerations, the terms occurring on the right-hand side of these bounds are all integrable, i.e.~decaying at rate $(1+s)^{-a}$ with $a > 1$ and, thus, controllable with standard $L^1$-$L^\infty$-estimates, except for the marginal terms
\begin{align*}
\frac{|A(s)|^3}{1+s}, \qquad \frac{|k| |A(s)|^3}{\sqrt{1+s}},
\end{align*}
arising on the right-hand side of the bounds on $\|\mathcal{E}_5(\cdot,s)\|_{L^\infty}$ and $|k| \left\|\partial_k \mathcal{E}_5(\cdot,s)\right\|_{L^\infty}$, respectively, cf.~Lemma~\ref{l:cubicresidualZc}. Since we have the freedom to choose between these bounds for each $s \in [0,t]$, we can, using Lemma~\ref{l:linear estimates}, estimate the marginal terms in~\eqref{e:duh_toy} in $L^1(\R)$ as
\begin{align*}
\int_0^{\frac{t}2} \frac{|A(s)|^3}{(t-s)\sqrt{1+s}} \de s + \int_{\frac{t}2}^t \frac{|A(s)|^3}{\sqrt{t-s}(1+s)} \de s,
\end{align*}
and in $L^\infty(\R)$ as
\begin{align*}
\int_0^t \frac{|A(s)|^3}{\sqrt{t-s}\sqrt{1+s}} \de s,
\end{align*}
which suffices to close the nonlinear argument, as we will see in the upcoming section.  

Next, we shift focus to the ordinary differential equation~\eqref{e:Aeq} for $A(t)$. Applying the projections $P_\pm(0)$ to~\eqref{e:Aeq}, we arrive at
\begin{align}  \label{e:Aeqpm}
\partial_t \left(P_\pm(0) A(t)\right) = \re^{\mp \ri (t+1)} P_\pm(0) Z_3^{\mathrm{res}}\big(\Zc(t),\Zc(t),\Zc(t)\big)(0) + \re^{\mp \ri (t+1)} P_\pm(0) \mathcal{E}_4(0,t),
\end{align}
for $t \in [0,T_{\max})$. As argued above, the nonlinear term $\re^{\mp \ri (t+1)} P_\pm(0)\mathcal{E}_4(0,t)$ is irrelevant in the sense that it can be controlled with standard iterative $L^1$-$L^\infty$-estimates by invoking Lemma~\ref{l:residualnonlZc}. The following result provides an approximation of the leading-order nonlinear term in~\eqref{e:Aeqpm}, which follows by expanding $Z_3^{\mathrm{res}}(\Zc(t),\Zc(t),\Zc(t))(0)$ using the decomposition~\eqref{e:decompZc}.

\begin{lemma} \label{l:leadingorderapproxZc}
There exists a $t$-independent constant $C > 0$ such that
\begin{align*}
\mathcal{E}_6(t) = P_{\pm}(0) Z_3^{\mathrm{res}}\big(\Zc(t),\Zc(t),\Zc(t)\big)(0) - \frac{\omega_\pm \re^{\pm \ri(t+1)}}{1+t} \left|\left(P_\pm(0) A(t)\right)_1\right|^2 \left(P_\pm(0) A(t)\right)_1 \begin{pmatrix} 1 \\ \ri\end{pmatrix},
\end{align*}
can be bounded by
\begin{align*}
\left|\mathcal{E}_6(t)\right| \leq C \left(\big\|\rc(t)\big\|_{L^\infty}\left( \big\|\rc(t)\big\|_{L^1} +  \frac{|A(t)|}{\sqrt{1+t}} \right)^2 + \frac{|A(t)|^3}{(1+t)^{\frac32}} \right),
\end{align*}
for $t \in [0,T_{\max})$, where the coefficients $\omega_\pm \in \C$ are given by
\begin{align*}
\omega_\pm = \frac{\mp \pi \ri \left(9\beta + 10\kappa^2\right)}{3 \sqrt{3\alpha^2 + 1 \mp 2 \ri \alpha}}.
\end{align*} 
\end{lemma}
\begin{proof}
First, we note that by Lemma~\ref{l:linear estimates} there exists a $t$-independent constant $C_1 > 0$ such that 
\begin{align} \label{e:sicest2}
\|\sic(t)\|_{L^1} \leq |A(t)| \left\|\re^{\hL(\cdot) (t+1)}\right\|_{L^1} \leq C_1 \frac{|A(t)|}{\sqrt{1+t}},
\end{align}
for $t \in [0,T_{\max})$. Next, we expand
\begin{align*}
\chi(k) Z_3^{\mathrm{res}}(\Zc(t),\Zc(t),\Zc(t))(k) &= \chi(k) Z_3^{\mathrm{res}}\big(\sic(t),\sic(t),\sic(t)\big)(k) + E_1(k,t),
\end{align*}
with remainder
\begin{align*}
E_1(k,t) = \chi(k) \left(Z_3^{\mathrm{res}}\big(\Zc(t),\Zc(t),\Zc(t)\big)(k) - Z_3^{\mathrm{res}}\big(\sic(t),\sic(t),\sic(t)\big)(k)\right),
\end{align*}
for $k \in \R$ and $t \in [0,T_{\max})$. The fact that $Z^3_{\mathrm{res}}$ is trilinear in combination with Lemma~\ref{l:Zcnonlinearest}, identity~\eqref{e:decompZc} and estimate~\eqref{e:sicest2} yield $t$-independent constants $C_1,C_2 > 0$ such that
\begin{align} \label{e:E1bound}
\begin{split}
\left\|E_1(\cdot,t)\right\|_{L^\infty} &\leq C_1\big\|\rc(t)\big\|_{L^\infty}\left( \big\|\rc(t)\big\|_{L^1} + \big\|\sic(t)\big\|_{L^1} \right)^2\\ 
&\leq C_2 \big\|\rc(t)\big\|_{L^\infty}\left( \big\|\rc(t)\big\|_{L^1} +  \frac{|A(t)|}{\sqrt{1+t}} \right)^2, 
\end{split}
\end{align}
for $t \in [0,T_{\max})$. 

We define truncations $\tl_\pm(k)$ of the eigenvalues $\lambda_\pm(k)$ by
\begin{align} \label{e:def_truncation}
\tl_\pm(k) = \pm \ri + \frac{1}{2}\left(-\alpha \pm \ri\right) k^2,
\end{align}
so that the expansion~\eqref{e:eigenvalue expansion} yields
\begin{align} \label{e:eigenvalueapprox}
\left|\lambda_\pm(k) - \tl_\pm(k)\right| \lesssim |k|^4,
\end{align}
for $k \in (-k_0,k_0)$. Clearly, it holds
\begin{align} \label{e:highfrequencyeigenvalues}
\Re(\tl_\pm(k)) \leq  -\frac14\alpha\left(k_0^2 + k^2\right),
\end{align}
for $k \in \R \setminus (-k_0,k_0)$. Moreover, the mean value theorem implies
\begin{align} \label{e:standardexpest}
\left|\re^z - 1\right| \leq \re^{|z|} - 1 \leq |z| \re^{|z|},
\end{align}
for $z \in \C$. Thus, combining estimates~\eqref{e:eigenvalueapprox},~\eqref{e:highfrequencyeigenvalues}, ~\eqref{e:standardexpest} and taking $k_0 > 0$ smaller if necessary, we find $t$- and $k$-independent constants $C_{1,2} > 0$ such that
\begin{align*}
&\left|P_{\pm}(k) \sic(k,t) - P_{\pm}(0) \re^{\tl_\pm(k) (t+1)} A(t)\right|\\
&\quad \ \leq \left(\left|P_\pm(k)\chi(k)\right| \left|\re^{(\lambda_\pm(k) - \tl_\pm(k))(t+1)} - 1 \right| + \left|P_\pm(k)\chi(k) - P_\pm(0)\right|\right) \re^{\Re(\tl_\pm(k)) (t+1)} |A(t)| \\
&\quad \ \leq C_1\left(|k|^4 (t+1) + |k|\right) \re^{-\frac38 \alpha k^2 (t+1)} |A(t)| \leq C_2 |k| \re^{-\frac14 \alpha k^2 (t+1)} |A(t)|,
\end{align*}
for $k \in (-k_0,k_0)$ and $t \in [0,T_{\max})$, and
\begin{align*}
\left|P_{\pm}(k) \sic(k,t) - P_{\pm}(0) \re^{\tl_\pm(k) (t+1)} A(t)\right| &= \left|P_{\pm}(0) \re^{\tl_\pm(k) (t+1)} A(t)\right|\\ &\leq C_1 \re^{-\frac14 \alpha k^2(t+1)} \re^{-\frac14 \alpha k_0^2(t+1)} |A(t)|,
\end{align*}
for $k \in \R \setminus (-k_0,k_0)$ and $t \in [0,T_{\max})$. Therefore, we obtain a $t$-independent constant $C > 0$ such that
\begin{align} \label{e:sicapprox}
\begin{split}
\left\|P_{\pm} \sic(t) - P_{\pm}(0) \re^{\tl_\pm(\cdot) (t+1)} A(t)\right\|_{L^p} &\leq \frac{C}{(1+t)^{\frac{1}{2} + \frac{1}{2p}}} |A(t)|, \\ 
\left\|P_\pm(0) \re^{\tl_\pm(\cdot) (t+1)} A(t)\right\|_{L^p} &\leq \frac{C}{(1+t)^{\frac{1}{2p}}} |A(t)|,
\end{split}
\end{align}
for $t \in [0,T_{\max})$ and $p = 1,\infty$. 

Recalling~\eqref{e:defZ3} and~\eqref{e:defZ32}, we further expand
\begin{align*}
&Z_3^{\mathrm{res}}\big(\V,\W,\widehat{Z}\big)(k)\\ 
&\quad \ = \sum_{j \in \mathcal{T}} P_{j_0}(k) Z_3\big(P_{j_1} \V, P_{j_2} \W, P_{j_3} \widehat{Z}\big)(k) \\
&\quad \ = \sum_{j \in \mathcal{T}} P_{j_0}(k) \int_\R \int_\R \widetilde{N}_3(k,l_1,l_2)\big(P_{j_1} \V, P_{j_2} \W, P_{j_3} \widehat{Z}\big) \de l_1 \de l_2\\
&\quad \ = \sum_{j \in \mathcal{T}}  \int_\R \int_\R  \left[P_{j_0}(k)N_3(k,l_1,l_2)\big(P_{j_1} \V, P_{j_2} \W, P_{j_3} \widehat{Z}\big) \phantom{\sum_{m \in \{\pm 1\}}  \frac{P_{j_0}(k)}{\phi^2_{(j_0,m,j_3)}}} \right.\\ 
&\qquad \ \quad + \sum_{m \in \{\pm 1\}}  \frac{P_{j_0}(k)}{\phi^2_{(j_0,m,j_3)}(k,l_2)} N_2(k,l_2)\big(\chi P_m N_2(\cdot,l_1-l_2)  \big(P_{j_1} \V, P_{j_2} \W \big), P_{j_3} \widehat{Z} \big) \\
&\left. \quad \qquad\ + \sum_{m \in \{\pm 1\}} \frac{P_{j_0}(k)}{\phi^2_{(j_0,j_1,m)}(k,l_1)} N_2(k,l_1)\big(P_{j_1} \V , \chi P_{m} N_2(\cdot,l_2)\big(P_{j_2} \W, P_{j_3} \widehat{Z}\big)\big) \right] \de l_1 \de l_2,
\end{align*}
for $\V,\W,\widehat{Z} \in C_c(-k_0,k_0) \cap X_0$. So, we readily observe that there exist smooth functions $\mathcal{K}_j \colon (-k_0,k_0)^2 \to \C^2$ such that
\begin{align*}
Z_3^{\mathrm{res}}\big(\V,\W,\widehat{Z}\big)(0) &= \sum_{j \in \mathcal{T}} \int_\R \int_\R \mathcal{K}_j(l_1,l_2) \left(P_{j_1}(-l_1) \V(-l_1)\right)_1 \left(P_{j_2}(l_1-l_2) \W(l_1-l_2)\right)_1\\ &\qquad \cdot \left(P_{j_3}(l_2) \widehat{Z}(l_2)\right)_1 \de l_1 \de l_2,
\end{align*}
for $\V,\W,\widehat{Z} \in C_c(-k_0,k_0) \cap X_0$ and $j \in \mathcal{T}$. With the aid of~\eqref{e:specprojatcrit} we compute
\begin{align*}
\mathcal{K}_j(0,0) &= \left(\beta + \sum_{m \in \{\pm 1\}} \frac{\kappa^2 (P_m(0) \mathbf{e}_2)_1}{\phi^2_{j_0,m,j_3}(0,0)} + 
\sum_{m \in \{\pm 1\}} \frac{\kappa^2 (P_m(0) \mathbf{e}_2)_1}{\phi^2_{j_0,j_1,m}(0,0)}\right)P_{j_0}(0)\mathbf{e}_2\\
&= \left(\beta + \sum_{m \in \{\pm 1\}} \frac{\kappa^2 m}{2\left(m + j_3 - j_0\right)} + \sum_{m \in \{\pm 1\}} \frac{\kappa^2 m}{2\left(m + j_1 - j_0\right)}\right)P_{j_0}(0)\mathbf{e}_2,
\end{align*}
for $j \in \mathcal{T}$, which implies
\begin{align}  \label{e:compKj}
P_\pm(0) \sum_{j \in \mathcal{T}} \mathcal{K}_j(0,0) = \left(3\beta + \frac{10}{3} \kappa^2\right) P_\pm(0) \mathbf{e}_2.
\end{align}
In addition, by smoothness of $\mathcal{K}_j$ we have
\begin{align*}
\left|\mathcal{K}_j(l_1,l_2) - \mathcal{K}_j(0,0)\right| \lesssim |l_1| + |l_2|,
\end{align*}
for $l_1,l_2 \in (-k_0,k_0)$ and $j \in \mathcal{T}$. Thus, combining the latter with Lemma~\ref{l:linear estimates},~\eqref{e:sicest2} and~\eqref{e:sicapprox} and setting
\begin{align*}
\mathcal{Z}(t) &= \sum_{j \in \mathcal{T}} \int_\R \int_\R \mathcal{K}_j(0,0) \left(P_{j_1}(0) \re^{\tl_{j_1}(l_1)(t+1)} A(t) \right)_1  \left(P_{j_2}(0) \re^{\tl_{j_2}(l_1-l_2)(t+1)} A(t) \right)_1 \\ &\qquad \cdot \left(P_{j_3}(0) \re^{\tl_{j_3}(l_2)(t+1)} A(t) \right)_1 \de l_1 \de l_2,
\end{align*}
we obtain $t$-independent constants $C_{1,2} > 0$ such that
\begin{align} \label{e:Z3resapprox}
\begin{split}
&\left|Z_3^{\mathrm{res}}\big(\sic(t),\sic(t),\sic(t)\big)(0) - \mathcal{Z}(t)\right|\\ 
&\quad \ \leq C_1 \left(\int_\R \int_\R \left(|l_1| + |l_2|\right) \left|\re^{\hL(l_1)(t+1)}\right|\left|\re^{\hL(l_1-l_2)(t+1)}\right|\left|\re^{\hL(l_2)(t+1)}\right| |A(t)|^3 \de l_1 \de l_2\right.\\ 
&\quad \ \qquad \left.+ \, \frac{1}{(1+t)^{\frac32}} |A(t)|^3\right) \\
&\quad \ \leq \frac{C_2}{(1+t)^{\frac{3}{2}}} |A(t)|^3,
\end{split}
\end{align}
for $t \in [0,T_{\max})$. Using the standard integral
\begin{align} \label{e:int_standard}
\int_\R \re^{-a^2 z^2 + bz} \de z = \frac{\sqrt{\pi}}{a}\re^{\frac{b^2}{4a^2}},    
\end{align}
for $a, b \in \C$ with $\Re(a^2) > 0$, cf.~\cite[Integral~3.323.2]{GRADRY}, and recalling $j_0 - j_1 - j_2 - j_3 = 0$, we compute
\begin{align} \label{e:compint}
\begin{split}
&\int_\R \int_\R \re^{\left(\tl_{j_1}(l_1) + \tl_{j_2}(l_1-l_2) + \tl_{j_3}(l_2)\right)(t+1)} \de l_1 \de l_2\\ 
&\quad \ = \frac{2 \pi \re^{(j_1+j_2+j_3) \ri (t+1)}}{(t+1)\sqrt{3 \alpha^2 - j_1 j_2 - j_1 j_3 - j_2 j_3 - 2 \ri \alpha \left(j_1 + j_2 + j_3\right) }}\\
&\quad \ = \frac{2\pi \re^{j_0 \ri (t+1)}}{(t+1)\sqrt{3\alpha^2 + 1 - 2 \ri \alpha j_0}},
\end{split}
\end{align}
for $j \in \mathcal{T}$ and $t \geq 0$. Thus, invoking~\eqref{e:conjugateids},~\eqref{e:compKj} and~\eqref{e:compint}, we calculate
\begin{align} \label{e:ecompint2}
P_\pm(0) \mathcal{Z}(t) &= \frac{2\pi\left(9\beta + 10 \kappa^2\right)\re^{\pm \ri(t+1)}}{3(t+1) \sqrt{3\alpha^2 + 1 \mp 2 \ri \alpha}} \left|\left(P_{\pm}(0) A(t)\right)_1\right|^2 \left(P_{\pm}(0) A(t)\right)_1 P_{\pm}(0) \mathbf{e}_2,
\end{align}
for $t \in [0,T_{\max})$. The desired estimate now follows from~\eqref{e:specprojatcrit},~\eqref{e:E1bound},~\eqref{e:Z3resapprox} and~\eqref{e:ecompint2}, which concludes the proof.
\end{proof}

Since $A(t)$ is real, we have $|(P_\pm(0)A(t))_1| = \frac12 |A(t)|$ by~\eqref{e:specprojatcrit}. Hence, to control the norm of $A(t)$, it suffices to control the norm of the first coordinate of $P_\pm(0)A(t)$. We express the first coordinate of $P_\pm(0)A(t)$ in polar coordinates by setting 
\begin{align*}
(P_\pm(0) A(t))_1 \deq r_\pm(t) \re^{\ri \psi_\pm(t)}.\end{align*} 
Inserting this into the first component of~\eqref{e:Aeqpm}, multiplying with $\smash{\re^{-\ri\psi_\pm(t)}}$ and taking real parts, we arrive at
\begin{align} \label{e:eqrpm}
r_\pm'(t) = \frac{\Re(\omega_\pm)}{1+t} r_\pm(t)^3 + \Re\left(\re^{-\ri \psi_\pm(t) \mp \ri (t+1)}\left(P_\pm(0) \mathcal{E}_4(0,t) + \mathcal{E}_6(t)\right)_1\right).
\end{align}
As mentioned before, we expect that $\rc(k,t)$ exhibits higher-order decay since it vanishes at the critical frequency $k = 0$. In fact, we find that the bound on $\mathcal{E}_6(t)$, established in Lemma~\ref{l:leadingorderapproxZc}, is integrable in $t$, which implies that $\mathcal{E}_6(t)$ can be controlled using standard $L^1$-$L^\infty$-estimates in our nonlinear argument.

Thus, the critical term in the equation~\eqref{e:eqrpm} is the cubic $\frac{\omega_\pm}{1+t} r_\pm(t)^3$, that cannot be controlled through standard $L^1$-$L^\infty$ estimates. In fact, the sign of $\Re(\omega_\pm)$ turns out to be important, as we now see by solving the simpler ODE
\begin{align} \label{e:eqrODE} 
r_\pm'(t) = \frac{\Re(\omega_\pm)}{1+t} r_\pm(t)^3,
\end{align}
By separation of variables, it possesses the nonnegative solution
\begin{align*}
r_\pm(t) = \frac{r_\pm(0)}{\sqrt{1-2\Re(\omega_\pm) r_\pm(0)^2 \log(1+t)}},
\end{align*}
which exists globally in forward time and decays with rate $(\log(1+t))^{-\frac12}$ if $\Re(\omega_\pm) < 0$, and blows up in finite time if $\Re(\omega_\pm) > 0$ and $r_\pm(0) \neq 0$. One readily observes that the sign of $\Re(\omega_\pm)$ is equal to that of $9\beta + 10\kappa^2 = \frac{1}{8\pi^2}( 3N'''(0) + 5N''(0)^2)$. That is, $\Re(\omega_\pm)$ is negative if and only if~\eqref{e:signcondition} is fulfilled. 

In the upcoming nonlinear iteration we proceed as in~\cite[Theorem~5.1]{DRS} and apply a nonstandard method to control the radii $r_\pm(t)$ in case $\Re(\omega_\pm) < 0$. That is, we mimic the separation-of-variables procedure for~\eqref{e:eqrpm} and obtain
\begin{align} \label{e:lead}
\frac12 \partial_t \left(\frac{1}{r_\pm(t)^2}\right) = -\frac{\Re(\omega_\pm)}{1+t} - \Re\left(\frac{\re^{-\ri \psi_\pm(t) \mp \ri (t+1)}}{r_\pm(t)^3} \left(P_\pm(0) \mathcal{E}_4(0,t) + \mathcal{E}_6(t)\right)_1\right).
\end{align}
Then, integrating~\eqref{e:lead} and using that the residual
\begin{align*}
-\Re\left(\frac{\re^{-\ri \psi_\pm(t) \mp \ri (t+1)}}{r_\pm(t)^3} \left( P_\pm(0) \mathcal{E}_4(0,t) + \mathcal{E}_6(t)\right)_1\right)
\end{align*}
exhibits growth at a rate strictly smaller than $\log(1+t)$, one finds that the leading-order behavior of $r_\pm(t)^{-2}$ is given by $-2\Re(\omega_\pm) \log(1+t)$ yielding that $r_\pm(t)$ decays at rate $\smash{(\log(1+t))^{-\frac12}}$.
In summary, we establish that the leading-order behavior of $r_\pm(t) = \frac12 |A(t)|$ is governed by the separable nonlinear ODE~\eqref{e:eqrODE}.

\section{Proof of main results}
\label{s:theorem-proofs}

In this section, we prove our main results, Theorems~\ref{thm:mainresult} and~\ref{thm:mainresult2}. The proof of Theorem~\ref{thm:mainresult} relies on a global nonlinear iteration argument in the variables $\Uc,\Us$ and $\rc$ using $L^1$-$L^\infty$-estimates. In this iteration argument, we obtain control on the amplitude $|A(t)|$ of the leading-order Gaussian term through the nonlinear ODE~\eqref{e:eqrpm} by mimicking a separation-of-variable procedure, as outlined at the end of the previous section. This allows us to exploit the fact that the critical resonant nonlinear term is of absorption type.

\begin{proof}[Proof of Theorem~\ref{thm:mainresult}] 
We close a nonlinear iteration argument, controlling the variables $\Uc(k,t) = \chi(k) \U(k,t)$, $\Us(k,t) = (1-\chi(k))\U(k,t)$, $\U(k,t)$, and the remainder 
\begin{equation*}
\rc(t) = \smash{\re^{\hL(k)(t+1)} \re^{-\hL(0)(t+1)}} \chi(k) \Zc(0,t) - \Zc(k,t),
\end{equation*}
where we recall that $\Zc(k,t)$ can be expressed in terms of $\Uc(k,t)$ through~\eqref{e:normalform}. 
\bigskip\\
\noindent\textbf{Template function.} By~\eqref{e:regularity2},~\eqref{e:regularity3},~\eqref{e:regularity4} and~\eqref{e:regularity5} the template function $\eta \colon [0,T_{\max}) \to \R$ given by
\begin{align*}
\eta(t) &= \sup_{0 \leq s \leq t} \left[(\log(2+s))^{\frac23} \left(\big\|\rc(s)\big\|_{L^\infty} + \big\|\partial_k \rc(s)\big\|_{L^1} + \sqrt{1+s} \big\|\rc(s)\big\|_{L^1} \right) + 
\big\|\U(s)\big\|_{L^\infty} \right. \\
&\left.\left.\phantom{\frac{\big\|\Uc\big\|_{L^\infty}}{(t)^{\frac13}}} + \big\|\partial_k \U(s)\big\|_{L^1} +  \sqrt{1+s} \left(\frac{\big\||\cdot|\Uc(s)\big\|_{L^\infty} + \big\||\cdot|\partial_k \Uc(s)\big\|_{L^1}}{(\log(2+s))^{\frac13}} + \big\|\U(s)\big\|_{L^1}
\right. \right.\right. \\
&\left.\left.\phantom{\frac{\big\|\Uc\big\|_{L^\infty}}{(t)^{\frac13}}} +   \big\|\Us(s)\big\|_{L^\infty} + \big\|\partial_k \Us(s)\big\|_{L^1} \right) +  (1+s)\left(\big\|\Us(s)\big\|_{L^1} + \frac{\big\||\cdot|\Uc(s)\big\|_{L^1}}{(\log(2+s))^{\frac13}}\right)\right],
\end{align*}
is well-defined, continuous, monotonically increasing and, if $T_{\max} < \infty$, then it satisfies
\begin{align} \label{e:blowupeta}
\lim_{t \uparrow T_{\max}} \eta(t) = \infty,
\end{align}
by~\eqref{e:blowup}. 

Clearly,~\eqref{e:Fdef} implies
\begin{align} \label{e:boundUcUS}
\big\|\partial_k^j |\cdot|^m \Uc(t)\big\|_{L^p} + \big\|\partial_k^j \Us(t)\big\|_{L^p} \lesssim \big\|\U(t)\big\|_{W^{j,p}},
\end{align}
for $j = 0,1$, $m = 0,1$, $p = 1,\infty$ and $t \in [0,T_{\max})$, where we use that $\chi$ has compact support. Moreover,~\eqref{e:normalform},~\eqref{e:Q2bound} and~\eqref{e:K3bound} and the fact that $\chi$ has compact support, afford the estimate
\begin{align} \label{e:boundZC}
\big\||\cdot|^m \partial_k^j \big(\Zc(t) - \Uc(t)\big)\big\|_{L^p} \lesssim \big\|\Uc(t)\big\|_{W^{j,1}} \big\|\Uc(t)\big\|_{L^p},
\end{align}
implying
\begin{align} \label{e:boundZC2}
\big\|\partial_k^j \Zc(t) \big\|_{L^p} \lesssim \big\|\Uc(t)\big\|_{W^{j,p}} \lesssim  \big\|\U(t)\big\|_{W^{j,p}},
\end{align}
for $m = 0,1$, $j = 0,1$, $p = 1,\infty$ and $t \in [0,T_{\max})$ with $\|\U(t)\|_{W^{j,1} \cap L^\infty} \leq 1$ Furthermore, thanks to~\eqref{e:specprojatcrit},~\eqref{e:HL0expr} and~\eqref{e:boundZC2} it holds
\begin{align} \label{e:boundrrho}
r_\pm(t) &= \frac12 |A(t)| = \frac12 \big|\Zc(0,t)\big| \lesssim \big\|\U(t)\big\|_{L^\infty}, 
\end{align}
for all $t \in [0,T_{\max})$ with $\|\U(t)\|_{W^{j,1} \cap L^\infty} \leq 1$. Finally, from~\eqref{e:Gaussest} we have
\begin{align} \label{e:boundrc}
\begin{split}
\big\|\partial_k^j \rc(t)\big\|_{L^p} &\leq \big\|\partial_k^j \sic(t)\big\|_{L^p} + \big\|\partial_k^j \Zc(t)\big\|_{L^p} \lesssim (1+t)^{-\frac{1}{2p} + \frac{j}{2}} |A(t)| + \big\|\partial_k^j \Zc(t)\big\|_{L^p}\\ 
&\lesssim (1+t)^{-\frac{1}{2p} + \frac{j}{2}} \big\|\U(t)\big\|_{L^\infty} + \big\|\U(t)\big\|_{W^{j,p}},
\end{split}
\end{align}
for $j = 0,1$, $p = 1,\infty$ and $t \in [0,T_{\max})$ with $\|\U(t)\|_{W^{j,1} \cap L^\infty} \leq 1$. Thus, using~\eqref{e:cond_initial3},~\eqref{e:boundUcUS} and~\eqref{e:boundrc}, we conclude that there exists a constant $K_0 \geq 1$, independent of $E_0$, such that
\begin{align} \label{e:eta0est} \eta(0) \leq K_0E_0.\end{align}
\smallskip\\
\noindent\textbf{Approach.} Our goal is to prove that there exist constants $C \geq K_0$ and $\eta_0 \in (0,1)$ such that for all $t \in [0,T_{\max})$ with $\eta(t) \leq \eta_0$ the key inequality
\begin{align}
\eta(t) \leq C\left(E_0 + \eta(t)^2\right), \label{e:etaest}
\end{align}
is fulfilled. Then, taking
$$\eps = \min\left\{\frac{1}{4C^2},\frac{\eta_0}{2C}\right\},$$
it follows by the continuity of $\eta$ that, provided $E_0 \in (0,\eps)$, we have $\eta(t) \leq 2CE_0 \leq \eta_0$ for all $t \in [0,T_{\max})$. Indeed, given $t \in [0,T_{\max})$ with $\eta(t) > 2CE_0$, there must, by continuity of $\eta$ and the fact that $\eta(0) < 2CE_0$ by~\eqref{e:eta0est}, exist $s \in (0,t]$ with $\eta(s) = 2CE_0 \leq \eta_0$. Now~\eqref{e:etaest} yields 
$$\eta(s) \leq C\left(E_0 + 4C^2E_0^2\right) < 2CE_0,$$
contradicting $\eta(s) = 2CE_0$. We conclude that, if the key inequality~\eqref{e:etaest} holds, then we have $\eta(t) \leq 2CE_0$ for all $t \in [0,T_{\max})$. In this case, we must have $T_{\max} = \infty$ by~\eqref{e:blowupeta}. We then infer $\eta(t) \leq 2CE_0$ for all $t \geq 0$, which leads to the estimates~\eqref{e:diffest} and~\eqref{e:diffest2}. Hence, all that remains is to prove the key inequality~\eqref{e:etaest} and derive the estimates~\eqref{e:diffest} and~\eqref{e:diffest2}. 
\bigskip\\
\noindent\textbf{Controlling $r_\pm(t)$.}  First of all, by~\eqref{e:boundrrho} there exists an $E_0$-independent constant $K_1 > 0$ such that
\begin{align} \label{e:rpm0}
r_\pm(0) &\leq K_1 E_0. 
\end{align}
Define the set
\begin{align*}
\mathcal{S} = \left\{t \in [0,T_{\max}) : r_\pm(t) \sqrt{\left(K_1E_0 + \eta(t)^{\frac32}\right)^{-2} - 2 \Re(\omega_\pm) \log(1+t)} \leq 2\right\}.
\end{align*}
By~\eqref{e:rpm0} it holds 
\begin{align} \label{e:rpm01}
r_\pm(0) \sqrt{\left(K_1E_0 + \eta(0)^{\frac32}\right)^{-2}} \leq 1,
\end{align}
implying $0 \in \mathcal{S}$. Our aim is to show that, we have, $\mathcal{S} = \{t \in [0,T_{\max}) : \eta(t) \leq \eta_0\}$, provided $\eta_0 > 0$ is sufficiently small (but independent of $E_0$). We argue by contradiction and assume that there exists $t \in [0,T_{\max}) \setminus \mathcal{S}$ with $\eta(t) \leq \eta_0$. So, by identity~\eqref{e:rpm01} and continuity of $r_\pm$ and $\eta$, there must exist $t_1 \in (0,t)$ such that
\begin{align*} 
r_\pm(s)\sqrt{\left(K_1E_0 + \eta(s)^{\frac32}\right)^{-2} - 2 \Re(\omega_\pm) \log(1+s)} \geq 1,\\
r_\pm(t_1)\sqrt{\left(K_1E_0 + \eta(t_1)^{\frac32}\right)^{-2} - 2 \Re(\omega_\pm) \log(1+t_1)} = 1, 
\end{align*}
for $s \in [t_1,t]$, implying
\begin{align} 
\begin{split}
r_\pm(s)^{-2} &\leq \left(K_1E_0 + \eta(s)^{\frac32}\right)^{-2} - 2 \Re(\omega_\pm) \log(1+s),\\
r_\pm(t_1)^{-2} &= \left(K_1E_0 + \eta(t_1)^{\frac32}\right)^{-2} - 2 \Re(\omega_\pm) \log(1+t_1), 
\end{split}
\label{e:aid}
\end{align}
for $s \in [t_1,t]$. Integrating~\eqref{e:lead} from $t_1$ to $t$ and employing~\eqref{e:aid}, we arrive at
\begin{align}
\begin{split}
r_\pm(t)^{-2} &= \left(K_1E_0 + \eta(t_1)^{\frac32}\right)^{-2} - 2\Re(\omega_\pm) \log(t+1)\\ 
&\qquad -\, 2\int_{t_1}^t \Re\left(\frac{\re^{-\ri \psi_\pm(s)\mp \ri (s+1)}}{r_\pm(s)^3} \left(P_\pm(0) \mathcal{E}_4(0,s) + \mathcal{E}_6(s)\right)_1\right) \de s. 
\end{split}
\label{e:INTr}
\end{align}
The next steps are devoted to estimating the integral term in~\eqref{e:INTr}. Lemmas~\ref{l:residualnonlZc} and~\ref{l:leadingorderapproxZc}, Young's inequality and estimate~\eqref{e:boundrrho} yield a $t$- and $E_0$-independent constant $C_1 > 0$ such that
\begin{align} \label{e:INT}
\begin{split}
|\mathcal{E}_4(0,s)| &\leq C_1\frac{\eta(s)^2 (\log(2+s))^{\frac13}}{(1+s)^{\frac32}},\\
|\mathcal{E}_6(s)| &\leq \frac{C_1}{1+s} \left(\frac{\eta(s)^3}{(\log(2+s))^{2}} + \frac{\eta(s) r_\pm(s)^2}{(\log(2+s))^{\frac23}}\right),
\end{split}
\end{align}
for each $s \in [t_1,t]$, where we use $\eta(t) \leq \eta_0 \leq 1$. 

Thus, using $\eta(t) \leq \eta_0 \leq 1$,~\eqref{e:aid} and~\eqref{e:INT}, we obtain $t$- and $E_0$-independent constants $C_2,C_3,C_4 > 0$ such that
\begin{align}\begin{split}
&\left|\int_{t_1}^t \Re\left(\frac{\re^{-\ri \psi_\pm(s)\mp \ri (s+1)}}{r_\pm(s)^3} P_\pm(0) \mathcal{E}_4(0,s)\right) \de s\right|\\ 
&\quad \ \leq C_2 \int_{t_1}^t \eta(s)^2 \frac{\left(\left(K_1E_0 + \eta(s)^{\frac32}\right)^{-2} + \log(1+s)\right)^{\frac32}(\log(2+s))^{\frac13}}{(1+s)^{\frac32}} \de s\\
&\quad \ \leq C_3 \sqrt{\eta(t)}\left(\left(K_1E_0 + \eta(t_1)^{\frac32}\right)^{-2} - 2 \Re(\omega_\pm)\log(1+t)\right)\\
&\quad \ \qquad \cdot \int_{t_1}^t \left(\frac{\eta(s)^{\frac32}(\log(2+s))^{\frac13}}{\left(K_1E_0 + \eta(s)^{\frac32}\right)(1+s)^{\frac32}} + \frac{(\log(2+s))^{\frac56}}{(1+s)^{\frac32}}\right) \de s \\
&\quad \ \leq C_4\sqrt{\eta_0}\left(\left(K_1 E_0 + \eta(t_1)^{\frac32}\right)^{-2} - 2 \Re(\omega_\pm)\log(1+t)\right).
\end{split}\label{e:E4est1}\end{align}
Next, we proceed with bounding the contribution from $\mathcal{E}_6$ in~\eqref{e:INTr}. We split the estimate in two parts corresponding to the two terms on the right-hand side of the bound~\eqref{e:INT} on $|\mathcal{E}_6(s)|$. The estimate associated with the first term is similar to~\eqref{e:E4est1}. That is, using $\eta(t) \leq \eta_0 \leq 1$ and~\eqref{e:aid}, we obtain $t$- and $E_0$-independent constants $C_2,C_3,C_4 > 0$ such that
\begin{align}\begin{split}
&\int_{t_1}^t \frac{\eta(s)^3}{r_\pm(s)^3(1+s)(\log(2+s))^2} \de s\\ 
&\quad \ \leq C_2 \int_{t_1}^t \eta(s)^3 \frac{\left(\left(K_1E_0 + \eta(s)^{\frac32}\right)^{-2} + \log(1+s)\right)^{\frac32}}{(1+s)(\log(2+s))^{2}} \de s\\
&\quad \ \leq C_3 \eta(t)\left(\left(K_1E_0 + \eta(t_1)^{\frac32}\right)^{-2} - 2 \Re(\omega_\pm)\log(1+t)\right)\\
&\qquad \quad \ \cdot \int_{t_1}^t \left(\frac{\eta(s)^{\frac32}}{\left(K_1E_0 + \eta(s)^{\frac32}\right)(1+s)(\log(2+s))^{2}} + \frac{1}{(1+s)(\log(2+s))^{\frac{3}{2}}}\right) \de s \\
&\quad \ \leq C_4\eta_0\left(\left(K_1E_0 + \eta(t_1)^{\frac32}\right)^{-2} - 2 \Re(\omega_\pm)\log(1+t)\right).
\end{split}\label{e:Aest1}\end{align}
We proceed with the estimate associated with the second term. Here, we use $\eta(t) \leq \eta_0$ and~\eqref{e:aid} to establish $t$- and $E_0$-independent constants $C_2,C_3 >0$ such that
\begin{align}\begin{split}
&\int_{t_1}^t \frac{\eta(s)}{r_\pm(s)(1+s)(\log(2+s))^{\frac{2}{3}}} \de s\\ 
&\quad \ \leq C_2 \eta(t) \int_{t_1}^t \frac{\sqrt{\log(1+s)}\sqrt{\left(K_1E_0 + \eta(s)^{\frac32}\right)^{-2} - 2 \Re(\omega_\pm)\log(1+s)}}{(1+s)(\log(2+s))^{\frac{2}{3}}\sqrt{\log(1+s)}} \de s \\
&\quad \ \leq C_2 \eta(t) \int_{t_1}^t \frac{\sqrt{\log(1+t)}\sqrt{\left(K_1E_0 + \eta(t_1)^{\frac32}\right)^{-2} - 2 \Re(\omega_\pm)\log(1+t)}}{(1+s)(\log(2+s))^{\frac{2}{3}}\sqrt{\log(1+s)}} \de s \\
&\quad \ \leq C_3\eta_0\left(\left(K_1E_0 + \eta(t_1)^{\frac32}\right)^{-2} - 2 \Re(\omega_\pm)\log(1+t)\right).
\end{split}\label{e:Aest2}\end{align}
Thus, combining~\eqref{e:INT},~\eqref{e:Aest1} and~\eqref{e:Aest2} and recalling $\eta_0 \leq 1$, we obtain a $t$- and $E_0$-independent constant $C_4 > 0$ such that
\begin{align} \label{e:Aestfinal}
\begin{split}
&\left|\int_{t_1}^t \Re\left(\frac{\re^{-\ri \psi_\pm(s) \mp \ri (s+1)}}{r_\pm(s)^3} \mathcal{E}_6(s)\right) \de s\right|\\ 
&\quad \ \leq C_4 \sqrt{\eta_0}\left(\left(K_1E_0 + \eta(t_1)^{\frac32}\right)^{-2} - 2 \Re(\omega_\pm)\log(1+t)\right).
\end{split}
\end{align}

Finally, we apply~\eqref{e:E4est1} and~\eqref{e:Aestfinal} to estimate the integral in~\eqref{e:INTr}. So, provided $\eta_0 > 0$ is smaller than $1/(8C_4)^2$, we arrive at the lower bound
\begin{align*}
r_\pm(t)^{-2} &\geq \left(1-4C_4\sqrt{\eta_0}\right)\left(\left(K_1 E_0 + \eta(t_1)^{\frac32}\right)^{-2} - 2 \Re(\omega_\pm) \log(t+1)\right)\\ 
&\geq \frac{1}{2}\left(\left(K_1 E_0 + \eta(t)^{\frac32}\right)^{-2} - 2 \Re(\omega_\pm) \log(t+1)\right),
\end{align*}
implying
\begin{align*}
r_\pm(t) \sqrt{\left(K_1 E_0 + \eta(t)^{\frac32}\right)^{-2} - 2 \Re(\omega_\pm) \log(t+1)} \leq \sqrt{2} < 2,
\end{align*}
which contradicts $t \notin \mathcal{S}$. We conclude that
\begin{align}
\begin{split}
r_\pm(t) &\leq \frac{2}{\sqrt{\left(K_1E_0 + \eta(t)^{\frac32}\right)^{-2} - 2 \Re(\omega_\pm) \log(1+t)}}\\ 
&\leq 2\min\left\{K_1E_0 + \eta(t)^{\frac32},\frac{1}{\sqrt{- 2 \Re(\omega_\pm) \log(1+t)}}\right\}, 
\end{split}
\label{e:restfin}
\end{align}
holds for all $t \in [0,T_{\max})$ with $\eta(t) \leq \eta_0$.
\newpage
\noindent\textbf{Bounds on $\rc(t)$.} Let $t \in [0,T_{\max})$ with $\eta(t) \leq \eta_0 \leq 1$ and $E_0 \in (0,\eps)$. We bound the residual $\rc(t)$ through its Duhamel formula
\begin{align}
\label{e:rcDuh}
\begin{split}
\rc(k,t) &= \re^{\hL(k)t} \rc(k,0) + \int_0^t \re^{\hL(k)(t-s)} \chi(k) \mathcal{E}_7(k,s) \de s,
\end{split}
\end{align}
which arises by integrating~\eqref{e:rceq}, where we denote
\begin{align*}
\mathcal{E}_7(k,s) = \mathcal{E}_4(k,s) - \re^{\hL(k)(s+1)} \re^{-\hL(0)(s+1)} \mathcal{E}_4(0,s) +  \mathcal{E}_5(k,s) .
\end{align*}

We start by bounding the nonlinear term in~\eqref{e:rcDuh}. The fact that $\mathcal{E}_5(0,t) = 0$, the mean value theorem, Lemma~\ref{l:cubicresidualZc} and estimates~\eqref{e:boundrrho} and~\eqref{e:restfin} yield a $t$- and $E_0$-independent constant $C_1 > 0$ such that
\begin{align} \label{e:rcnonl}
\begin{split}
|\mathcal{E}_5(k,s)| &\leq \|\mathcal{E}_5(\cdot,s)\|_{L^\infty} \leq C_1 \frac{E_0 + \eta(s)^2}{(1+s) (\log(2+s))^{\frac23}},\\
|\partial_k \mathcal{E}_5(k,s)| &\leq \|\partial_k \mathcal{E}_5(\cdot,s)\|_{L^\infty} \leq C_1 \frac{E_0 + \eta(s)^2}{\sqrt{1+s} (\log(2+s))^{\frac23}},\\
|\mathcal{E}_5(k,s)| &\leq |k| \|\partial_k \mathcal{E}_5(\cdot,s)\|_{L^\infty} \leq C_1 |k| \frac{E_0 + \eta(s)^2}{\sqrt{1+s} (\log(2+s))^{\frac23}},
\end{split}
\end{align}
hold for $s \in [0,t]$ and $k \in \R$, where we use $\eta(t) \leq \eta_0 \leq 1$ and $E_0 \in (0,\eps)$. Similarly, we employ the mean value theorem, identity~\eqref{e:HL0expr}, and Lemmas~\ref{l:linear estimates} and~\ref{l:residualnonlZc} to establish a $t$- and $E_0$-independent constant $C_1 > 0$ such that
\begin{align*}
\mathcal{E}_8(k,s) \deq \mathcal{E}_4(k,s) - \re^{\hL(k)(s+1)} \re^{-\hL(0)(s+1)} \mathcal{E}_4(0,s)
\end{align*}
can be bounded as
\begin{align} \label{e:rcnonl2}
\begin{split}
\left|\chi(k) \mathcal{E}_8(k,s)\right| &\leq \|\mathcal{E}_4(\cdot,s)\|_{L^\infty} \leq C_1 \frac{\eta(s)^2 (\log(2+s))^{\frac13}}{(1+s)^{\frac32}},\\
\left|\partial_k \left(\chi(k)\mathcal{E}_8(k,s)\right)\right| &\leq \|\chi \mathcal{E}_4(\cdot,s)\|_{W^{1,\infty}} + \left\|\partial_k \re^{\hL(\cdot)(s+1)}\right\|_{L^\infty} \|\chi \mathcal{E}_4(\cdot,s)\|_{L^\infty}\\ 
&\leq C_1 \frac{\eta(s)^2 (\log(2+s))^{\frac13}}{1+s}, \\
\left|\chi(k) \mathcal{E}_8(k,s)\right| &\leq |k| \left\|\partial_k \left(\chi \mathcal{E}_8(\cdot,s)\right)\right\|_{L^\infty} \leq C_1 |k| \frac{\eta(s)^2 (\log(2+s))^{\frac13}}{1+s},
\end{split}
\end{align}
for $s \in [0,t]$ and $k \in \R$, where we use $\eta(t) \leq \eta_0 \leq 1$ and $\mathcal{E}_8(0,s) = 0$. 

In case $t \geq 1$, we apply Lemma~\ref{l:linear estimates}, employ the estimates~\eqref{e:rcnonl} and~\eqref{e:rcnonl2}, and use that $\chi$ is supported on $(-k_0,k_0)$ to arrive at
\begin{align} \label{e:nlest1}
\begin{split}
&\left\|\int_0^t |\cdot|^j \re^{\hL(\cdot)(t-s)}  \chi \mathcal{E}_7(\cdot,s)\de s \right\|_{L^1}\\ 
&\quad \ \leq C_2\left(\int_0^{\frac{t}{2}} \int_{-k_0}^{k_0} \frac{E_0 + \eta(s)^2}{\sqrt{1+s} (\log(2+s))^{\frac23}} |k|^{1+j} \re^{-\frac12 \alpha k^2 (t-s)} \de k \de s\right.\\ 
&\quad  \ \qquad \left. + \, \int_{\frac{t}{2}}^t \int_{-k_0}^{k_0} \frac{E_0 + \eta(s)^2}{(1+s) (\log(2+s))^{\frac23}} |k|^j \re^{-\frac12 \alpha k^2 (t-s)} \de k  \de s\right)\\
&\quad \ \leq C_3\left(\int_0^{\frac{t}{2}} \frac{E_0+\eta(s)^2}{\sqrt{t-s}(1+t-s)^{\frac{1+j}{2}}\sqrt{1+s} (\log(2+s))^{\frac23}} \de s\right. \\
&\quad \ \qquad \left. + \, \int_{\frac{t}{2}}^t \frac{E_0 + \eta(s)^2}{\sqrt{t-s}(1+t-s)^{\frac{j}{2}}(1+s) (\log(2+s))^{\frac23}} \de s\right)\\
&\quad \ \leq C_4\frac{E_0 + \eta(t)^2}{(1+t)^{\frac{1+j}{2}} (\log(2+t))^{\frac23- j}},
\end{split}
\end{align}
for $j = 0,1$ and some $t$- and $E_0$-independent constants $C_2,C_3, C_4 > 0$. In case $t \leq 1$, we use Lemma~\ref{l:linear estimates}, the fact that $\chi$ has compact support and estimates~\eqref{e:rcnonl} and~\eqref{e:rcnonl2}, to establish the short-time bound
\begin{align} \label{e:nlest2}
\begin{split}
\left\||\cdot|^j \int_0^t \re^{\hL(\cdot)(t-s)}  \chi \mathcal{E}_7(\cdot,s)\de s \right\|_{L^1}
&\leq C_2\int_0^t \int_\R \frac{E_0 + \eta(s)^2}{(1+s) (\log(2+s))^{\frac23}} \re^{-\frac12 \alpha k^2 (t-s)} \de k  \de s\\
&\leq C_3\int_0^t  \frac{E_0 + \eta(s)^2}{\sqrt{t-s}} \de s\\ &\leq C_4\frac{E_0 + \eta(t)^2}{(1+t)^{\frac{1+j}{2}} (\log(2+t))^{\frac23- j}},
\end{split}
\end{align}
for $j = 0,1$ and some $t$- and $E_0$-independent constants $C_2,C_3, C_4 > 0$. Furthermore, Lemma~\ref{l:linear estimates}, the fact that $\chi$ is supported on $(-k_0,k_0)$ and estimates~\eqref{e:rcnonl} and~\eqref{e:rcnonl2} yield
\begin{align} \label{e:nlest3}
\begin{split}
&\left\||\cdot|^j \partial_k \int_0^t \re^{\hL(\cdot)(t-s)}  \chi \mathcal{E}_7(\cdot,s)\de s \right\|_{L^1}\\
&\quad \ \leq C_2 \int_0^t \int_{-k_0}^{k_0} \frac{\left(E_0 + \eta(t)^2\right)\left(1 + |k|^2 (t-s)\right)}{\sqrt{1+s} (\log(2+s))^{\frac23}} |k|^j\re^{-\frac12 \alpha k^2 (t-s)} \de k \de s\\ 
&\quad \ \leq C_3\int_0^t \frac{E_0+\eta(s)^2}{\sqrt{t-s} (1+t-s)^{\frac{j}{2}} \sqrt{1+s} (\log(2+s))^{\frac23}} \de s\\ 
&\quad \ \leq C_4\frac{E_0 + \eta(t)^2}{(1+t)^{\frac{j}{2}}(\log(2+t))^{\frac23-j}},
\end{split}
\end{align}
for $j = 0,1$ and some $t$- and $E_0$-independent constants $C_2,C_3, C_4 > 0$. Similarly, the $L^\infty$-estimate
\begin{align*}
\left\|\int_0^t |\cdot|^j \re^{\hL(\cdot)(t-s)}  \chi \mathcal{E}_7(\cdot,s)\de s \right\|_{L^\infty}
&\leq C_2\int_0^t \frac{E_0+\eta(s)^2}{\sqrt{t-s} (1+t-s)^{\frac{j}{2}}\sqrt{1+s} (\log(2+s))^{\frac23}}  \de s\\ &\leq C_3\frac{E_0 + \eta(t)^2}{(1+t)^{\frac{j}{2}}(\log(2+t))^{\frac23-j}},
\end{align*}
follows for some $t$- and $E_0$-independent constants $C_2,C_3 > 0$. 

We proceed with bounding the linear term in~\eqref{e:rcDuh}. By Lemma~\ref{l:linear estimates}, estimate~\eqref{e:boundrc} the mean value theorem and the facts that $\rc(0,0) = 0$ and that $\rc(\cdot,0)$ is supported on $(-k_0,k_0)$, there exist $t$- and $E_0$-independent constants $C_2,C_3 > 0$ such that for $p = 1,\infty$ and $j = 0,1$, we have 
\begin{align} \label{e:linst1a}
\begin{split}
\left\||\cdot|^j \re^{\hL(\cdot)t} \rc(0)\right\|_{L^p} &\leq C_2\left\|\rc(0)\right\|_{L^p} \leq  C_3 E_0,\\
\left\||\cdot|^j \partial_k \re^{\hL(\cdot)t} \rc(0)\right\|_{L^1} &\leq C_2 \left\|\rc(0)\right\|_{W^{1,1}} \leq C_3 E_0,
\end{split}
\end{align}
for $t\leq 1$, and 
\begin{align} \label{e:linst1b}
\begin{split}
\left\||\cdot|^j \re^{\hL(\cdot)t} \rc(0)\right\|_{L^1} &\leq C_2 \int_\R |k|^{1+j} \re^{-\frac12 \alpha  k^2 t} \left\|\partial_k \rc(0)\right\|_{L^\infty} \de k \leq C_3 \frac{E_0}{t^{1+\frac{j}{2}}},\\
\left\||\cdot|^j \re^{\hL(\cdot)t} \rc(0)\right\|_{L^\infty} &\leq C_2 \sup_{k \in \R} |k|^{1+j} \re^{-\frac12 \alpha  k^2 t} \left\|\partial_k \rc(0)\right\|_{L^\infty} \leq C_3 \frac{E_0}{t^{\frac{1+j}{2}}},\\
\left\||\cdot|^j \partial_k \re^{\hL(\cdot)t} \rc(0)\right\|_{L^1} &\leq C_2\Big( 
t \left\||\cdot|^{1+j} \re^{-\frac12 \alpha |\cdot|^2 t} \rc(0)\right\|_{L^1}\\ 
&\qquad \, + \left\| |\cdot|^j \re^{-\frac12 \alpha |\cdot|^2 t} \partial_k \rc(0)\right\|_{L^1}\Big) \leq C_3 \frac{E_0}{t^{\frac{1+j}{2}}},
\end{split}
\end{align}
for $t\geq 1$.

All in all, applying~\eqref{e:nlest1},~\eqref{e:nlest2},~\eqref{e:nlest3} and~\eqref{e:linst1a}-\eqref{e:linst1b} to bound the right-hand side of~\eqref{e:rcDuh} yields a $t$- and $E_0$-independent constant $C_* > 0$ such that
\begin{align} \label{e:finalrc1}
\big\|\rc(t)\big\|_{L^1} \leq \frac{C_*\left(E_0 + \eta(t)^2\right)}{\sqrt{1+t} (\log(2+t))^{\frac23}}, \qquad \big\|\rc(t)\big\|_{L^\infty}, \big\|\partial_k \rc(t)\big\|_{L^1} \leq\frac{C_*\left(E_0 + \eta(t)^2\right)}{(\log(2+t))^{\frac23}},
\end{align}
and
\begin{align} \label{e:finalrc2}
\begin{split}
\big\||\cdot| \partial_k^j \rc(t)\big\|_{L^1} \leq\frac{C_*\left(E_0 + \eta(t)^2\right)(\log(2+t))^{\frac13}}{(1+t)^{1-\frac{j}{2}}},\\
\big\||\cdot| \rc(t)\big\|_{L^\infty} \leq\frac{C_*\left(E_0 + \eta(t)^2\right)(\log(2+t))^{\frac13}}{(1+t)^{\frac{1}{2}}},
\end{split}
\end{align}
for $j = 0,1$.
\bigskip\\
\noindent\textbf{Bounds on $\Us(t)$.} Let $t \in [0,T_{\max})$ with $\eta(t) \leq \eta_0 \leq 1$ and $E_0 \in (0,\eps)$. We derive estimates on $\Us(t)$ by bounding the right-hand side of its Duhamel formula~\eqref{e:Duhdamp}. First, we apply Lemma~\ref{l:noncriticalbounds} and obtain a $t$- and $E_0$-independent constant $C_1 > 0$ such that
\begin{align} \label{e:nonlest00}
\left\|\partial_k^j \mathcal{F} \mathcal{N}\left(\mathcal{F}^{-1} \U(s)\right)\right\|_{L^p} \leq \frac{C_1 \eta(s)^2}{(1+s)^{\frac{1}{2p} + \frac12 (1-j)}},
\end{align}
for $j = 0,1$, $p = 1, \infty$ and $s \in [0,t]$, where we use $\eta(t) \leq 1$. Next, we evoke Lemma~\ref{l:linear estimates} and estimate~\eqref{e:nonlest00} to arrive at
\begin{align} \label{e:nlest5}
\begin{split}
&\left\|\partial_k^j \int_0^t \re^{\hL(\cdot)(t-s)} (1-\chi) \mathcal{F} \mathcal{N}\left(\mathcal{F}^{-1} \U(s)\right) \de s \right\|_{L^p}\\ 
&\quad \ \leq C_2 \int_0^t \re^{-\theta(t-s)} \left\|\mathcal{F} \mathcal{N}\left(\mathcal{F}^{-1} \U(s)\right)\right\|_{W^{j,p}} \de s \\
&\quad \ \leq C_3 \int_0^t \frac{\re^{-\theta(t-s)} \eta(s)^2}{(1+s)^{\frac{1}{2p} + \frac12 (1-j)}} \de s \leq C_4 \frac{\eta(t)^2}{(1+t)^{\frac{1}{2p} + \frac12 (1-j)}},
\end{split}
\end{align}
for some $t$- and $E_0$-independent constants $C_2,C_3, C_4 > 0$. On the other hand, Lemma~\ref{l:linear estimates} readily yields the linear estimate
\begin{align} \label{e:linst2}
\begin{split}
\left\|\partial_k^j \re^{\hL(\cdot)t} (1-\chi)\U(0)\right\|_{L^p} \leq C_2 \re^{-\theta t} \|\U(0)\|_{W^{j,p}} \leq C_2  \re^{-\theta t} E_0,
\end{split}
\end{align}
for some $t$- and $E_0$-independent constant $C_2 > 0$. Thus, applying~\eqref{e:nlest5} and~\eqref{e:linst2} to the right-hand side of~\eqref{e:Duhdamp} we obtain a $t$- and $E_0$-independent constant $C_* > 0$ such that
\begin{align} \label{e:finalUs}
\big\|\Us(t)\big\|_{L^1} \leq C_*\frac{E_0 + \eta(t)^2}{1+t}, \qquad \big\|\Us(t)\big\|_{L^\infty}, \big\|\partial_k \Us(t)\big\|_{L^1} \leq C_*\frac{E_0 + \eta(t)^2}{\sqrt{1+t}}.
\end{align}

\noindent\textbf{Bounds on $\Uc(t)$ and $\U(t)$.} Let $t \in [0,T_{\max})$ with $\eta(t) \leq \eta_0 \leq 1$ and $E_0 \in (0,\eps)$. Applying Lemma~\ref{l:linear estimates} and using that $\chi$ is supported on $(-k_0,k_0)$, we deduce
\begin{align*}
\left\||\cdot|^m \partial_k^j \left(\chi \re^{\hL(\cdot)(s+1)}\right)\right\|_{L^1} &\lesssim \int_{-k_0}^{k_0} |k|^m (1 + j |k| (s+1)) \re^{-\frac{1}{2} \alpha k^2 (s+1)} \de k \lesssim \frac{1}{(1+s)^{\frac{1}{2} + \frac{m-j}{2}}},\\
\left\||\cdot|^m \left(\chi \re^{\hL(\cdot)(s+1)}\right)\right\|_{L^\infty} &\lesssim \sup_{k \in (-k_0,k_0)} |k|^m \re^{-\frac{1}{2} \alpha k^2(s+1)} \lesssim \frac{1}{(1+s)^{\frac{m}{2}}}
\end{align*}
for $s \geq 0$, $p = 1,\infty$, $j = 0,1$ and $m = 0,1$. So, recalling~\eqref{e:defAsic},~\eqref{e:decompZc} and~\eqref{e:boundrrho} and employing the estimates~\eqref{e:restfin},~\eqref{e:finalrc1} and~\eqref{e:finalrc2}, we obtain a $t$- and $E_0$-independent constant $C_1 > 0$ such that
\begin{align*}
\big\||\cdot|^m \Zc(t)\big\|_{L^p} &\leq \big\||\cdot|^m\sic(t)\big\|_{L^p} + \big\||\cdot|^m \rc(t)\big\|_{L^p}\\ 
&\leq |A(t)| \left\||\cdot|^m \left(\chi \re^{\hL(\cdot)(t+1)}\right)\right\|_{L^p}  + \big\||\cdot|^m \rc(t)\big\|_{L^p} \\
&\leq C_1 \frac{\left(E_0 + \eta(t)^{\frac32}\right) (\log(2+t))^{\frac{m}{3}}}{(1+t)^{\frac{1}{2p} + \frac{m}{2}}},
\end{align*}
and
\begin{align*}
\big\||\cdot|^m \partial_k \Zc(t)\big\|_{L^1} &\leq \big\||\cdot|^m\partial_k \sic(t)\big\|_{L^1} + \big\||\cdot|^m\partial_k \rc(t)\big\|_{L^1}\\ 
&\leq |A(t)| \left\||\cdot|^m \partial_k \left(\chi \re^{\hL(\cdot)(t+1)}\right)\right\|_{L^1}  + \big\||\cdot|^m \partial_k \rc(t)\big\|_{L^1} \\
&\leq C_1 \frac{\left(E_0 + \eta(t)^{\frac32}\right) (\log(2+t))^{\frac{m}{3}}}{(1+t)^{\frac{m}{2}}},
\end{align*}
for $p = 1,\infty$ and $m = 0,1$, where we use $\eta(t) \leq \eta_0 \leq 1$. Combining the latter with~\eqref{e:boundZC} we find $t$- and $E_0$-independent constants $C_2,C_3,C_4 > 0$ such that  
\begin{align} \label{e:finalUC1}
\begin{split}
\big\||\cdot|^m \Uc(t)\big\|_{L^p} &\leq C_2\left(\big\||\cdot|^m \Zc(t)\big\|_{L^p} + \big\|\U(t)\big\|_{L^1} \big\|\U(t)\big\|_{L^p}\right)  \\
&\leq C_3 \left(\frac{\left(E_0 + \eta(t)^{\frac32}\right) (\log(2+t))^{\frac{m}{3}}}{(1+t)^{\frac{1}{2p} + \frac{m}{2}}} + \frac{\eta(t)^2}{(1+t)^{\frac12 + \frac1{2p}}}\right)\\
&\leq C_4 \frac{\left(E_0 + \eta(t)^{\frac32}\right) (\log(2+t))^{\frac{m}{3}}}{(1+t)^{\frac{1}{2p} + \frac{m}{2}}},
\end{split}
\end{align}
and
\begin{align}  \label{e:finalUC2}
\begin{split}
\big\||\cdot|^m \partial_k \Uc(t)\big\|_{L^1} &\leq C_2\left(\big\||\cdot|^m \partial_k \Zc(t)\big\|_{L^1} + \big\|\U(t)\big\|_{W^{1,1}} \big\|\U(t)\big\|_{L^1}\right)  \\
&\leq C_3 \left(\frac{\left(E_0 + \eta(t)^{\frac32}\right) (\log(2+t))^{\frac{m}{3}}}{(1+t)^{\frac{m}{2}}} + \frac{\eta(t)^2}{\sqrt{1+t}}\right)\\
&\leq C_4 \frac{\left(E_0 + \eta(t)^{\frac32}\right) (\log(2+t))^{\frac{m}{3}}}{(1+t)^{\frac{m}{2}}},
\end{split}
\end{align}
for $p = 1,\infty$ and $m = 0,1$, where we use $\eta(t) \leq \eta_0 \leq 1$. Finally, identity~\eqref{e:Fdef} and estimates~\eqref{e:finalUs},~\eqref{e:finalUC1} and~\eqref{e:finalUC2} yield
\begin{align} \label{e:finalU}
\big\|\U(t)\big\|_{L^1} \leq \frac{C_*\left(E_0 + \eta(t)^{\frac32}\right)}{\sqrt{1+t}}, \qquad \big\|\U(t)\big\|_{L^\infty}, \big\|\partial_k \U(t)\big\|_{L^1} \leq C_*\left(E_0 + \eta(t)^{\frac32}\right),
\end{align}
for some $t$- and $E_0$-independent constant $C_*> 0$. 
\bigskip\\
\noindent\textbf{Proof of key inequality.} By the estimates~\eqref{e:finalrc1},~\eqref{e:finalUs},~\eqref{e:finalUC1},~\eqref{e:finalUC2} and~\eqref{e:finalU} there exist $t$- and $E_0$-independent constants $C_0,C > 0$ such that
\begin{align*}
\eta(t) \leq C_0\left(E_0 + \eta(t)^{\frac32}\right) \leq C_0\left(E_0 + C_0^{\frac32}\left(E_0 + \eta(t)^{\frac32}\right)^{\frac32}\right) \leq C\left(E_0 + \eta(t)^2\right)
\end{align*}
for $t \in [0,T_{\max})$ with $\eta(t) \leq \eta_0 \leq 1$ and $E_0 \in (0,\eps)$. Therefore, we have established the key inequality~\eqref{e:etaest}, which, as argued before, implies $T_{\max} = \infty$ and $\eta(t) \leq 2CE_0$ for all $t \geq 0$. 
\bigskip\\
\noindent\textbf{Proof of the estimates~\eqref{e:diffest} and~\eqref{e:diffest2}.} 
The estimate~\eqref{e:diffest} immediately follows from the fact that $\eta(t) \leq 2CE_0$ for $t \geq 0$ upon taking $M_0 \geq 2C$. So, all that remains is to establish the pointwise bound~\eqref{e:diffest2}.  Recalling~\eqref{e:Fdef},~\eqref{e:normalform} and~\eqref{e:decompZc}, while using Propositions~\ref{p:IBP} and~\ref{p:IBP2}, we obtain a $t$- and $E_0$-independent constant $C_1 > 0$ such that
\begin{align*}
\left\|\U(t) - \sic(t)\right\|_{L^1} \leq \frac{C_1 E_0}{\sqrt{1+t} \log(2+t)^{\frac23}} 
\end{align*}
for $t \geq 0$. On the other hand,~\eqref{e:sicapprox},~\eqref{e:boundrrho} and~\eqref{e:restfin} yield a $t$- and $E_0$-independent constant $C_2 > 0$ such that
\begin{align*}
\left\|P_\pm \sic(t) - P_\pm(0) \re^{\tl_\pm(\cdot) t} A(t)\right\|_{L^1} \leq \frac{C_2}{(1+t)\sqrt{\log(2+t)}}
\end{align*}
for $t \geq 0$, where $\tl_\pm$ is given by~\eqref{e:def_truncation}. Using the latter two estimates and the reverse triangle inequality, we find a $t$- and $E_0$-independent constant $C_3 > 0$ such that
\begin{align}
 \label{e:final_est_pointwise}
\begin{split}
&\left|U(x,t)\right| - \left|P_+(0) \mathcal{F}^{-1}\left(\re^{\tl_+(\cdot) t}\right)(x) + P_-(0)\mathcal{F}^{-1}\left(\re^{\tl_-(\cdot) t}\right)(x)\right||A(t)|\\ 
&\quad \ \leq \left\|\U(t) - \sic(t)\right\|_{L^1} + \left\|P_+ \sic(t) - P_+(0) \re^{\tl_+(\cdot) t} A(t)\right\|_{L^1}\\ 
&\qquad \quad \ + \, \left\|P_- \sic(t) - P_-(0) \re^{\tl_-(\cdot) t} A(t)\right\|_{L^1} \\
&\quad \ \leq \frac{C_3}{\sqrt{1+t} \log(2+t)^{\frac23}}
\end{split}
\end{align}
for $x \in \R$ and $t \geq 0$. Finally, using the standard integral~\eqref{e:int_standard}, we compute
\begin{align*}
\mathcal{F}^{-1}\left(\re^{\tl_\pm(\cdot) t}\right)(x) = \frac{1}{2\pi} \int_\R \re^{\tl_\pm(k) t + \ri k x} \de k = \frac{\re^{\pm \ri (1+t) - \frac{x^2}{2(\alpha \mp \ri) (1+t)}}}{\sqrt{2\pi(\alpha \mp \ri)(1+t)}}
\end{align*}
for $x \in \R$ and $t \geq 0$. Combining the latter with~\eqref{e:boundrrho},~\eqref{e:restfin} and~\eqref{e:final_est_pointwise}, we arrive at the pointwise bound~\eqref{e:diffest2}, which finishes the proof. 
\end{proof}

Finally, we state the proof of Theorem~\ref{thm:mainresult}, which relies on an iteration argument in the variables $\U$, $\Us$ and $\Zc$ using $L^1$-$L^\infty$-estimates. 

\begin{proof}[Proof of Theorem~\ref{thm:mainresult2}] We close a nonlinear argument controlling the variables $\U(k,t)$, $\Zc(k,t)$ and $\Us(k,t) = (1-\chi(k))\U(k,t)$.
\bigskip\\
\noindent\textbf{Template function.} Using~\eqref{e:blowup},~\eqref{e:regularity2},~\eqref{e:regularity32} and~\eqref{e:regularity42} we observe that the template function $\eta \colon [0,T_{\max}) \to \R$ given by
\begin{align*}
\eta(t) &= \sup_{0 \leq s \leq t} \left[\big\|\U(s)\big\|_{L^\infty} + \big\|\Zc(s)\big\|_{L^\infty} + \sqrt{1+s}\left(\big\|\U(s)\big\|_{L^1} + \big\|\Zc(s)\big\|_{L^1} + \big\|\Us(s)\big\|_{L^\infty}\right)\right. \\
&\qquad \left. + \, (1+s)\big\|\Us(s)\big\|_{L^1}\right],
\end{align*}
is well-defined, continuous, monotonically increasing and, if $T_{\max} < \infty$, then we have
\begin{align} \label{e:blowupeta2}
\lim_{t \uparrow T_{\max}} \eta(t) = \infty.
\end{align}
Analogous to the proof of Theorem~\ref{thm:mainresult}, we establish the estimates~\eqref{e:boundUcUS},~\eqref{e:boundZC} and~\eqref{e:boundZC2} for $j,m = 0$. Combining these with~\eqref{e:cond_initial3} we find an $E_0$-independent constant $K_0 \geq 1$ such that
\begin{align} \label{e:eta0est2}
\eta(0) \leq K_0 E_0.
\end{align} 

\noindent\textbf{Key inequality.} Our aim is to establish a constant $C \geq K_0$ such that for each $t \in [0,T_{\max}) \cap [0,T_\eps]$ with $\eta(t) \leq 1$ the key inequality
\begin{align} \label{e:etaest2}
\eta(t) \leq C\left(E_0 + \eta(t)^2 \log(2+t)\right),
\end{align}
is satisfied. Then, taking
$$\eps < \frac{1}{4C^2}, \qquad M_0 = 2C,$$
it follows, by the continuity of $\eta$, that, provided $E_0 \in (0,\eps)$, it holds $\eta(t) \leq M_0E_0 = 2CE_0 \leq 1$ for all $t \in [0,T_{\max}) \cap \big[0,T_\varepsilon\big]$. Indeed, given $t \in [0,T_{\max}) \cap \big[0,T_\varepsilon\big]$ with $\eta(t) > 2CE_0$, there must, by continuity of $\eta$ and the fact that $\eta(0) < 2CE_0$ by~\eqref{e:eta0est2}, exist $s \in (0,t]$ with $\eta(s) = 2CE_0 \leq 1$. We arrive at the contradiction 
$$\eta(s) \leq C\left(E_0 + 4C^2E_0^2\log(2+s)\right) < 2CE_0,$$
by applying estimate~\eqref{e:etaest2} and the fact that we have $4C^2 E_0\log(2+s) < 1$ for $s \leq T_\eps = \re^{\eps/E_0} - 2$. We conclude that, if~\eqref{e:etaest2} holds, then $\eta(t) \leq 2CE_0$, for all $t \in [0,T_{\max}) \cap \big[0,T_\varepsilon\big]$, which implies by~\eqref{e:blowupeta2} that $T_{\max} > T_\varepsilon$. Consequently, we have $\eta(t) \leq M_0E_0$ for all $t \in \big[0,T_\varepsilon\big]$, which readily yields the desired estimates. So, all that remains is to establish the key inequality~\eqref{e:etaest2}.
\bigskip\\
\noindent\textbf{Bounds on $\Zc(t)$.} Let $t \in [0,T_{\max})$ with $\eta(t) \leq 1$ and $E_0 \in (0,\eps)$. We establish estimates on $\Zc(t)$ by estimating the terms on the right-hand side of its Duhamel formulation~\eqref{e:DuhZc}. First, Lemmas~\ref{l:Zcnonlinearest} and~\ref{l:residualnonlZc} and the estimate~\eqref{e:boundUcUS} with $j,m=0$ yield a $t$- and $E_0$-independent constant $C_1 > 0$ such that
\begin{align} \label{e:nonlest100}
\begin{split}
\left\|\mathcal{E}_4(\cdot,s)\right\|_{L^\infty}, \left\|Z_3^{\mathrm{res}}(\Zc(s),\Zc(s),\Zc(s))\right\|_{L^\infty} &\leq C_1\frac{\eta(s)^2}{1+s},
\end{split}
\end{align}
for $s \in [0,t]$, where we use $\eta(t) \leq 1$. On the one hand, Lemma~\ref{l:linear estimates} and estimate~\eqref{e:nonlest100} yield
\begin{align} \label{e:nlest15}
\begin{split}
&\left\|\int_0^t \re^{\hL(\cdot)(t-s)} \chi \left(Z_3^{\mathrm{res}}(\Zc(s),\Zc(s),\Zc(s)) + \mathcal{E}_4(\cdot,s)\right) \de s \right\|_{L^1}\\ 
&\quad \ \leq C_2 \int_0^t \frac{\eta(s)^2}{1+s} \int_\R \re^{-\frac12 \alpha k^2 (t-s)} \de k \de s \\
&\quad \ \leq C_3 \int_0^t \frac{\eta(s)^2}{\sqrt{t-s} (1+s)} \de s \leq C_4 \frac{\eta(t)^2\log(2+t)}{\sqrt{1+t}},
\end{split}
\end{align}
and
\begin{align} \label{e:nlest16}
\begin{split}
&\left\|\int_0^t \re^{\hL(\cdot)(t-s)} \chi \left(Z_3^{\mathrm{res}}(\Zc(s),\Zc(s),\Zc(s)) + \mathcal{E}_4(\cdot,s)\right) \de s \right\|_{L^\infty}\\ 
&\quad \ \leq C_2 \int_0^t \frac{\eta(s)^2}{1+s} \de s \leq C_3 \eta(t)^2\log(2+t)
\end{split}
\end{align}
for some $t$- and $E_0$-independent constants $C_2,C_3, C_4 > 0$. On the other hand, Lemma~\ref{l:linear estimates}, estimate~\eqref{e:boundZC2} with $j = 0$ and the fact that $\Zc(\cdot,0)$ is supported on $(-k_0,k_0)$, afford the linear estimates
\begin{align} \label{e:linst12}
\begin{split}
\left\|\re^{\hL(\cdot)t} \Zc(0)\right\|_{L^1} &\leq C_2 \int_\R \re^{-\frac12 \alpha k^2 t} \de k \big\|\U(0)\big\|_{L^\infty} \leq C_3 \frac{E_0}{\sqrt{t}},\\
\left\|\re^{\hL(\cdot)t} \Zc(0)\right\|_{L^p} &\leq C_2 \big\|\U(0)\big\|_{L^p} \leq C_3 E_0,
\end{split}
\end{align}
for $p = 1,\infty$ and some $t$- and $E_0$-independent constants $C_2,C_3 > 0$. Thus, applying~\eqref{e:nlest15},~\eqref{e:nlest16} and~\eqref{e:linst12} to the right-hand side of~\eqref{e:DuhZc} we obtain a $t$- and $E_0$-independent constant $C_* > 0$ such that
\begin{align} \label{e:finalZc}
\big\|\Zc(t)\big\|_{L^1} \leq C_*\frac{\eta(t)^2\log(2+t)}{\sqrt{1+t}}, \qquad \big\|\Zc(t)\big\|_{L^\infty} \leq C_*\eta(t)^2\log(2+t).
\end{align}
\smallskip\\
\noindent\textbf{Bounds on $\Us(t)$.} Let $t \in [0,T_{\max})$ with $\eta(t) \leq 1$ and $E_0 \in (0,\eps)$. We obtain bounds on $\Us(t)$, which satisfies the Duhamel formula~\eqref{e:Duhdamp}. First,  Lemma~\ref{l:noncriticalbounds} provides a $t$- and $E_0$-independent constant $C_1 > 0$ such that
\begin{align} \label{e:nonlest200}
\left\|\mathcal{F} \mathcal{N}\left(\mathcal{F}^{-1} \U(s)\right)\right\|_{L^p} \leq \frac{C_1 \eta(s)^2}{(1+s)^{\frac12 + \frac{1}{2p}}},
\end{align}
for $p = 1, \infty$ and $s \in [0,t]$, where we use $\eta(t) \leq 1$. Analogous to the proof of Theorem~\ref{thm:mainresult} we derive the estimates~\eqref{e:nlest5} and~\eqref{e:linst2} for $j = 0$ and $p = 1,\infty$ using Lemma~\ref{l:linear estimates} and estimate~\eqref{e:nonlest200}. Thus, applying~\eqref{e:nlest5} and~\eqref{e:linst2} to the right-hand side of~\eqref{e:Duhdamp} we arrive at 
\begin{align} \label{e:finalUs2}
\big\|\Us(t)\big\|_{L^1} \leq C_*\frac{E_0 + \eta(t)^2}{1+t}, \qquad \big\|\Us(t)\big\|_{L^\infty} \leq C_*\frac{E_0 + \eta(t)^2}{\sqrt{1+t}}.
\end{align}
for some $t$- and $E_0$-independent constant $C_* > 0$.
\bigskip\\
\noindent\textbf{Bounds on $\U(t)$.} Let $t \in [0,T_{\max})$ with $\eta(t) \leq 1$ and $E_0 \in (0,\eps)$. We apply estimate~\eqref{e:boundZC} with $j = 0$ and the bound~\eqref{e:finalZc} to find $t$- and $E_0$-independent constants $C_2,C_3,C_4 > 0$ such that  
\begin{align} \label{e:finalUC11}
\begin{split}
\big\|\Uc(t)\big\|_{L^p} &\leq C_2\left(\big\|\Zc(t)\big\|_{L^p} + \big\|\U(t)\big\|_{L^1} \big\|\U(t)\big\|_{L^p}\right)  \\
&\leq C_3 \left(\frac{\left(E_0 + \eta(t)^2\right) \log(2+t)}{(1+t)^{\frac{1}{2p}}} + \frac{\eta(t)^2}{(1+t)^{\frac12 + \frac1{2p}}}\right)\\
&\leq C_4 \frac{\left(E_0 + \eta(t)^2\right) \log(2+t)}{(1+t)^{\frac{1}{2p}}},
\end{split}
\end{align}
for $p = 1,\infty$ where we use $\eta(t) \leq 1$. Finally, estimates~\eqref{e:finalUs2} and~\eqref{e:finalUC11} imply
\begin{align} \label{e:finalU1}
\big\|\U(t)\big\|_{L^1} \leq C_* \frac{\left(E_0 + \eta(t)^2\right) \log(2+t)}{\sqrt{1+t}}, \quad \  \big\|\U(t)\big\|_{L^\infty} \leq C_*\left(E_0 + \eta(t)^2\right)\log(2+t),
\end{align}
for some $t$- and $E_0$-independent constant $C_*> 0$.
\bigskip\\
\noindent\textbf{Proof of key inequality.} The key inequality~\eqref{e:etaest2} follows readily by combining the estimates~\eqref{e:finalZc},~\eqref{e:finalUs2} and~\eqref{e:finalU1}, which concludes the proof.
\end{proof}

\bibliographystyle{abbrv}
\bibliography{KGbib}

\end{document}